UNIVERSIDAD AUTONOMA DE MADRID

**MÁSTER EN FORMACIÓN DE PROFESORADO DE EDUCACIÓN SECUNDARIA Y BACHILLERATO**

Título: La educación matemática en Rusia y los Círculos Matemáticos

Autor: Daniel Grilo Bartolomé

Director: Dr. Carlo Giovanni Madonna

TRABAJO DE FIN DE MÁSTER

Curso: 2014/2015

# La educación matemática en Rusia y los Círculos Matemáticos

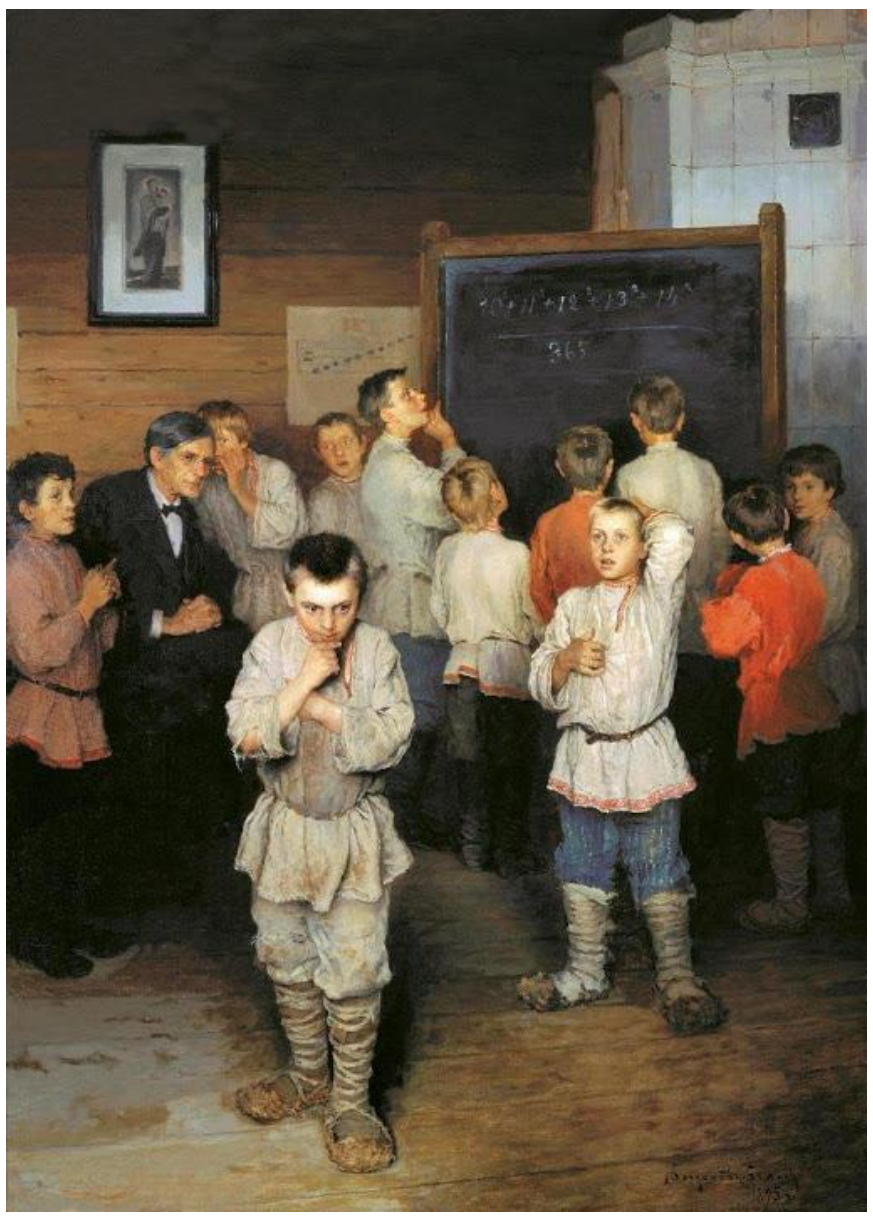

"***Aritmética mental en la escuela pública de Rachimsky***" de N. P. Bogdanov-Belsky (1895)

(El Dibujante 2.0: un ejercicio complicado, 2015)

El problema escrito en la pizarra, consiste en calcular mentalmente el valor de: $\frac{10^2+11^2+12^2+13^2+14^2}{365}$

Already before the war, Mekh-Mat operated **mathematical kruzhoks** for Moscow's elementary school students. In a kruzhok for seventh graders I tried to teach the differential and integral calculus through problems. The students learned the material after a minimal amount of training and learned it far better than high school students. (Landis, 1993, p. 57)

*Ya antes de la Guerra, existían los **Círculos Matemáticos** en la Facultad de Mecánica y Matemáticas dirigidos a los estudiantes de primaria y secundaria de Moscú. En un Círculo para alumnos de séptimo grado, enseñé cálculo diferencial e integral, a través de problemas matemáticos. Los alumnos aprendieron los contenidos después de un cierto entrenamiento y lo hicieron bastante mejor que los alumnos de secundaria en sus clases ordinarias.* (Traducción del autor de Landis, 1993, p. 57)

Once the remarkable Soviet geometer *Boris Nikolaevich Delone*, appearing before his pupils, expressed the idea that the only difference between a great scientific discovery and a good Olympiad problem is that solving the **Olympiad** problem requires 5 hours but finding a powerful scientific result requires an investment of 5000 hours. (Tikhomirov, 1993, pp. 106-107)

*En una ocasión, el memorable geómetra soviético Boris Nikolaevich Delone, apareciendo frente a sus alumnos, les expresó la idea de que la única diferencia entre un gran descubrimiento científico y la resolución de un dificultoso problema propuesto en las **Olimpiadas**, consiste en que éste último conlleva 5 horas mientras que un importante hallazgo científico requiere de 5000*. (Traducción del autor de Tikhomirov, 1993, pp. 106-107)

Try not to lecture. Even though introducing new theory and techniques is an integral part of **math circles**, your sessions should be as interactive as possible. Score yourself: 1 point per minute you talk; 5 points per minute a student talks; 10 points per minute you argue with a student; 50 points per minute the students argue among themselves. (Davis, 2014)

*Intenta no impartir clase. A pesar de que introducir teoría y métodos sea una parte fundamental de los **Círculo Matemáticos**, tus sesiones deberían ser lo más interactivas posibles. Puntúate a ti mismo: date 1 punto por cada minuto que hables, 5 puntos por cada minuto que hable un estudiante; 10 puntos por cada minuto que discutas con un estudiante; 50 puntos por cada minuto que los estudiantes discutan entre ellos.* (Traducción del autor de Davis, 2014)

# Índice

**Págs**

# Introducción

El presente trabajo pretende dar respuesta a la pregunta: ¿cómo debería producirse un proceso de enseñanza-aprendizaje de las matemáticas? La pregunta de carácter general, puede acotarse un poco más, si se piensa, en todos aquellos chicos y chicas, que se sienten atraídos hacia esta antigua ciencia, y que quieren contar con el mejor aprendizaje posible, es decir, aquel que es más significativo y constructivo, y que pone el acento en la forma de proceder que se da dentro de la investigación matemática. Pero no sólo para ellos, pues una vez concebido el proceso didáctico que se considere apropiado, éste puede ser trasladado al nivel en el que se quiera hacer uso de él y con los objetivos que se hayan propuesto conseguir.

La pregunta se podría haber enfrentado buscando una respuesta en la innovación educativa, sin embargo, se ha creído más oportuno centrarla en una cuestión histórica, es decir, cabe establecer una segunda pregunta: ¿dónde, históricamente, se ha atendido a la cuestión de la enseñanza de las matemáticas con especial interés y relevancia, y en qué ha consistido? La preocupación por el estudio de las matemáticas como las matemáticas mismas, es tan antigua al menos como la civilización humana. Pero en la búsqueda de la respuesta a la pregunta, tiene mayor sentido fijarse en la historia de la educación matemática contemporánea, más concretamente, en el s. XX. Dentro de este siglo, llama la atención la enorme tradición de la educación matemática de Rusia (perteneciente a la Unión Soviética entre 1922 y 1991), donde la atención que se le prestaba, trascendía a la enseñanza proporcionada en las escuelas ordinarias, y el interés hacia su estudio estaba comúnmente extendido dentro de su sociedad. Todo ello dio lugar, a la formación de muchos e importantes matemáticos durante este periodo, lo cual es un hecho ampliamente conocido, pero no lo es tanto en qué consistió la experiencia educativa matemática en la que se formaron.

Hablar de la educación matemática en Rusia remite a hablar, inevitablemente, de los Círculos Matemáticos, y viceversa, pues la metodología que les es característica está implícita o se utiliza en buena manera en el resto de experiencias innovadoras que se desarrollaron fundamentalmente de los años 30 a los años 80. Los dos primeros capítulos constituyen una indagación histórica en la educación matemática en Rusia. El primero de ellos, permite posicionar al lector dentro del contexto histórico donde se desarrollan los programas de educación matemática (extra-curriculares y extra-ordinarios) expuestos en el segundo capítulo y cuya lectura debe realizarse entendiéndolos no como objetos independientes, sino de manera conjunta, es decir como propuestas todas ellas interrelacionadas, bajo una manera característica de entender la educación matemática y donde el proceso didáctico de los Círculo Matemáticos, juega un papel fundamental por el influjo y utilización en el resto de estos programas. De esta forma, esta fórmula educativa, resultó ser predominante y fundamental en la educación matemática de aquellos jóvenes que quisieran, por diferentes motivos, que las matemáticas tuvieran un papel importante en su actividad intelectual, colaborando en buena medida en que la Unión Soviética se convirtiera en una potencia mundial en este aspecto.

Una vez comprendida la importancia y los fundamentos de los programas de la educación matemática en Rusia, la investigación se centra, en los dos siguientes capítulos, en aquello que es de interés por sus aplicaciones didácticas con el objetivo de trasladarlas a la enseñanza actual, que no es otro que el proceso didáctico inherente a los Círculos Matemáticos. La tradición de los Círculos Matemáticos se recoge hasta hoy en día, siendo los Estados Unidos, los nuevos impulsores de los mismos. También siguen desarrollándose, aunque no con la misma importancia, en la actual Rusia. En el tercer capítulo se analizan las experiencias actuales que se consideran más relevantes, desde el punto de vista de sus objetivos y organización. En el cuarto capítulo se proporciona una guía metodológica sobre cómo debe organizarse un Círculo Matemático, cual es su carácter y filosofía, y que criterios deben tenerse en cuenta a la hora de establecer un proceso didáctico dentro de él. Por último, y dado que la cuestión, aun partiendo de la respuesta de una pregunta de investigación histórica, se quiere conducirla a una cuestión práctica, se ha elaborado en forma de anexo, una propuesta consistente en cómo desarrollar sesiones al estilo de un Círculo Matemático. Estas sesiones están pensadas para proponerlas dentro de un aula de secundaria, y también, para hacerlo extra-curricularmente, organizando para ello un Círculo Matemático en una Escuela o Instituto de Secundaria. Las sesiones han sido denominadas, sesiones de investigación y sesiones de problemas, por ser precisamente estas dos tareas, los procesos de investigación matemática mediante la formulación y la respuesta a problemas de diferente temática y grado de dificultad, lo que más caracteriza a un Círculo Matemático.

## Estado de la cuestión

Al emprender un trabajo sobre la educación matemática en Rusia y los Círculos matemáticos la primera dificultad se encuentra en el desconocimiento de la lengua rusa que impide el acceso a la mayor parte de la bibliografía sobre el tema. Ello limita el campo de estudio restringiéndolo únicamente a la bibliografía escrita en lengua inglesa, no muy abundante aunque de gran interés que aborda desde el relato histórico hasta la cuestión didáctica de aplicación actual. Referencias fundamentales son las publicaciones “Russian Mathematics Education” (Karp y Vogeli, 2010a, 2010b), “Socialist Mathematics Education”, (Swetz, 1978) y toda la serie de 17 volúmenes que constituye la “MSRI Mathematical Circles Library” (han sido utilizados por el autor Fomin, Genkin, y Itenberg, 1996; Stankova y Rike, 2008, 2015; Vandervelde, Yashchenko, 2010; 2009; Sally y Sally, 2011; Zvonkin, 2011; Dorichenko, 2012; Rozhkovskaya, 2014)

La bibliografía en castellano es muy reducida, con un papel de puesta en conocimiento sobre alguno de los aspectos del tema al lector, pero sin llegar a desarrollarlo en ningún caso. Algunas publicaciones al respecto son “Poder y matemáticas en la antigua Unión Soviética” (Criado, 2011), “Andrey Kolmogorov: El último gran matemático universal” (Gordienko, 2011), “Los laboratorios matemáticos y el proyecto MateMatita” (Madonna y Gómez Esquinas, 2015), y la traducción al castellano del libro perteneciente a la serie “MSRI Mathematical Circles Library” titulado “Círculos Matemáticos” (Fomin, Genkin, y Itenberg, 2012).

Por lo tanto, el tema tal y como aquí se propone, no ha sido tratado por ninguna publicación en toda su extensión, es decir, desde los orígenes históricos hasta los objetivos y organización de las experiencias más recientes y la descripción del proceso didáctico, para la formulación de una propuesta de carácter práctico. De esta forma se ha conseguido desarrollar de forma extensiva el tema y acercarlo al lector en castellano. En este sentido, se puede calificar de original el tratamiento otorgado.

El trabajo constituye a la vez una invitación para seguir desarrollando el tema. Se consideraría de especial interés la realización de un análisis de la práctica didáctica desde la experiencia directa a través de la observación y la participación en los diferentes Círculos Matemáticos que existen en la actualidad, los cuales se encuentran fundamentalmente en Rusia y Estados Unidos. Desarrollando futuras investigaciones en este sentido, podría dar lugar a la obtención de nuevas y relevantes conclusiones que ayudasen en la mejora de los procesos didácticos matemáticos.

# 1. Resumen histórico de la educación matemática en Rusia durante sus etapas moderna y contemporánea.

A lo largo de este punto se exponen los principales hitos del desarrollo de la educación matemática rusa durante sus etapas moderna y contemporánea, desde el siglo XVIII hasta llegar al Currículo de 1981, con el fin de establecer un punto de partida para entender el contexto y los planteamientos de la educación matemática del sistema educativo en la Unión Soviética entre 1930 y 1985, objeto de la investigación del presente trabajo.

La educación matemática en Rusia ha contado en su desarrollo moderno con multitud de cambios y avances, que perseguían la constante preocupación por mejorarla. La importancia de la enseñanza de las matemáticas dentro de la educación rusa ha sido prácticamente desde sus inicios de primer orden hasta llegar a la Unión Soviética, donde se le ha prestado una importancia aún mayor dotándola de una red educativa de carácter extra-ordinario y extra-curricular que acompañaba y ampliaba la educación ordinaria, aunque fundamentalmente también se ocupaba de cubrir otros aspectos como desarrollar las habilidades características de un matemático que se dedica a la investigación. Es por ello, que la implicación de la comunidad matemática Universitaria desempeñaría un papel fundamental en su formalización.

## 1.1. La educación matemática en Rusia durante el s. XVIII.

La historia moderna de la educación matemática rusa comienza bajo el absolutismo de Pedro I. Este zar, quien a su vez fue instruido en matemáticas, supo ver la importancia de esta disciplina en la instrucción técnica y militar. Decidido, entre otras cosas, a mejorar las fuerzas navales rusas, impulsó la creación en Moscú de las primeras **escuelas de matemáticas y de navegación** en 1701**,** que constituirían el primer intento de creación de una educación civil, pública, nacional y de carácter gratuita. Estas escuelas formaron no sólo a la población militar sino también a la civil en diferentes tipos de oficios, como ingenieros, arquitectos y futuros profesores, entre otros. La escuela se trasladó en 1715 a San Petersburgo, y pasó a integrarse dentro de la **Academia del cuerpo de las fuerzas navales**. Para entonces había formado a unos 500 alumnos. El programa educativo fue impulsado en un primer momento por Fargwarson, profesor de la Universidad de Aberdeen, e incluía aritmética, álgebra, geometría plana y esférica, y trigonometría. Realizó también labores como profesor y autor de libros de texto. Debido a los consecuentes problemas con el inglés por parte de los estudiantes, Magnitsky, personaje ilustrado de la época, fue invitado como profesor, y acabó imponiéndose como líder de la escuela. Será él el autor de varios libros de texto de matemáticas, como el titulado ***"Aritmética"*** y publicado en 1703, libro que funcionaba como una enciclopedia matemática de la época, ya que Magnitsky había compilado diferentes conocimientos matemáticos partiendo de otras publicaciones en distintos idiomas que él dominaba (Polyakova, 2010, p. 8).

En su intento de educar a las clases bajas, Pedro I, fundó también las llamadas **escuelas de aritmética** en 1714, las cuales se distribuyeron por las diferentes provincias rusas y cuyos profesores habían sido formados en las escuelas de navegación y matemáticas de Moscú. En 1727, estas escuelas, contaban con unos 2000 estudiantes procedentes de diferentes estratos sociales (Karp, 2014, p. 305). La asistencia a ellas era obligatoria y tenían un régimen muy estricto, por lo que las familias solían ser reticentes a enviar a sus hijos, prefiriendo que estuvieran vinculados a los negocios familiares (Polyakova, 2010, p. 10). Con el tiempo algunos estratos (nobleza y clero) se vieron libres del carácter de obligatoriedad, pero decidieron enviar a sus hijos a otras escuelas militares y religiosas. Como consecuencia de ello, el número de estas escuelas se redujo drásticamente, de 42 a 28, asistiendo únicamente a ellas miembros de las clases bajas. En 1744, fueron absorbidas por las ya existentes **escuelas militares garrison** fundadas también por Pedro I en 1716. En las escuelas de aritmética se estudiaba aritmética y geometría. El profesor se hacía cargo de 20-30 alumnos que a su vez, tenían que estudiar también otras asignaturas. No había metodología alguna ni libros de texto que sirvieran de guía. La educación se basaba principalmente en la memorización de reglas y la resolución de problemas (Polyakova, 2010, p. 10).

En estos inicios del siglo XVIII también aparecieron en Rusia las primeras escuelas de educación secundaria con una carácter generalista, los denominados **Gymnasium** (Gimnasio) fundados por extranjeros: El **Gimnasio Ernst Glück** en Moscú (1703-1715) y el **Gimnasio de la Academia de San Petersburgo**, academia de ciencias establecida en 1724, y que abrió en 1726, concebido para preparar a los estudiantes para su etapa universitaria. Se invitó a participar en ellos a diferentes científicos de renombre de la época, entre ellos el alemán Leonhard Euler, que contribuiría de forma importante en el desarrollo de la educación matemática rusa, no sólo en contenido sino también en metodología (Polyakova, 2010, p. 12). Estos Gimnasios se apartaban de la concepción de las otras escuelas, las cuales constituyeron un primer intento, aunque fallido, de educar a toda la población civil incluidas las clases bajas. Contrariamente, en ellos la educación contaba con un carácter elitista y de proyección internacional. El Gimnasio de la **Academia de San Petersburgo** estuvo abierto durante todo el s. XVIII cerrando finalmente a principios del s. XIX.

Leonhard Euler introdujo, a lo largo de todo el siglo XVIII, importantes novedades metodológicas importadas de Europa, donde las matemáticas se formulaban sistemáticamente y en base a demostraciones, y que fueron adaptadas a la educación matemática rusa, con la creación de textos en ruso. Es el caso del importante texto del propio Euler titulado ***"Manual de aritmética para el Gimnasio de la Academia Imperial de Ciencias"*** escrito en alemán y traducido al ruso en 1740, que constituyó el primer manual accesible con contenido matemático de alto nivel académico (Polyakova, 2010, p. 13). También publicó después, en 1770, su ***"Introducción al álgebra"***, tanto para su uso en el Gimnasio como en la Universidad. Los recientes logros matemáticos fueron también incluidos en los libros de texto, como la trigonometría moderna de Euler, y sus trabajos en cálculo diferencial e integral. Otro cambio metodológico importante introducido por Euler consistió en el desglose de las matemáticas en diferentes disciplinas constituidas en asignaturas para su estudio, unido a la introducción de ciertos principios didácticos: sistematización, fundamentación científica y accesibilidad de la exposición matemática (Polyakova, 2010, p. 14).

En los inicios del s. XVIII, se incrementó la demanda de una educación avanzada para la población, en parte porque había cada vez más gente con una educación matemática elemental. Además, la idea de que los estudios de matemáticas eran especialmente necesarios empezó a alcanzar una posición relevante, y matemáticos de gran de calado se empezaron a hacer cargo de la formación matemática. Es el caso del **sistema educativo profesional**, el cual estaba formado por las **escuelas militares** (cuerpos de cadetes navales y de la armada), **escuelas técnicas militares**, y otras, donde en todas ellas, la enseñanza matemática era considerada de primer orden. La enseñanza matemática en los cuerpos de cadetes navales era de considerable calidad, por contar precisamente con importantes pedagogos y matemáticos provenientes de la escuela de Euler, como Nicolas Fuss, asistente de Euler y autor de varios libros de texto matemáticos. En 1760 la Academia ofrecía también un curso de matemáticas superiores complementando al curso de matemáticas elementales ya ofrecido. A este curso, a finales del s. XVIII, se le añadirían estudios de geometría analítica y de análisis, y sería impartido por el propio Nicolas Fuss (quien también fue profesor del cuerpo de cadetes de la armada). Por otra parte, en las escuelas técnicas militares, bajo la dirección de M. I. Mordvinov, se incluían la enseñanza de mecánica, aritmética, introducción al álgebra y geometría. Uno de sus mejores profesores fue Y. P. Kozel'skij, quien expresó mediante ensayos filosóficos la importancia en la conexión entre la teoría y la práctica. Otro de sus profesores, N. V. Vereschagin, discípulo de Kozel'skij, sería uno de los primeros en enseñar geometría analítica en Rusia, y tuvo una gran influencia en el desarrollo de la educación matemática durante toda esta segunda mitad del s. XVIII. De esta forma, se formalizó una instrucción matemática avanzada cuyo peso recaería también en la educación universitaria. Sin ir más lejos, M. V. Lomonosov (1711-1765), considerado el padre de la ciencia rusa, estudió inicialmente en el **Gimnasio de la Academia de San Petersburgo** y posteriormente en la Universidad.

En 1755, se fundó la **Universidad de Moscú** por iniciativa de M. V. Lomonosov, donde en un principio sólo se estudiaban matemáticas como asignaturas auxiliares para otros estudios. Habrá que esperar hasta 1759 para que D. S. Anichkov presida una facultad de matemáticas ofreciendo un curso en matemáticas puras en un ciclo de dos años. En el primer año se enseñaba aritmética y geometría y, durante el segundo año, geometría y trigonometría. Más adelante el curso se amplió con una duración total de tres años y fue impartido por V. K. Arshenevskij incluyendo los estudios de álgebra en el tercer año. También se llegó a impartir un curso en un ciclo de tres años de matemáticas aplicadas. En general el nivel será inferior que en las academias y las escuelas profesionales de la época.

Lo más interesante, aconteció a finales de los años 60, con la aparición de cursos para la preparación pedagógica de futuros profesores de Gimnasios, constituyéndose finalmente en un seminario conocido como el **seminario de los profesores**. En 1755, se fundaron dos Gimnasios en Moscú y poco después otro en Kazan. Posteriormente se fundaron escuelas internados y escuelas para los hijos de la nobleza, todas ellas apoyadas por la Universidad de Moscú que les proporcionaría profesores, libros de textos e incluso dinero en determinadas ocasiones (Polyakova, 2010, p. 19). En general, el nivel matemático de los Gimnasios era bastante bajo, y tan sólo la Escuela Internado de Moscú para hijos de la nobleza contaba con un nivel aceptable, al orientar sus estudios hacia un carácter práctico en vez del de humanidades, característico de los Gimnasios (Polyakova, 2010, pp. 19-20).

Llegados a este punto se puede señalar que la forma en la que se organizaba la educación en Rusia distaba bastante de ser un sistema organizado. I. I. Shuvalov, interceptor de M. V. Lomonosov frente a la emperatriz Isabel Petrovna (reinado: 1741-1762) en la fundación de la Universidad de Moscú, ya había pensado en organizar la educación en un sistema conformado por escuelas de educación básica y Gimnasios, pero finalmente sus planes no se llevaron a cabo. Será la emperatriz Caterina II (reinado: 1762-1796), seguidora de los enciclopedistas franceses e interesada en temas de educación, quien afronte el problema con la creación de una estructura a nivel nacional de escuelas públicas de educación básica, para lo que formará una comisión de investigación liderada por el científico y matemático Franz Aepinus a favor del modelo austriaco divido en tres niveles (“trivial”, “real” y “normal”). La reforma será llevada a cabo por la **Comisión de Escuelas Nacionales** compuesta por los rusos P. I. Zavadovsky, F. Aepinus, P. I. Pastukhov y el serbio F. Yankovich de Mirievo como consejero experto, enviado por el emperador austriaco José II en respuesta a la solicitud de Caterina II (Madariaga (de), 1979, p. 383). El reciente modelo ruso se inspiró en el modelo educacional del imperio Austriaco, imperio que se asemejaba al Ruso en su estructura económica y social y cuyo modelo educativo era moderno, utilitario y estaba abierto a toda la población, a la vez que laico y controlado por el estado, razón por la cual debieron de decidir imitarlo (Hans, 2012, p. 14). La comisión se encargaría de diseñar la estructura nacional de escuelas, la formación de los profesores y de proporcionar los libros de texto al alumnado, los cuales son atribuidos a M. E. Golovin (asistente de Euler) (Madariaga (de), 1979, p. 384). F. Yankovich publicará un ***“Manual para profesores de Escuelas Públicas”*** en 1783 constituyendo el primer tratado metodológico en la historia de la educación de Rusia basado en la obra de Abbot Felbiger, y divido en cuatro partes: métodos de enseñanza, las asignaturas a enseñar, el carácter del profesor y el funcionamiento de la escuela (Madariaga (de), 1979, pp. 385-386).

La estructura de las escuelas públicas quedó formalizada finalmente en 1786 a través del ***“Estatuto Ruso de Educación Nacional”***, obligando a cerrar a las iniciativas extranjeras y privadas y absorbiendo al resto de escuelas, que se dividieron en escuelas de primaria con dos cursos ubicadas en las principales ciudades de los distritos y, escuelas de secundaria con cuatro cursos ubicadas en las capitales, con una concepción laica, científica y utilitaria, y a las que podía acudir toda la población sin distinción de clase ni género, con el fin de proporcionarles acceso a la Universidad. La educación matemática de las escuelas de primaria, incluía cálculo escrito y mental en el primer curso ampliándose con “la regla de tres” en el segundo curso. En las escuelas de secundaria, se revisarían los conocimientos adquiridos en primaria durante los dos primeros cursos; en el tercer curso se incluirían fracciones, números decimales, y diferentes tipos de ejercicios incluyendo “la regla de tres”; en el cuarto curso se estudiaría geometría plana. Hacia 1800, habría en todas las escuelas de Rusia en torno a 20000 estudiantes, cifra que puede parecer escasa, pero hay que tener en cuenta que no se habían creado escuelas en los pueblos y que no era fácil conseguir la asistencia de todos los chicos y chicas en edad escolar (Karp, 2014, p. 308). No obstante y a pesar de no conseguir llegar a toda la población, el s. XVIII termina para Rusia, con la formalización y puesta en marcha de la primera estructura educativa de carácter público.

### 1.2. La educación matemática en Rusia durante el s. XIX.

En 1802 se fundó el **Ministerio de Educación** y en 1804 se aprobó un primer Estatuto que establecía la organización de las escuelas, las cuales estarían supeditadas al control de la Universidad, dentro de cada una de las seis regiones académicas en que se había dividido Rusia. La estructura de escuelas quedó dividida en: **escuelas provinciales** ("parish") con un curso, **escuelas de distrito** ("uyezd") con dos cursos y **Gimnasios** con cuatro cursos, aunque durante algún tiempo coexistiría la vieja estructura con la nueva. De esta forma, la Universidad pasó a establecer cierta dirección sobre la educación primaria y secundaria y pudo implementar en ellas su método científico. La estructura educativa era ampliamente aceptada por la población a la que además tenían acceso gratuito sin restricción de clase hasta que en a principios de los años 20 la gente se vio obligada a pagar para acceder a la educación (Polyakova, 2010, p. 23).

Durante este periodo el debate de la educación se extendió a la elección entre un sistema de enfoque "clásico", es decir basado en el estudio de lenguas clásicas, y un sistema de enfoque "práctico o real", basado en el estudio de matemáticas y ciencias naturales. El debate continuaría a lo largo del s. XIX y el posterior s. XX, y sería fruto de continuas reformas y contrarreformas. La educación en los Gimnasios tenía un enfoque práctico donde las matemáticas asumían un papel prioritario con una dedicación de 6 horas de clases a la semana durante los tres primeros cursos, además, durante el tercer y cuarto curso se impartía también estadística, con 2 horas y 4 horas de duración semanal, respectivamente. Finalmente, el gobierno, acabó adoptando nuevamente un sistema de educación gratuita y sin restricciones de clase (Polyakova, 2010, p. 23). El contenido de la programación fue adoptado por dictámenes de la institución académica, aunque se mantuvo el uso generalizado de los mismos libros de texto, los cuales se mejoraron. Para ello se creó un comité, con N. I. Fuss formando parte de él, que publicó en 1805 los dos primeros volúmenes de ***"Curso de Matemáticas"*** escritos por T. F. Osipovsky y utilizados en los Gimnasios por la gran calidad de contenido y la sencillez de su exposición, a pesar de que el nivel de esta publicación estaba dirigido a la educación universitaria (Polyakova, 2010, p. 24). Es por esto último que el comité terminó reutilizando el material existente, especialmente dirigido al nivel de los cursos de los Gimnasios, para publicar ***"Fundamentos de las Matemáticas Puras",*** que sería utilizado como libro de texto entre 1814 y 1828. El desarrollo metodológico también continuó vigente con las primeras publicaciones metodológicas matemáticas por parte de *S.* E. Gur'yev.

En 1828 se aprobó un nuevo estatuto que introdujo los 3 años de educación básica inicial también en los Gimnasios, pasando a tener un total de 7 cursos. Si se cursaba en las escuelas provinciales y de distrito, no era posible el acceso a los Gimnasios. Después se introdujeron de nuevo limitaciones de clase, lo que permitió la entrada, tan sólo a los hijos de los nobles y de los mercaderes del primer gremio (Polyakova, 2010, p. 25). A su vez, los Gimnasios quedaron divididos en dos tipos, los que incluían el estudio del latín y el griego, y los que incluían el estudio tan solo del latín, lo que se puede hacer corresponder con el sistema de enfoque "clásico" y el sistema de enfoque "práctico" respectivamente (Polyakova, 2010, p. 25). Este estatuto, incluía también por primera vez un **plan de estudio**. En los Gimnasios con un sistema práctico se estudiaba aritmética en los cursos primero y segundo, introducción al álgebra incluyendo ecuaciones de segundo grado en el tercer curso, una segunda y última parte de álgebra y una primera parte de geometría en el cuarto curso, una segunda y última parte de geometría en el quinto curso, geometría descriptiva y aplicaciones del álgebra a la geometría en el sexto curso, y una última parte de aplicaciones del álgebra a la geometría incluyendo secciones cónicas en el séptimo curso. En los Gimnasios con un sistema clásico el estudio de las matemáticas se redujo a un total de 15 horas a la semana, contabilizando los 7, frente a las 46 que se dedicaban al estudio del latín y el griego.

El equilibro entre los dos sistemas, se mantuvo hasta 1834, cuando con la llegada de S. S. Uvarov al puesto de primer ministro, se impuso como dominante el sistema clásico, lo que hizo perder su posición privilegiada a las matemáticas, cuyo programa se redujo a 20 horas a la semana, distribuidas en los 4 últimos cursos. En 1846, F. I. Busse amplió en contenido el curso de aritmética y de álgebra, e introdujo la trigonometría en el programa con el objetivo de desarrollar la parte aplicada de las matemáticas. En 1848, a raíz de las revoluciones sucedidas en Europa, el ministro de educación trató de eliminar la "libertad de pensamiento" en la última etapa de los Gimnasios y en la Universidad (Polyakova, 2010, p. 27), por lo que terminó imponiendo el sistema práctico frente al clásico e incrementado a 30 la cantidad de horas a la semana de matemáticas, distribuidas a lo largo de esos 4 últimos cursos, que se vio reducida a 22,5 horas semanales en 1852. Elaboró un nuevo programa, aunque no muy diferente al de 1846, que imponía nuevamente un enfoque aplicado de la enseñanza de las matemáticas. Aparecieron nuevos libros de texto escritos por F. I. Busse , P. S. Gur'yev y V. Y. Bunyakovsky principalmente, que incluían libros de ejercicios, que junto a los libros del profesor constituyeron un primer conjunto metodológico en la historia de la enseñanza matemática (Polyakova, 2010, p. 28). D. M. Perevoshchikov y N. L. Lobachevsky, rectores de las Universidades de Moscú y Kazan, respectivamente, tuvieron una importante influencia y publicaron también libros de textos usados en los Gimnasios.

En 1858, P. L. Chebysev propuso un nuevo plan para regular la enseñanza matemática, consistente en la partición de los 7 cursos de los Gimnasios en una primera etapa de carácter general para los primeros 5 cursos y una última etapa de especialización para los 2 últimos cursos. Finalmente, el *Comité Académico*, rechazó esta bifurcación. En 1860, se propuso una nueva regulación decantándose por un sistema práctico, con 27,5 horas a la semana, distribuidas para todos los cursos, dedicados a la enseñanza de matemáticas y 24 horas a la semana para la enseñanza de física y ciencias naturales. En 1862 y 1864, se volvió a modificar, conllevando un decremento de la carga horaria en matemáticas. En 1864, tuvo lugar una conferencia en Odessa sobre educación, para discutir los problemas de la educación y, en concreto, sobre la sobrecarga de contenido en el programa, lo que condujo a adoptar una reducción del contenido del programa y la propuesta de nuevos principios metodológicos para la enseñanza de las matemáticas. Por último todo ello desembocó en nuevas regulaciones, apareciendo tres posibles elecciones de enseñanza: la clásica, la clásica con sólo una lengua, y la práctica; que garantizaban la igualdad de derechos a todos los graduados (Polyakova, 2010, p. 33). En 1865, P. L. Chebyshev realizó algunas modificaciones, lo que supuso que los profesores pudieran elegir sobre su programación con la autorización de un comité facultativo, rompiendo así la línea de uniformidad que hasta ahora se había venido dando. La carga lectiva semanal, resultó ser de 22 horas para las enseñanzas académicas y de 25 horas para la enseñanza práctica.

En 1871 se aprobaron nuevas regulaciones, todos los Gimnasios pasaron a tener únicamente una enseñanza del sistema clásico, y la enseñanza práctica se impartiría en las llamadas escuelas "prácticas" (real). La enseñanza se amplió a 8 cursos lectivos. En los Gimnasios las asignaturas de ciencias se llevaban 37 horas semanales (repartidas en los 8 cursos) frente a las 49 del latín o las 36 del griego. En las escuelas "prácticas" los estudiantes podían elegir, finalizado el sexto curso, entre tres opciones: comercial, técnica-mecánica o química. En estas escuelas no se estudiaban lenguas clásicas, las matemáticas y el dibujo técnico eran asignaturas importantes, y una vez graduados no se permitía el acceso a la Universidad. Dentro de estas regulaciones se aprobaron Currículos para todas las asignaturas por primera vez en la historia educativa rusa. El Currículo de matemáticas se aprobó en 1872, a cargo de P. L. Chebyshev. En 1888 se introdujeron nuevas regulaciones en las escuelas "prácticas" y en 1891 para los Gimnasios, pero apenas afectaron al programa de matemáticas. En 1890, por primera vez en la historia de la educación rusa, se incluyeron en el Currículo referencias metodológicas, en las que se mencionaba, en el Currículo de matemáticas, la importancia de una clara exposición teórica, la cual debía ser expuesta mediante ejemplos prácticos que facilitaran su comprensión, a la vez que el desarrollo de capacidades de cálculo, por lo que la visión de las matemáticas aplicadas se vio relegada una vez más. El programa de 1890 se mantuvo hasta 1917, y como resultado, el estudio clásico de los Gimnasios se vino a menos ya que el griego dejó de ser obligatorio, el número de horas de clases de latín se redujo, el examen escrito de latín se eliminó y el número total de horas semanales dedicadas a las matemáticas se estableció en 30 (distribuidas en todos los cursos).

Al final del siglo, el Ministerio convocó a la opinión pública para evaluar los problemas en la educación secundaria, lo que conllevó el desarrollo de una nueva estructura para las escuelas: Gimnasios con dos lenguas clásicas, Gimnasios con una lengua clásica, escuelas "prácticas" y un cuarto tipo de escuelas. El objetivo de la enseñanza de las matemáticas se concretó en: "las matemáticas serán adoptadas como una ciencia en todo su derecho y como un método científico de indagación" (Traducción del autor del texto citado en Polyakova, 2010, p. 35). Este último sistema de 1890 suele ser conocido bajo el nombre de **sistema clásico internacional de enseñanza matemática**, caracterizado por integrar de forma actualizada las matemáticas del s. XVII y anteriores, la división en diferentes áreas y cursos para su estudio, y la inclusión de principios metodológicos. Sistema que tiene sus ventajas pero también sus deficiencias (Polyakova, 2010, p. 36). En 1897, 58.092 estudiantes asistieron a los Gimnasios y 24.279 a las escuelas "prácticas", aunque ahora la población del país era de 125 millones (Karp, 2014, p. 310), por lo que una vez terminado el s. XIX estas cifras continuaban siendo significativamente pequeñas. Como resultado de este sistema, tan sólo el 40-50% de la población matriculada se graduaba, lo que llevo a que se generase la opinión extrema de eliminar completamente la enseñanza de las matemáticas en la educación básica (Polyakova, 2010, p. 37).

## 1.3. La educación matemática en Rusia (URSS) en el s. XX hasta 1985.

### 1.3.1. La educación matemática en Rusia (URSS) durante los inicios del s. XX hasta 1931: la reforma en contra del sistema clásico.

En los inicios del s. XX, las ideas reformistas no se hacen esperar. A finales del s. XIX ya se habían inaugurado con la publicación en diferentes revistas, de trabajos donde se criticaba el sistema educativo y se pedía una renovación de los programas que incluyeran nuevos cursos de geometría, y la inclusión del análisis matemático y la geometría analítica, entre otras demandas, que eran ignoradas por el Currículo de matemáticas de 1890. Las demandas eran de calado internacional. En 1897, en el Congreso Matemático Internacional celebrado en Zurich, Felix Kein, propuso que la reforma de la educación matemática era necesaria y expuso cuales debían ser sus principios en el **Congreso de Merano** en 1905, lo que incluía la integración de las matemáticas puras y aplicadas en el Currículo. Esta nueva visión reformista, alcanzó mayor peso en Rusia con la publicación en 1911 de la revista ***"Educación Matemática"*** que se convirtió en centro de discusiones, y sobre todo con los **Congresos de profesores de matemáticas** celebrados, el primero en San Petersburgo en 1911-1912, con la participación de 1217 profesores y 71 presentaciones; y el segundo en Moscú en 1912-13, con la participación de 1200 profesores y 22 presentaciones (Polyakova, 2010, pp. 38-39). Se discutieron en ellos, temas como la necesidad de incluir el estudio de funciones en el Currículo, los problemas en referencia al estudio del cálculo elemental, la importancia de la geometría visual y el método de laboratorio, la utilización de la historia en la exposición, la utilización de juegos para desarrollar las capacidades intelectuales, etc. Gracias a los congresos, se establecieron las bases sobre las que debía construirse la reforma que modernizase la educación matemática rusa. Las publicaciones académicas sobre metodología matemática se siguieron sucediendo durante los años 1912-1915 con diversos contenidos provenientes del mismo, lo que sin duda sentó las bases para la reforma llevada a cabo tras la revolución de 1917.

El descontento en los inicios del s. XX con el sistema educativo era generalizado, y como el ministro de educación, Tolstoy (1858-1916), de estos primeros años, menciona, todo el mundo parece estar de acuerdo en su mala calidad, pero todos dan diferentes motivos para ello (Karp, 2010a, p. 45). Los profesores de Universidad se quejaban de la falta de preparación de los estudiantes graduados en los gimnasios. El panorama se presenta por tanto ante la necesidad de un reforma profunda, y ésta se terminará realizando como un cambio absoluto contra el sistema hasta ahora clásico de enseñanza, coincidente con la coyuntura de las dos revoluciones de febrero y octubre de 1917 y la consecuente llegada al poder del **Partido Comunista**, encabezado por Lenin, que implantará a partir de ahora un modelo de educación que trata de dar respuesta a la ruptura con el antiguo régimen monárquico y la creación de la nueva sociedad comunista.

En los años previos a la revolución, los partidarios de la “educación libre” encabezada por L. Tostoy, y también, los partidarios de la “educación realista” se habían interesado por la escuela progresiva norteamericana de John Dewey y sus “escuelas nuevas” basadas en pedagogías activas ("learning by doing") (Mchitarjan, 2009, p. 166). La reforma supuso, una vez aprobado el decreto de 1918, la unificación de todos los tipos de escuelas y Gimnasios en la **escuelas de trabajo**, en las cuales se establecían dos etapas educativas: una primera etapa de 5 cursos (edades de 8 a 13) y una segunda etapa de 4 cursos (edades de 13 a 17). Los exámenes y las tareas para casa se prohibieron, y en general se rechazó "la vieja disciplina, que destruía la vida en la escuela y el libre desarrollo de la personalidad de los chicos" (Traducción del autor del texto citado en Karp, 2012). La adopción de esta nueva y rompedora educación, venía en parte influida por la escuela progresiva norteamericana, cuya forma de entender la educación se ajustaba a los ideales comunistas y fueron tomadas como referencia por los pedagogos rusos, al menos en la primera mitad de los años 20, pues después, estos mismos pedagogos la rechazaron argumentando que no se ponía al servicio de la “lucha de clases”, al priorizar la educación en sí misma en el individuo para dotarlo de autonomía, sin ninguna intencionalidad política o que sirviera a intereses de Estado (Mchitarjan, 2009). Para sus escuelas, se terminó optando por tanto por el estudio-trabajo en torno a temas (“complexes”), en vez de asignaturas, el método de laboratorio, el trabajo cooperativo, y una conexión permanente con la experiencia en todos los procesos de aprendizaje tal y como proponía Dewey. Integrantes de la escuela progresiva norteamericana, como el propio Dewey, Kilpatrick, Counts, y otros, visitarían la Unión Soviética y expresarían su admiración por el diseño de su sistema educativo.

Trabajar en torno a temas ("complex-based approach") significaba que una vez seleccionado un tema ("El servicio postal", "La casa", "El trabajo agrícola y su organización", etc.), éste debía desarrollarse a través del uso de las diferentes áreas de conocimiento. A. V. Lankov autor de ***"Mathematics and the Complex"*** recomendaba el uso de 4 temas durante el cuarto curso para cada uno de los 4 periodos y 3 temas de desarrollo más corto. Narazenko (1926), profesor de provincia, narra su experiencia de cómo podían desarrollarse las matemáticas a través de un tema. De este modo, se hicieron muchos esfuerzos para desarrollar el Currículo de matemáticas (1918), especialmente O. A. Volger, encargado de su transformación, y partidario de no enseñar matemáticas, ni ninguna otra materia. Según se cita en Karp (2012) "utilizarlas directamente aplicándolas bajo la idea de una interdependencia funcional, acostumbró, por ejemplo, a los alumnos al empleo de ecuaciones desde el primer año, de la misma forma que las propias matemáticas se habían desarrollado históricamente" (p. 556) [traducción del autor]. Sin embargo no tardaron en aparecer opiniones en contra, como la de I. I. Chistiakov, perteneciente al círculo de matemáticos de Moscú y participante en los Congresos de Profesores de Matemáticas, que defendía el valor del razonamiento matemático y la imposibilidad de ser enseñado de esta manera. A. V. Lankov también afirmó que la enseñanza de las matemáticas estaba desconectada de los temas, y varios matemáticos en general argumentaron sobre la dificultad de su conexión, como K. Rashevsky que exponía que salvo la aritmética, las matemáticas no tenían prácticamente conexión con la vida cotidiana, que la enseñanza en torno a temas iba en contra del método de laboratorio (que no se podía aplicar), y que era imposible desarrollar habilidades matemáticas a través del uso de temas; por lo que en general proponía que se enseñara de forma independiente y que sólo se utilizaran temas en determinados casos específicos donde pudiera resultar interesante, posibilidad albergada únicamente en la etapa de primaria y bajo ningún concepto en la de secundaria (Karp, 2012, p. 556).

A pesar de ello, el programa de 1920 recogía que: "el significado de las matemáticas residía en el hecho de que proporcionaba métodos que eran insustituibles en su aplicación a la realidad, y que el valor de las matemáticas como ciencia abstracta, como una forma de gimnasia mental para los jóvenes, carecía de valor alguno" (Traducción del autor del texto citado en Karp, 2010a, p. 47). El Currículo de 1925 señalaba, entre otras aportaciones que: "las matemáticas en sí mismas no tienen ningún valor educativo en las escuelas; las matemáticas tan sólo son importantes en el grado de ayuda que proporcionan para resolver problemas prácticos [...] en que puede ayudar al alumno en la lucha diaria y en la construcción de la vida" y "la adquisición de procesos de lectura, escritura, habilidades de cálculo mental y mediciones tienen que estar unidas al estudio de las prácticas reales; la aritmética y la lengua rusa no deben existir en las escuelas como asignaturas separadas" (Traducción del autor de Karp, 2012, p. 556). Varios autores, como Korolev, Korneichik, Ravkin, y Nikitin, comentan que en la práctica, especialmente en los cursos superiores, no se siguieron las recomendaciones, y que la enseñanza de las matemáticas se siguió realizando de forma separada. El Currículo de 1925, permitía en la segunda etapa clases separadas de matemáticas, pero motivadas por el desarrollo de la propia extensión del tema utilizado, y que podían ser programadas como el profesor estimara conveniente.

El **Plan Dalton**, que utilizaba también un método de laboratorio, y que había comenzado a aplicarse por la educadora progresiva norteamericana Helen Parkhurst en 1916, se utilizó también en las escuelas para desarrollar la metodología. I. B. Charnetsky (profesor, 1929), explicaba como las clases fueron sustituidas por laboratorios y como el Currículo era organizado mediante tareas con las que cada grupo o alumno procedía a su propio ritmo. Exponía que las consecuencias de este método eran: "una mayor autonomía, mayor capacidad para investigar, y que desarrollaba en el alumno, las tan deseadas, combinación y conversión entre las diferentes áreas matemáticas" (Traducción del autor de Karp, 2012, p. 557). El congreso poco después, reconoció los beneficios del método de laboratorio y lo comenzó a recomendar, pero procediendo primero bajo un periodo de prueba a modo de preparación. Se desconoce hasta que punto llegó a implementarse pero los comentarios de la gente de aquella época manifestaban que “a pesar de que el resultado fue malo, el pseudo-método basado en temas (“pseudo-complex-based approach”) fue utilizado prácticamente en todos los sitios” (Traducción del autor de Karp, 2012, p. 558). Por otra parte, y a pesar de que la prioridad de este nuevo sistema no era el de proporcionar una educación matemática sólida, hay que recordar que con él fue con el que se educaron los creadores del satélite Sputnik (Karp, 2010a, p. 52). En cuanto al número de estudiantes de primaria y secundaria, durante esta etapa se pasó de 7,8 millones en 1914 a 20 millones en 1931.

### 1.3.2. La educación matemática en Rusia (URSS) de 1931 a los años 50: la vuelta al sistema clásico.

Entre 1931 y 1936, y poco después de que Iósif Stalin se hubiera hecho con el poder absoluto de la Unión Soviética (1927), se aprobaron una serie de resoluciones que deshicieron la gran reforma que se había llevado a cabo en el sistema educativo después de la revolución de 1917 declarándolas “distorsiones izquierdistas”. En ellas se argumentaba que el sistema educativo no era capaz de proporcionar un conocimiento básico y general adecuado ni los conocimientos científicos necesarios para ingresar en escuelas técnicas e instituciones de educación superior (Citado en Karp, 2012, p. 558). Todo ello supuso una vuelta al sistema clásico, eliminando el sistema basado en el estudio por temas y la utilización del método de laboratorio (los propios educadores habían comenzado una lucha por su desaparición), con la consecuente imposición de control en el trabajo de los profesores y del desarrollo del curso diario de la actividad escolar. Consecuentemente se incrementó considerablemente la rigidez del recién sistema donde además se hacía ineludible el estricto cumplimiento del Currículo, que entre otras medidas hacía obligatorio el uso del libro de texto estipulado para cada asignatura en toda la Unión Soviética. Los libros de texto que se utilizaron fueron los redactados por Andrey Kiselev en el periodo anterior a la revolución. También se centralizó la evaluación de los alumnos, con la redacción desde Moscú de exámenes iguales para toda la población. Otra de las cuestiones que se tuvo muy en cuenta, fue la lucha contra el fracaso escolar que rondaba el 20% y que en amplio número se refería a la enseñanza matemática. Al mismo tiempo se trató de evitar que el nivel de los cursos estuviera devaluado, no correspondiendo con el grado del curso en el que él alumno se encontraba. Para ello la estructura de las clases y el comportamiento del profesor se controlaron de forma exhaustiva, y se dictaron indicaciones concretas al respecto, como por ejemplo: “proponiendo una revisión de la tarea para casa durante los primeros 10-15 minutos, llamando al alumno al encerado y examinando su cuaderno para señalarle sus errores” (Traducción del autor del texto citado en Karp, 2010a, p. 59). Además, a los alumnos aventajados se les recomendaba trabajo extra-curricular (Círculos Matemáticos, Olimpiadas en Leningrado en 1934, publicaciones), y a los alumnos más retrasados se les obligaba en general a asistir a clases de refuerzo. A las clases con un nivel bajo se les asignaba tarea para casa extra, el repaso obligatorio de todos los ejercicios de los exámenes, incluso se les prolongaba con una hora extra la clase, además de organizar y supervisar el estudio a través de las organizaciones con carácter socio-político de los grupos **Jóvenes Pioneros** y **Komsomol**. La consecuencia de esta presión sobre el estudio de los alumnos fue que muchos de ellos tenían que dedicar del orden de 5 a 6 horas diarias de trabajo en casa (Citado en Karp, 2010a, p. 77).

La importancia de los psicólogos y pedagogos de la etapa anterior fue completamente desechada y en su lugar se puso el acento en la preparación matemática de los profesores de forma continuada. Los profesores además, estaban sometidos a un estricto control de su trabajo, primeramente por parte del director y del subdirector que debían observar de forma constante las clases de los profesores para dejar constancia de ellas en un informe; segundo, eran sometidos al control, ellos y la escuela, por parte de los matemáticos especialistas en metodología; y por último al control de los inspectores, que no tenían por qué tener conocimientos matemáticos pero que podían hacerse acompañar por alguien que si los tuviera. Los inspectores se encargaban de criticar duramente la falta de contenido matemático en las clases y el mal uso del tiempo, con lo que no se conseguía el aprendizaje necesario. Si se consideraba que el profesor no estaba lo suficientemente preparado, se le obligaba a realizar cursos formativos o incluso podía ser despedido.

A partir de los años 40 aumentó la preocupación por la importancia en la resolución de problemas con un nivel alto en dificultad, de la manera más óptima posible, lo cual puso el foco en la exposición de las explicaciones y no simplemente en las soluciones. La educación matemática asociada a la resolución de problemas resultantes de supuestos prácticos, estaba asociada a la llamada **educación politécnica,** la cual había sido prácticamente eliminada a partir de 1932. A partir de 1945 comenzó otra vez a tenerse en cuenta, hasta el punto de llegar a plantearse como diseñar los cursos de matemáticas de una forma más "politécnica", la cual tendría su influencia sobre el planteamiento de la enseñanza matemática unos años más tarde, a finales de los años 50 y principios de los 60. No obstante, a lo largo de este periodo se incrementó de forma constante el nivel de exigencia, por ejemplo, entre 1940 y 1949 un miembro del **Instituto para la educación continua de los profesores de Leningrado**, describe como los problemas de geometría eran resueltos en un primer momento sin ninguna justificación, y como después, la propia explicación determinaba la calidad de la solución. Se pasó así de una calificación A para una solución correcta pero sin explicación, a una calificación C. Algunos matemáticos profesionales, como Kolmogorov, habían venido involucrándose en la educación escolar cada vez más, pero finalmente el empleo de estos matemáticos como profesores no dio buen resultado debido a su actitud negativa frente a las metodologías y a su falta de conocimiento de las escuelas de secundaria, por lo que se dejó de llamarles (Citado en Karp, 2010a, p. 65). En 1943, se instauró la separación por géneros en el sistema educativo, la cual se aboliría tras la muerte de Stalin.

En resumen, se volvió durante este periodo al sistema clásico pre-revolucionario, aunque esta vez se hacía posible el acceso a gran parte de la población, y la línea del sistema educativo iba en consonancia con el régimen totalitario de Stalin, adoptado medidas centralizadoras, represivas y de control estricto, tal y como hemos visto. Por su parte, la enseñanza de las matemáticas se encontró en una situación privilegiada frente a otras asignaturas, ya que por sus propias características se vieron libres de tener que someterse al control ideológico, además de ser consideradas como un saber de primer orden necesario para el desarrollo industrial, tecnológico y militar al que se sometió a la Unión Soviética durante estos años (Criado, 2011). La forma en la que se aplicó este nuevo programa puede ser resumido por lo dicho por los participantes en un congreso sobre cuestiones metodológicas en 1952: "al principio, el seguimiento del programa era opcional, después se convirtió en algo obligatorio, y finalmente el programa se seguía sin más" (Traducción del autor del texto citado en Karp, 2010a, p. 58).

### 1.3.3. La educación matemática en Rusia (URSS) de los años 60 a los años 70: las reformas de Kruschev y de Kolmogorov.

A finales de los años 50, y con Nikita Kruschev como dirigente de la Unión Soviética tras la muerte de Stalin en 1953, se produjo una primera reforma conocida como "**la reforma de Kruschev"** (1958) que iba a establecer los primeros pasos hacia la otra reforma, de mucha mayor importancia, que acontecería a partir de 1966: **"la reforma de Kolmogorov"**. La reforma de Kruschev puso el énfasis en un mayor desarrollo de la llamada educación politécnica, que trataba de conciliar en igualdad de importancia el trabajo tanto intelectual como físico y que tenía sus bases en la filosofía marxista-leninista, cuyo fin era el de vincular la escuela a los trabajos que luego se realizarían en la sociedad. La reforma de Kolmogorov, por el contrario trataba de ir mucho más allá, y quería lograr una modernización del sistema educativo, precisamente vinculando esta educación politécnica (conservando su carácter práctico) con la educación académica de modo que propiciara el desarrollo científico. El germen de esta idea hay que verlo en las importantes reformas llevadas a cabo tras la revolución de 1917, no en vano Kruschev sometió a la Unión Soviética durante su mandato al conocido proceso de desestalinización.

La reforma de Kruschev introdujo como obligatorio un periodo de 8 cursos en lugar de los 7 hasta el momento, y a su terminación los jóvenes se deberían incorporar al trabajo productivo industrial o agrícola. La educación secundaria se podía realizar pero combinándola con el trabajo productivo (Shabanowitz, 1978, pp. 42-43). También, durante los años 1962 a 1967 se prolongó la educación secundaria con el curso de grado 11. El resto de reformas fueron poco significativas, consistiendo en la aprobación de nuevos libros de texto de geometría (I. N. Nikitin) y álgebra (A. N. Barsukov), y la inclusión de la asignatura de trigonometría en el Currículo. La importancia de las matemáticas estaba totalmente asentada durante estos últimos años 50 y quedaba reforzada por importantes logros tecnológicos como el lanzamiento del **satélite Sputnik** en 1957 o el **primer vuelo espacial tripulado** en 1961. Es por ello que se debió de empezar a considerar la necesidad de prestar una educación adecuada para los alumnos más talentosos o dotados no sólo en matemáticas y en ciencias, sino también en música y artes en general, lo que supuso la aparición de las primeras **escuelas especializadas**. La propuesta se desarrolló en profundidad con la reforma de Kolmogorov, y permitió establecer la individualización y diferenciación dentro del sistema educativo que en la etapa anterior había permanecido como una “educación de las masas” bajo la mirada de un control estricto.

En 1962 se realizó una competición sobre la redacción de libros de texto de matemáticas, al que se presentaron 86 grupos de autores. El primer premio lo consiguió el libro de álgebra para cursos superiores de E. S. Kochetkova que además introducía elementos de cálculo. Otros grupos recibieron también premios y reconocimiento, y de alguna manera se estaban fraguando los inicios de una nueva reforma. Se abrió también durante estos años un debate entre profesores universitarios donde venía a exponerse la necesidad de actualizar los cursos de matemáticas. Se publicaron numerosos artículos en revistas sobre educación matemática y las principales observaciones que se propusieron giraban en torno a la necesidad de introducir en los cursos superiores nuevos temas de actualidad como elementos de cálculo, geometría analítica, álgebra vectorial y transformaciones geométricas. Además aparecieron publicaciones sobre como introducirlos, como la publicación a gran escala para el curso de grado 9 en 1964 a cargo V. G. Boltyansky e I. M. Yaglom, sobre transformaciones geométricas y vectores.

Finalmente en Diciembre de 1967, el **Ministerio de Educación de la Unión Soviética**, adoptó un nuevo Currículo que constituyó una reforma de especial importancia para la enseñanza de las matemáticas, afectando a todo el sistema educativo (la reforma de Kolmogorov). El Currículo lo redactaron científicos, especialistas en metodología y profesores. Entre los matemáticos que intervinieron estaban V. G. Boltyansky, I. M. Yaglom, A. I. Markushevich y A. N. Kolmogorov, y la reforma suele ser conocida por el nombre de este último por la influencia tan grande que tuvo en el desarrollo de la educación matemática durante estas décadas para la Unión Soviética, además de ser uno de los principales matemáticos del s. XX. Para elaborar el Currículo se creó un Comité Central presidido por el prestigioso y conocido matemático A. I. Markushevich, cuyas buenas relaciones con A. N. Kolmogorov, al frente del comité de matemáticas, contribuyeron a su desarrollo. Se puso el foco en eliminar todo aquello que había quedado atrasado y en dotar a la enseñanza del mayor rigor científico posible, a su vez motivado por la necesidad del desarrollo tecnológico de la Unión Soviética en su lucha en la Guerra Fría frente a los Estados Unidos. Un miembro de la Academia de Ciencias de la República Socialista de Ucrania señalaba que:

> La preocupación por los problemas en la educación matemática en todo el mundo no es casual; ha surgido como consecuencia de los importantes cambios en nuestra consideración hacia el desempeño de la ciencia en el progreso social. Es característico de nuestro tiempo el rápido desarrollo de los fundamentos del conocimiento científico, el veloz cambio de los avances técnicos, la amplia estandarización del trabajo intelectual y físico, y la necesidad de las matemáticas no sólo en la ciencia, sino también de la mayor parte de las actividades de carácter práctico.
> (Traducción del autor del texto citado en Shabanowitz, 1978, p. 47)

El primer documento, ***“El ámbito del conocimiento matemático para la escuela de primaria (8 cursos)”*** fue publicado en 1965, y en él se presentaban los hitos que debían alcanzar los alumnos una vez terminada esta etapa, lo que permitió establecer una discusión de carácter público en torno a ellos. Kolmogorov, por su parte se encargó de la sección de aritmética y álgebra, y también preparó otro para los cursos de grado 9 y 10 pero no llegó a publicarse. Poco después se publicó ***"el borrador del Currículo y la programación para las escuelas de secundaria"***, que después de ser publicado a diferentes escalas (folletos, revistas) y de su discusión pública, terminó por publicarse en su versión definitiva en 1968.

El Currículo introdujo nuevos temas como elementos de cálculo diferencial, transformaciones geométricas, vectores y coordenadas, aunque el principal cambio consistió en someter a la enseñanza a un proceso de "bourbakización", es decir, el reemplazo de la enseñanza intuitiva y con cierta formalización principalmente a través de ejemplos, por una visión puramente formal basada en la teoría conjuntista, lo cual condujo a la introducción del método axiomático y a la exposición puramente lógica durante las clases (Gordienko, 2011, p. 10). Precisamente debido al exceso de formalización, tanto estudiantes como algunos profesores se sentían perdidos, y ello condujo a grandes problemas en la enseñanza (Gordienko, 2011, p. 11). En cualquier caso el nuevo currículo introdujo cambios significativos en cuanto a la estructuración y selección de los contenidos: se redujo el curso elemental de 4 a 3 años; se introdujo una única asignatura bajo el título de "matemáticas" para los cursos de grado 1 a grado 5; se ofrecieron asignaturas optativas a los estudiantes; se eliminaron la estadística, la probabilidad y los números complejos; se eliminó todo el contenido que tradicionalmente se presentaba de una forma aislada; algunos temas típicos del Currículo de matemáticas se presentaron utilizando nuevos métodos. Todos estos cambios venían acompañados de la necesidad de una notación más rigurosa y de textos matemáticos que cumplieran con una exposición adecuada. A partir de su aprobación, Kolmogorov y otros autores publicaron artículos en revistas para dar a conocer las ideas introducidas. Algunas editoriales también publicaron libros y folletos bajo el título de ***"La nueva escuela de matemáticas"***.

En esta etapa de modernización, que constituía una enseñanza más académica y científica, aunque sin olvidar el carácter práctico de la educación politécnica, se recuperaba también la importancia de las cuestiones psicopedagógicas puestas en valor en los años 20 y que habían sido olvidadas concentrándose tan sólo en la formación del profesorado. De esta forma se centraba la atención en el alumno, al que además se le iban a proporcionar gran cantidad de opciones para su adecuada formación, como la elección de asignaturas optativas además de los proyectos ya iniciados algunos años antes consistentes en la creación de escuelas internado especializadas en física y matemáticas, de escuelas especializadas en matemáticas y física, y de escuelas de matemáticas por correspondencia, junto a la ya larga tradición de enseñanza extra-curricular. Todo ello dotaría a la Unión Soviética en los años 70 de un programa para el estudio de las matemáticas sin precedente, dando lugar a un auténtico boom de las matemáticas hasta el final de la década. El nuevo Currículo se implementó gradualmente avanzando por cursos hasta alcanzar su total desarrollo en el año 1975, con un seguimiento y evaluación continua así como las revisiones necesarias.

### 1.3.4. La educación matemática en Rusia (URSS) durante los años 80: la contra-reforma y el Currículo de 1981.

El Currículo de 1968, aunque llegó a implantarse en su totalidad, nunca estuvo carente de críticas y contó siempre con su grupo de opositores. Pero fue al poco tiempo de ser implantado en su totalidad (1978), coincidente con la problemática sobre como redactar los exámenes de ingreso en la Universidad, cuando surgió una oposición de contrarreforma liderada por los matemáticos I. M. Vinogradov, A. N. Tikhonov y L. S. Pontryagin, que tenían discrepancias sobre la visión de la educación matemática de Kolmogorov, además de que las relaciones entre ellos no parecían estar en el mejor de los momentos. (Abramov, 2010, p. 126) La idea de una contrarreforma se vio reforzada por la intervención del ministro de educación de la RSFSR (Rusia), A. I. Danilov, que salió en su defensa frente al **Ministerio de Educación de la Unión Soviética**.

En Mayo de 1978, durante una discusión sobre cuestiones que afectaban a la educación en la **Academia de Ciencias** y donde estaban presentes ambos bloques de matemáticos con visiones diferentes, entre los que se encontraba Kolmogorov, se propuso establecer una contrarreforma y se estableció una reunión para su discusión en Diciembre de 1978. En ella expusieron sus visiones diversos matemáticos entre los que se encontraban Kolmogorov, crítico con algunas cosas, y L. S. Pontryagin y A. N. Tikhonov que arremetieron contra el Currículo de Kolmogorov (Abramov, 2010, p. 127). Finalmente la resolución adoptada fue que se reconocía como insatisfactorio el actual estado del Currículo y de los libros de texto empleados, lo que concluyó en la formación de un comité sobre educación matemática presidido por I. M. Vinogradov, que fue apoyado por el **Ministerio de Educación de la RSFSR** (Rusia) para la elaboración de un nuevo Currículo y de libros de texto. En 1982, el **Partido Comunista**, influenciado por I. M. Vinogradov y A. N. Tikhonov, dio paso a que algunos libros de texto fueran modificados. L. S. Pontryagin por su parte se encargó de publicitar la contrareforma a través de la publicación en 1979 en el periódico **Industria Socialista** (Sotsialisticheskaya industriya) de un primer artículo titulado ***"Ética y Aritmética"***, y de un segundo artículo en 1980, que tuvo mayor repercusión, en **Comunista** (Kommunist), el principal periódico ideológico del **Partido Comunista**, donde la reforma era elevada a una cuestión de carácter ideológico. En el primero de ellos se acusaba indirectamente a Kolmogorov de irresponsabilidad e inmoralidad (Abramov, 2010, p. 128). Finalmente en una reunión en 1980, L. S. Pontryagin reconoció que la contrarreforma no podía ser implantada ya que la renovación de los libros era, por el momento, inviable. M. A. Prokofiev, ministro de educación de la Unión Soviética, también entendió que los cambios no eran posibles ya que no había una alternativa seria preparada, a pesar de que en 1979 se habían publicado dos Currículos alternativos a cargo de I. M. Vinogradov y de A. N. Tikhonov (Abramov, 2010, p. 128). Sí se realizaron algunos cambios de libros de texto, introduciendo en 1980 los escritos por A. V. Pogorelov sobre geometría elemental y recomendados además por el propio Kolmogorov. El grupo opositor del Currículo de Kolmogorov terminó por dividirse y sus miembros se dedicaron a dirigir la publicación de libros de texto. A. N. Tikhonov pasó a encabezar la dirección de los grupos de autores de libros de texto creados en el **Ministerio de educación de la RSFSR** (Rusia). El comité de educación matemática presidido por I. M. Vinogradov, y después de su muerte en 1982 por L. S. Pontryagin, apoyó por su parte a otros grupos de autores. A. D. Aleksandrov, sustituyó a Kolmogorov como director del Consejo Metodológico Científico y se convirtió también en el director de un grupo de autores de libros de texto. Todo esto dio lugar a que se empezara a considerar la posibilidad de utilizar diferentes libros de texto frente a la tradición mantenida hasta ahora del uso de un único libro de forma extensiva (Abramov, 2010, p. 131).

En 1981 se redactó un nuevo Currículo elaborado por V. V. Firsov, N. N. Reshetnikov y Alexander Abramov, que incluía prácticamente el mismo cuerpo de teoría que el de 1968, ya que era coincidente en un 90% con las otras dos propuestas de contrarreforma que se habían elaborado (Abramov, 2010, p. 132). Por otra parte se habilitaba la posibilidad del uso de diferentes libros de texto, pero se debía impartir obligatoriamente el programa propuesto, lo cual se explicitaba a través de la sección "Contenidos de la Educación". Para el uso de los libros de texto además se redactaba una propuesta didáctica específica por áreas, incluida en la sección "Planteamiento de la asignatura". El Currículo se redactó por primera vez por áreas de conocimiento matemático, estableciendo una propuesta metodológica para cada una de ellas.

## 2. La educación matemática extra-curricular y extra-ordinaria en la Unión Soviética entre 1930 y 1985.

Entre 1930 y 1985 la educación extra-curricular y la educación extra-ordinaria, que tomaría en cierta medida el relevo de la extra-curricular a partir de ***"la reforma de Kolmogorov"*** **de** 1967, aunque ambas coexisten, dotó a la Unión Soviética de una formación matemática de especial relevancia y de un interés generalizado hacia las matemáticas por parte de su sociedad. Los mejores matemáticos rusos nacidos después de 1930 fueron educados en ella, y es conocida su amplia repercusión internacional, todo ello gracias en parte a la gran implicación de la mayoría de los matemáticos rusos, en la formación de las nuevas generaciones a través de nuevas fórmulas educativas más allá de la educación ordinaria, que trataban de profundizar en la adquisición del conocimiento matemático por parte de los estudiantes (Skvortsov, 1978, p. 1). La situación durante estos años, relativa a esta educación, ha sido descrita:

> Existía una buena tradición en la enseñanza de las matemáticas. En general, el nivel de la educación matemática soviética era superior al de Estados Unidos. Además, para los alumnos talentosos de primaria y secundaria, había un sistema especial de entrenamiento matemático a través de Círculos Matemáticos supervisados por los profesores y estudiantes de la Universidad. Cada año se organizaban las Olimpiadas Matemáticas para varios niveles. Siempre recuerdo de manera muy agradable uno de estos Círculos (perteneciente a la Universidad de Moscú) y las Olimpiadas Matemáticas de Moscú. Los departamentos de matemáticas de algunas Universidades son dignos de ser recordados. Por ejemplo, en los años 60, en el departamento de matemáticas de la Universidad de Moscú, la concentración de matemáticos de primer orden del s. XX no se podía comparar con ningún otro lugar en el mundo (A. Kolmogorov, I. Gelfand, L. Pontryagin, P. Aleksandrov, A. Tikhonov, L. Ljusternik, V. Arnold, Yu. Manin, S. Novikov, Ya. Sinai, R. Dobrushin y muchos otros). Había un sistema de seminarios presididos por científicos de este calibre, y muchos estudiantes participaban en ellos. Era una oportunidad excepcional de conocer a tales figuras legendarias y de verte involucrado en la investigación desde las primeras etapas.
> (Traducción del autor de Polyak, 2002, p. 2)

### 2.1. La educación matemática extra-curricular.

La educación matemática extra-curricular en Rusia (Unión Soviética) cuenta con una larga tradición y puede verse como una subcultura, ya que trata de formar a sus propios participantes con el objetivo, al margen del de su educación, de que tomen el relevo y se ocupen en el futuro de seguir desarrollando estas actividades, incluso expandiéndose hacia otras regiones, creando nuevos nodos (su actividad solía concentrarse en las ciudades más importantes como Moscú y San Petersburgo) que ayuden a ampliar la red (Saul y Fomin, 2010, pp. 225, 229). Es por tanto fruto, tal y como mencionan Fomin, Genkin, y Itenberg (2012) de unas circunstancias culturales dignas de ser consideradas, y estudiadas, y de las cuales se puede aprender mucho. Estaba orientada a jóvenes con un especial interés hacia las matemáticas, que con independencia de sus buenos resultados en la educación ordinaria, querían disfrutar "haciendo matemáticas" en un ambiente distendido y participativo a través de la realización de actividades cuidadosamente diseñadas tanto en sus contenidos como en su forma de llevarse a cabo. Son varios los formatos bajo los que se desarrolló, y aunque cada uno de ellos tiene sus propias características, todos ellos giran en torno a la resolución de problemas, cuestión que ha tenido una gran importancia histórica en el transcurso de la elaboración de la propia matemática, y que dentro de la investigación matemática rusa ha constituido también una de las principales preocupaciones (Karp y Leikin, 2010, pp. 228, 422). A través de estas propuestas, se trata principalmente de fomentar en los jóvenes estudiantes, la adquisición de habilidades que les permitan enfrentarse a la resolución de problemas, capacitándoles y haciéndolos más autónomos y seguros de sí mismos frente al estudio. Dentro de sus posibles formatos encontramos: los Círculos Matemáticos, las Olimpiadas Matemáticas, los Campamentos Matemáticos y los Festivales Matemáticos, aunque todo ello debe entenderse conjuntamente, tal y como veremos más adelante, ya que estaba estrechamente interrelacionado. Había también otro tipo de actividades paralelas con carácter divulgativo y formativo que trataban de despertar en los alumnos de las escuelas ordinarias su interés hacia el "mundo matemático" e incentivar su participación en los Círculos Matemáticos. Es el caso de las publicaciones en revistas de gran tirada y precio económico, el tablón matemático que había en todas las escuelas y las conferencias por parte de matemáticos dedicados a la investigación. Más tardíamente (años 70 y 80) se realizarían también experiencias muy importantes como la publicación mensual de **La Revista Kvant**, **La Competición Lomonosov** y **El Torneo de Ciudades** creados por Nikolai Konstantinov o la **Universidad para Judios** de Bella Subbotovskaya con una metodología próxima a los Círculos. También cabe considerar toda la literatura matemática que se generó en el desarrollo de estas actividades como las series ***"Conferencias populares de matemáticas"*** y ***"La biblioteca del Círculo Matemático"*** entre otras, enfocadas hacia el fomento del trabajo autónomo del lector y utilizadas también en la programación de las sesiones de los Círculos Matemáticos (Wirszup, 1963, p. 200).

### 2.1.1. Los Círculos Matemáticos.

Los inicios de los denominados **Círculos Matemáticos** se sitúan antes de la revolución de 1917, aunque es a partir de esta fecha cuando comenzaron a cobrar una mayor importancia y obtener gran popularidad frente a la devaluación de la educación ordinaria, que se vio sometida a una dura crítica (Marushina y Pratusevich, 2010, p. 375). A partir de los años 30, con la reestructuración de la educación y la vuelta al sistema clásico, la educación extra-curricular pasó a ser objeto de interés por parte de eminentes matemáticos de la época (A. N. Kolmogorov, B. N. Delone, L. A. Lyusternik) que se hicieron cargo de ella y en adelante llegaría a alcanzar una repercusión sobresaliente ya que proporcionaba un papel determinante en la formación de los futuros matemáticos.

Los Círculos Matemáticos heredan la tradición de los **kruzhok** (кружок) o Círculos, que consistían en la reunión de forma libre de un grupo de gente en torno a un tema de interés (política, literatura, ciencia) con una periodicidad normalmente semanal e impulsados por alguien. En herencia de ese carácter contaban con cierta independencia de las instituciones oficiales educativas y sus impulsores no recibían remuneración alguna, aunque sí contaban a su disposición con las aulas de las escuelas y universidades. Se formaron Círculos para todos los cursos (grados 5 a 10) en las principales ciudades y eran conducidos por profesores de Universidad, estudiantes graduados o en sus últimos años de estudios. En las ciudades de Moscú y San Petersburgo se contaba con una extensa red de Círculos en los años 50 y 60 que era muy popular, y seguirían funcionando hasta los años 70 y siguientes, pero con una importancia menor debido a la aparición de nuevos programas como las escuelas internado de matemáticas y física (Sossinsky A., 2010, p. 191).

Según Wirszup (1963):

> Cada Círculo se reúne una vez a la tarde cada semana y cuenta con 15 integrantes. Bajo la dirección de un profesor, los estudiantes aprenden a leer matemáticas de forma crítica e independiente. Discuten asuntos matemáticos haciendo uso del lenguaje matemático, amplían los conocimientos adquiridos en clase, y, mediante la resolución de problemas de cierto grado de interés y dificultad, poco a poco van adquiriendo una posición activa frente a las matemáticas. (p. 194) [traducción del autor]

La formación de los Círculos podía realizarse de diferentes maneras, pero en ningún caso consistía en tener que pasar algún tipo de prueba o examen. Cuando se formaba el Círculo podían acudir todos los alumnos que estuvieran interesados en participar, y a partir de entonces, con el tiempo el Círculo se iría conformando con los alumnos realmente interesados ya que los que no, lo acabarían abandonando. Una vez formado, el profesor conductor del Círculo podía invitar a participar directamente en ellos a los alumnos, o también podían ser invitados a causa de haber obtenido buenos resultados en las Olimpiadas Matemáticas. La formación de un Círculo solía durar varios años, lo ideal era que durase hasta que los participantes se hubieran graduado en sus escuelas para poder entrar en la Universidad. En general, la situación del Círculo dependía del profesor conductor del mismo (Marushina y Pratusevich, 2010, p. 397). En los Círculos además del profesor-guía solía haber asistentes que podían ser estudiantes interesados en estas prácticas, pero que no contaban con ninguna experiencia, de tal forma que su participación no sólo ayudaba en el desarrollo de la dinámica del Círculo sino que les proporcionaba la formación necesaria para en el futuro poder abrir su propio Círculo (Marushina y Pratusevich, 2010, p. 399).

Para el desarrollo de las sesiones el profesor o guía utilizaba material que incluía contenidos que ofrecer a sus alumnos pero también metodología de planificación que le ayudasen a conducir la sesión. Un investigador asociado a la sección de matemáticas del Instituto de métodos de enseñanza de la Academia de Ciencias de la Pedagogía de la RSFSR (Rusia) establecía, de este modo, los siguientes objetivos para un Círculo Matemático:

i. Proporcionar un mayor nivel en la cultura matemática de los estudiantes a través de un estudio más activo y profundo de varias áreas de las matemáticas y su historia, y mediante la familiarización con algunos problemas de la matemática contemporánea, etc.
ii. Desarrollar las habilidades lógicas y la visión espacial de los estudiantes a través de la solución de problemas de interés y puzles, analizando sofismas y paradojas, etc.
iii. Crear un grupo de estudiantes activos y entusiastas que apoyen al profesor en su trabajo diario. Este grupo podría ayudarle en diferentes tareas, como ayudar a los estudiantes más retrasados, preparar soporte visual, ayudarle en la conducción de proyectos de mayor envergadura (tardes matemáticas para la escuela y el tablón matemático, por ejemplo), etc.
iv. Finalmente, descubrir, a lo largo de las sesiones del Círculo, estudiantes con habilidades especiales para las matemáticas. La actividad del Círculo debería ayudar a los estudiantes a desarrollar sus talentos naturales y a adquirir ciertas habilidades imprescindibles para el trabajo científico en el futuro.

(Traducción del autor del texto citado en Wirszup, 1963, p. 195)

Dentro de las sesiones del Círculo se incluían conferencias, sobre diferentes temas, por parte de matemáticos dedicados a la investigación, y en especial, se destinaba una importante parte de las sesiones a la preparación de las Olimpiadas Matemáticas, que constituían el punto fuerte dentro de las actividades extra-curriculares del alumnado de secundaria.

David Shkliarskii, un joven talentoso matemático que murió en la Segunda Guerra Mundial en 1942, es una figura clave dentro de la formación de los Círculos. Entre 1938 y 1941, transformó la estructura de los Círculos Matemáticos de Moscú, reemplazando la vieja práctica consistente en la redacción de artículos por la resolución sistemática de problemas de diferentes grados dificultad (Citado en Marushina y Pratusevich, 2010, p. 398). Los participantes del Circulo, comenzaron a trabajar, bajo este nuevo sistema, individual y colectivamente. El método de Shkliarskii, resultó lo suficientemente probado después de las cuartas Olimpiadas en 1938, donde su grupo ganó 12 premios, incluyendo 4 primeros premios, lo que hizo que fuera adoptado por los demás Círculos (Kukushkin, 1996, pp. 556-557). La figura es tan conocida en la cultura rusa que a la pregunta de “¿dónde puedo encontrar este problema?” se responde popularmente: “mira en Shkliarskii” (Traducción del autor de Kukushkin, 1996, pp. 556-557).

Como ejemplo de referencia de uno de estos Círculos se puede tomar el **Círculo Matemático de Moscú** que fue creado en 1935 por los conocidos matemáticos L. G. Shnirelman, L. A. Lyusternik y I. M. Gelfand. El segundo domingo de cada mes durante la tarde, profesores de la Universidad de Moscú y de otros institutos impartían conferencias de matemáticas para alumnos de cursos de grados 7-8 y grados 9-10 por separado y cada domingo a la tarde el Círculo se reunía dirigido por estudiantes del Departamento de Mécanica y Matemáticas de la Universidad, aunque en sus orígenes las sesiones fueron conducidas por sus propios fundadores (Wirszup, 1963, p. 198). Las sesiones del Círculo se desarrollaron en un principio en secciones temáticas (álgebra, geometría, teoría de números, etc.) pero a partir de finales de los años 50 se comenzaron a impartir bajo un carácter general, debido según un reportaje del profesor E. B. Dynkin (1958) a que "la formación redujo la adquisición de conocimiento matemático por parte de los participantes, desmotivándoles en su trabajo" (Citado en Wirszup, 1963, p. 200). El principal objetivo del Círculo consistía en conducirles hacia una relación activa con las matemáticas, para ello, el profesor-guía del Círculo introducía un tema y a continuación proporcionaba un problema que sus integrantes tenían que solucionar con un tiempo estimado de entre 5 a 10 minutos. A continuación se discutía y explicaba su solución, y el profesor-guía continuaba con la exposición del tema seleccionado (Wirszup, 1963, p. 200). En 1956-57 la **Universidad de Moscú** llegó a organizar 15 Círculos en las escuelas de secundaria de la región de Moscú.

### 2.1.2. Las Olimpiadas Matemáticas.

Las Olimpiadas Matemáticas de la Unión Soviética son el origen de las competiciones matemáticas en el resto de países. Su gran éxito fue debido a la implicación de numerosos matemáticos profesionales, más de 100, que diseñaron muchos, diversos e interesantes problemas por cursos o grados para los estudiantes de primaria y secundaria, en general a partir de 4º curso o grado (Slinko, 1997, p. 7). El formato de su organización ha variado, según el número y la forma de agrupar los problemas por un lado, y respecto a la edad de participación, aunque en general, a partir de 1952, los participantes compiten dentro de un mismo rango establecido por curso o grado (Kukushkin, 1996, p. 554). Lo más interesante, al margen de las mismas Olimpiadas, es el hecho de que se promovió todo un sistema de entrenamiento enfocado a las matemáticas y que el interés por ellas se incrementó, aunque normalmente sean consideradas principalmente como un sistema para detectar el "talento" matemático.

Las Olimpiadas Matemáticas tienen su origen en competiciones individuales en la resolución de problemas de matemática elemental en Rumanía, Hungría y Rusia, a finales del s. XIX (Skvortsov, 1978, p. 351). Hungría, fue el primer país, donde este tipo de competiciones fueron organizadas a escala nacional, bajo el nombre de **Competición Eötvös** en 1894, y que permitían a los 10 primeros el libre acceso a la Universidad sin necesidad de pasar las pruebas, además de que sus nombres fuesen publicados y recibieran reconocimiento nacional. A partir de ese año se celebraron anualmente (salvo en algunos años debido a la primera y segunda guerra mundial), y terminaría por servir de modelo a otras iniciativas y países en todo el mundo. Ese mismo año, también en Hungría, apareció una **Revista Matemática para Secundaria** en la que se compilarían los problemas utilizados en la Competición. El nombre de la competición se cambiaría en 1949 pasándose a denominar **Competición Kürschák** (Suranyi y Halmos, 1978, p. 291).

La organización de las primeras Olimpiadas Matemáticas rusas, tienen lugar en 1933 en Georgia, cuando los profesores Vashakidze y Petrakovskaja organizaron una competición en los colegios a nivel regional. En 1934, se organizarían a nivel local en Leningrado (actual San Petersburgo) y en Tbilisi, y en 1935 en Moscú y en Kiev (Slinko, 1997, p. 3). Las Olimpiadas fueron creadas bajo la predominante opinión de los educadores y matemáticos soviéticos de que en los últimos cursos o grados de educación secundaria, cuando los jóvenes comienzan a tener 13-14 años, es cuando mejor se pueden desarrollar las habilidades matemáticas, ya que de no hacerlo después sería mucho más difícil de conseguir (Skvortsov, 1978, p. 352). Este talento matemático a veces no se percibe bajo la enseñanza ordinaria, por ello, las Olimpiadas Matemáticas jugaban un papel fundamental en el diagnóstico de la posesión de estas habilidades. Con el establecimiento de esta competición se incrementó la actividad matemática en torno a ella, por ejemplo, después de la Olimpiada Matemática celebrada en Moscú, aumentaron considerablemente el número de Círculos Matemáticos (Saul y Fomin, 2010, p. 233). El establecimiento de las Olimpiadas Matemáticas en estas ciudades hizo también que en los años siguientes se propagara a otras, y a finales de los años 50 constituían una institución local en muchas ciudades. Los matemáticos Delone y Tartakovsky (Leningrado), Shnirelman y Lyusternik (Moscú) y Kravchuk (Kiev), fueron quienes en un primer momento se implicaron en la organización de estas Olimpiadas locales, entre otras cosas, elaborando problemas.

La tradición de las Olimpiadas en la Unión Soviética se caracteriza porque lo que se evalúa es la calidad en la formulación de la respuesta, que en función de ella recibe una cantidad de puntos, permitiendo al participante optar a una primera, segunda, tercera posición o mención, en función de los puntos adquiridos, por lo que varios participantes pueden ocupar un mismo puesto. Lo que se valora, por encima de que la respuesta sea correcta, es sobre todo el conjunto de razonamientos que se plantean para llegar a ella. En muchas ocasiones, por falta de tiempo no se puede llegar a desarrollar del todo, pero sin embargo es igualmente premiado atendiendo a su calidad. De hecho, en una ocasión se concedió un primer premio a un participante que no había conseguido ninguna respuesta correcta, pero cuyos métodos y estrategias para resolver los problemas eran realmente increíbles (Moser, 1967, p. 888). La organización de las Olimpiadas solía ir unida al desarrollo de un festival, donde otras actividades, como conferencias, excursiones, conciertos, presentaciones y análisis de las soluciones tenían lugar, además de constituir un lugar de encuentro y amistad entre sus participantes que contaban con intereses comunes. Aunque solían desarrollarse siempre durante unas determinadas jornadas, en realidad abarcaban un propósito mucho mayor o de conjunto, que se engloba dentro de toda la educación de carácter extra-curricular y que A. N. Kolmogorov definía con precisión en un folleto titulado ***"La profesión del matemático"*** (O professii matematika,1959):

> La enseñanza en las clases de carácter obligatorio está diseñada para la adquisición por parte de todos los alumnos de un sólido conocimiento de las matemáticas. Independientemente de ello, probar las habilidades que uno tiene mediante la resolución de problemas de gran dificultad y familiarizarse con el método científico para solucionarlos, así como con la aplicación de las matemáticas en la ciencia y en la tecnología, sólo se puede conseguir en un Círculo Matemático.
>
> …las Olimpiadas deberían corresponderse con el punto culminante de la actividad anual, y no como una celebración aislada…
>
> …el éxito en las Olimpiadas matemáticas se logra… cuando las actividades en los Círculos están bien organizadas…
>
> Los problemas proporcionados en los Círculos y en las Olimpiadas algunas veces tienen un carácter artificioso e incluso cómico. No hay problema con ello si los problemas escogidos entrañan un esfuerzo mental considerable, que requiera de una cierta madurez en la producción matemática.
> (Traducción del autor de Wirszup, 1963, pp. 196-197)

La RSFSR (Rusia), por su parte, publicaría en 1956 unas instrucciones dirigidas a la preparación y conducción de las Olimpiadas y competiciones en general, tanto de matemáticas como de física, historia, geografía y astronomía, que establecían los siguientes objetivos principales:

1. Descubrir a los alumnos con especial habilidad matemática y facilitarles sus futuros estudios.
2. Elevar el nivel de la educación matemática de los estudiantes, e indirectamente el de los profesores.

3. Estimular el gusto por las matemáticas en los jóvenes.
(Traducción del autor de Wirszup, 1963, p. 196)

Las Olimpiadas Matemáticas se iban a organizar además de a nivel local a una escala cada vez mayor. En 1961, se organizaron también en Moscú, las primeras **Olimpiadas Matemáticas de toda Rusia**, por parte A. N. Kolmogorov (quién dirigió por un tiempo el comité matemático de organización central), en las que se invitaba a participar a los ganadores del resto de Olimpiadas locales. Lavrent'ev repitió la idea en 1962 organizando las **Olimpiadas Matemáticas de toda Siberia**. El acontecimiento tuvo mucho éxito en la RSFRSR (Rusia) así que se repitió durante los siguientes años en Moscú y después se trasladó a otras ciudades (en 1966 en Voronezh). En ellas también participaban en ocasiones estudiantes de otras Repúblicas Soviéticas. El éxito acontecido se propagó incluso hacia otras disciplinas (física, biología, lingüística, geografía e informática) que comenzaron a desarrollar sus propias Olimpiadas. La participación de otros estudiantes provenientes de otras Repúblicas Soviéticas era cada vez más habitual, hasta que en 1967, y tras dar su apoyo el Ministerio de Educación de la Unión Soviética a la organización de las Olimpiadas, pasaron a denominarse **Olimpiadas Matemáticas de toda la Unión Soviética**, y en ese año se celebraron en Tblisi, capital de la República de Georgia, y después su sede se emplazaría de modo cambiante en las diferentes Repúblicas. Las Olimpiadas se celebraban a lo largo del año en 4 etapas:

1. A nivel de escuelas (grados 4-10), donde todos los alumnos pueden participar y se intenta incentivarles para que lo hagan. Los problemas los diseñaban empleados de Institutos educativos.

2. A nivel de ciudad (grados 5-10), donde siguen pudiendo participar todos los alumnos. Los problemas los diseñaba la Comisión de las Olimpiadas de cada ciudad.

3. A nivel regional (grados 7-10). Los problemas los diseñaba la Comisión de las Olimpiadas. En esta etapa tan sólo un participante pasaba a la siguiente. Hasta 1975 había tan sólo 3 niveles, y a la prueba final acudían los 4 finalistas de todas las Olimpiadas regionales celebradas en la Unión Soviética y los ganadores de la edición anterior. El número de personas involucradas en las Olimpiadas contando los miembros de los equipos, líderes, miembros del comité organizador, y profesores llegó a alcanzar la cifra de 800. Para reducir ese número se incrementó el número de exámenes (de 3 a 4), y a partir de ese año hubo que superar una prueba de corte en cada una de los distritos o zonas (Kukushkin, 1996, p. 554).

4. A nivel de distrito o zona (grados 8-10), donde toda Rusia quedaba dividida en 4 regiones más las de las ciudades de Moscú y Leningrado, que tenían normas diferentes, ya que el nivel y la participación en ellas era mucho mayor. También se invitaba a participar a los ganadores de la edición pasada. Hasta 1992, este examen era el último para las **Olimpiadas Matemáticas de toda Rusia**, después se desarrollarían en 5 exámenes.

Al examen final perteneciente a las fase de las **Olimpiadas Matemáticas de toda la Unión**, acudían los ganadores de las 6 zonas, donde los participantes tenían que resolver problemas especialmente difíciles, algunos de los cuales eran tan sólo resueltos por 1 o 2 participantes (Kukushkin, 1996, p. 560).

La última etapa de la competición se organizaba durante 2 días, y sólo los estudiantes de los 3 últimos cursos o grados (8-10), con edades entre los 15 y los 17 años, podían participar a este nivel. Cada día debían resolver 4 problemas en 4 horas, por lo que se necesitaban 24 problemas nuevos cada año, y los responsables de su diseño eran los miembros del **Comité Metodológico**, unos 20-30 matemáticos profesionales o educadores. El Comité estaba a cargo del Ministerio de Educación a partir de 1979 lo que hizo que perdiera su independencia, y que el clima informal y amigable que las caracterizaba se diluyera, aunque el excelente nivel matemático se mantuvo (Sossinsky A., 2010, p. 207). Los miembros de este Comité formaban la mitad del jurado, y la otra mitad estaba formada por matemáticos profesionales o educadores locales. El Comité, también realizaba la importante tarea de diseñar los problemas que servirían a modo de ejemplo y entrenamiento previos a la competición, un total de entre 20 y 25 problemas que el Ministerio de Educación se encargaba de distribuir por las Repúblicas Socialistas, y que cada República utilizaba como creía conveniente (Slinko, 1997, p. 5). Las Olimpiadas suscitaban una gran actividad e interés en torno suyo en toda la comunidad educativa, lo sucedido era cubierto por periódicos afines y las revistas **Kvant** y **Mathematika v Shkole** (Matemáticas en la Escuela) publicaban los problemas con sus soluciones (Skvortsov, 1978, p. 360). Al final de la competición, los miembros del jurado solían analizar los problemas y se organizaba una ceremonia de entrega de premios. Además, los ganadores tenían preferencia de entrada en cualquier Universidad siempre y cuando cumplieran los requisitos en el resto de asignaturas. Los seleccionados para representar a la Unión Soviética en las **Olimpiadas Internacionales** no tenían que realizar el examen de entrada (Skvortsov, 1978, p. 359).

Al margen de estas Olimpiadas de carácter oficial, hubo otras iniciativas paralelas, por parte de Universidades e Institutos técnicos, de desarrollar de modo independiente sus propias Olimpiadas con la intención de atraer a los estudiantes más capaces a ellas. La primera Olimpiada de este tipo fue convocada en 1962 por el **Instituto de Física y Matemáticas (MIMP)**, y se desarrollaron en las escuelas pertenecientes al Instituto presentes en 58 ciudades (Kukushkin, 1996, p. 554).

Las primeras **Olimpiadas Internacionales** se celebraron en Bucarest (Rumanía) entre el 20 y 30 de Julio de 1959, y fueron organizadas conjuntamente por la Sociedad de Matemáticas y Física y el Ministerio de Educación de este país. Se invitó a participar a 9 Repúblicas Socialistas de las cuales aceptaron 6 (Turner, 1985, p. 333). En los años siguientes se irían sumando no sólo países socialistas sino también otros países Europeos. En 1977, 21 países estaban representados: Algeria, Austria, Bélgica, Bulgaria, Gran Bretaña, Hungría, República Democrática Alemana (RDA), República Federal Alemana (RFA), Cuba, Mongolia, Los Países Bajos, Polonia, Rumanía, Unión Soviética, Estados Unidos, Finlandia, Francia, Italia, Checoslovaquia, Suecia y Yugoslavia. Al principio no había ningún comité organizador lo que supuso que finalmente en 1980 las Olimpiadas no se celebrasen (Turner, 1985, p. 333). A partir de 1981 cada país participante contaría con un Comité organizador y representantes del Ministerio de Educación. A ellos se les ofrecía un conjunto de problemas para que seleccionaran los mejores que se remitían al jurado internacional que se encargaba de seleccionar los definitivos, para lo que se requerían 2 o 3 días. Los problemas definitivos se editaban y publicaban, inicialmente, en cuatro idiomas: inglés, alemán, ruso y francés. Posteriormente, para que pudieran ser repartidos también a otros participantes cuyas lenguas maternas no eran esas, se traducían a tantos idiomas como fueran necesarios. Después de los dos días de participación, se corregían los problemas y el jurado se reunía para designar a los premiados. El primer, segundo y tercer premio eran especiales para las soluciones mejores y más creativas. Oficialmente la competición se realizaba de forma individual y no entre equipos aunque se conformaban grupos de 8 integrantes, a los que se les asignaba un número, para organizar la competición. En 1982 este número sería de 4, y a partir de ese año hasta la actualidad el número es de 6 (Turner, 1985, p. 333). Aunque no se competía por países, de modo extra-oficial se realizaba una clasificación entre ellos, que permitía no sólo establecer un ranking, sino principalmente poder comparar el nivel matemático entre los diferentes países participantes. A los participantes se les ofrecía la posibilidad de viajar, relacionarse y conocer los países donde se realizase la competición (Skvortsov, 1978, p. 370).

### 2.1.3. Los Campamentos Matemáticos y los Festivales Matemáticos.

En la Unión Soviética los jóvenes (mediante su participación en los Jóvenes Pioneros) asistían a campamentos de verano que eran subvencionados por el Estado a través de diferentes organizaciones. Los campamentos solían durar 3 o 4 semanas y en ellos se tenía cubierta toda la estancia, incluyendo comida y alojamiento. Los campamentos ofrecían un programa de entretenimiento y actividades deportivas, y existían también campamentos que ofrecían además de otros programas un curso intensivo de matemáticas o únicamente de matemáticas. Estos **Campamentos Matemáticos**, prolongaban la actividad de los Círculos Matemáticos durante el verano, es decir, mantenían en ellos el mismo formato y la única diferencia consistía en que se desarrollaba de manera intensiva (hasta 6 horas diarias dedicadas a las matemáticas) y en un ambiente totalmente distinto más propio de la época estival. La forma de organizarlos solía partir de una invitación por parte de los profesores-guías de los Círculos Matemáticos directamente a los jóvenes participantes en ellos, una vez que se había decidido un lugar y una fecha donde celebrar el campamento. A parte de resolver problemas según la metodología de los Círculos, a veces también se introducían nuevos temas a través de presentaciones teóricas. Un ejemplo actual podemos encontrarlo en la **Escuela nº 57 de Moscú**, cuyo programa se puede consultar en la URL: http://sch57.ru/mathcamp.en.html/

Otra forma usual de fomentar el interés por las matemáticas en la Unión Soviética consistía en la realización de **Festivales Matemáticos**, encuentros normalmente de un día, donde se daba lugar a un programa lúdico a la vez que con contenido matemático (conferencias, resolución de problemas, etc.). Un ejemplo popular de ellos son los **Festivales para Jóvenes Matemáticos** celebrados durante varios años en Batumi, gracias a la iniciativa de Medea Zhgenti, profesor local, y al apoyo del consejo editorial de la **Revista Kvant** (Marushina y Pratusevich, 2010, p. 404). En este festival era usual la redacción de artículos o informes por parte de los participantes, tradición heredada de los inicios de los Círculos Matemáticos, y que sería reutilizada en este tipo de encuentros. Los artículos eran leídos por parte de los participantes a un jurado formado principalmente por miembros del consejo editorial de la **Revista Kvant**. La forma de organizar un Festival podía ser variada pero siempre requería de una importante colaboración por parte de todos los organizadores.

### 2.1.4. Conferencias de matemáticos dedicados a la investigación.

La implicación, durante de la Unión Soviética, de los matemáticos vinculados con las Universidades y dedicados a la investigación, en la educación matemática de las nuevas generaciones fue fundamental para la consecución de los logros que se alcanzarían años más tarde. Parte de esta implicación, consistía en las conferencias impartidas por matemáticos dedicados a la investigación en las Universidades de Moscú, Leningrado, etc., a estudiantes de secundaria y participantes en los Círculos Matemáticos. En la **Universidad de Moscú**, multitud de jóvenes y sus profesores asistían dos veces al mes a estas conferencias impartidas por matemáticos de reputación internacional como P. S. Alexandrov, A. O. Gelfond, I. M. Yaglom y muchos otros (Wirszup, 1963, p. 198). Las conferencias se organizaban por niveles y varias de ellas se publicarían en diferentes volúmenes en la serie titulada ***"Conferencias populares de Matemáticas"***, y que formarían parte importante de las publicaciones matemáticas con carácter extracurricular en toda esta etapa de la Unión Soviética, que además serían utilizadas de forma habitual en los Círculos Matemáticos. Su publicación comenzó en 1950 y en 1961 contaba con 35 volúmenes (Wirszup, 1963, p. 199). Como ya hemos comentado, estas conferencias también se incluían dentro de algunas de las sesiones de los Círculos Matemáticos, y se centraban en torno a algún tema. Algunos ejemplos de ellas, fueron las impartidas en 1955-56 para los cursos o grados 7 y 8, dentro de los Círculos vinculados a la Universidad de Moscú: "Métodos para la suma de series (G. Ts. Tumarkin), Construcciones geométricas mediante el uso exclusivo del compás (V. A. Tonyan), Ordenadores electrónicos (M. A. Kartsev), Simetrías y movimientos en el plano (E. A. Morozova), Resolución de problemas matemáticos mediante métodos mecánicos (V. A. Uspenskii), Vectores y sus aplicaciones en geometría elemental (S. Ya. Khavinson), Teoremas de Configuración y Construcción de problemas (L. A. Skornyakov) y El significado geométrico de las ecuaciones e inecuaciones algebraicas de primer grado (S. Ya. Khavinson)" (Wirszup, 1963, p. 198).

### 2.1.5. La revista Kvant.

A. N. Kolmogorov e I. K. Kikoin fundaron esta revista en 1969, con el apoyo de Petrovsky por parte de la **Academia de Ciencias**. Su intención era la de popularizar aún más si cabe las matemáticas y sobre todo proporcionar abundantes problemas interesantes y estimulantes para los estudiantes que se encontraran en áreas rurales o pequeñas ciudades lejos de los principales núcleos matemáticos que constituían las grandes ciudades. Tuvo una tirada de gran importancia, pasando de las iniciales 100.000 unidades, hasta las 370.000 durante el boom de las matemáticas y la física a principios de los 70, y retornando a las 200.000 a partir de 1977 una vez que el boom había terminado, hasta la década de los 80 (Sossinsky A. B., 1993, p. 239). En ella participarían renombrados matemáticos con interés en la educación de los jóvenes como los propios Kolmogorov y Kikoin, S. Fomin, Fuchs, Gelfand, Kirillov, M. Krein, Pontryaguin, Rokhlin, Tikhomirov, Viro, Zalgaller, Arnold, Rokhlin, Migdal, Frank-Kamenetski, muchos de ellos formadores de ganadores en las Olimpiadas Matemáticas, que a su vez participarían también en la revista. Es conocida la anécdota ocurrida en 1980, cuando L. S. Pontryagin en su lucha por el cambio de la educación matemática llevada a cabo durante esos años, intentó hacerse fallidamente con el control de la revista.

La revista se organizaba en diferentes secciones donde se incluían artículos de temática muy variada, y que tenían un tratamiento interactivo con el lector. La parte matemática era editada por A. B. Sossinsky y L. Makar-Limanov. Una de las secciones más populares era la sección de problemas dirigida por Nikolai Vasiliev e ilustrada por E. Nazarov, en la que se publicaban problemas (con sus soluciones) que habían sido propuestos en los exámenes de entrada a las principales Universidades, y que ayudó en gran manera en la alta difusión de la revista (Sossinsky A., 2010, p. 202). Otra sección importante, realizada por Nikolai Vasiliev, era en la que se planteaba a los lectores la resolución de 4 problemas, cuya solución se publicaba en el número siguiente. Desde la revista se invitaba a los lectores a participar enviando sus soluciones, las cuales eran corregidas, pudiendo optar a un premio anual como mejor solucionador de problemas. Muchos de estos ganadores, confesaron a A. B. Sossinsky años más tarde, que de no haber sido por la **Revista Kvant** y las Olimpiadas Matemáticas, no se hubieran dedicado profesionalmente a las matemáticas (Sossinsky A. B., 1993, p. 240). Es el caso de R. Bezrukavnikov (actualmente profesor en el Instituto Tecnológico de Massachusetts, MIT), I. Itenberg (profesor en la Universidad de Estrasburgo), I. Arzhantsev (profesor asociado en la Universidad Estatal de Míchigan, MSU), y G. Perelman (ganador de la medalla Fields en 2006) (Sossinsky A., 2010, p. 203). Todos los números de la revista original, en ruso, se puede encontrar en: http://kvant.mccme.ru/ y parte de un número de esta revista, en su publicación estadounidense, bajo el nombre de **Revista Quantum**, se puede ver en: http://www.nsta.org/quantum/qsampler.pdf/

### 2.2. La educación matemática extra-ordinaria.

La educación de carácter extra-ordinario ejerció una gran influencia en toda la Unión Soviética e internacionalmente, tanto en los modos de organizar la educación como en sus procesos didácticos (Karp, 2010b, p. 280-281). Sus orígenes se sitúan a partir de la **reforma de Kruschev** de 1958, cuando se produjo un aumento por la preocupación de la formación más especializada y el descubrimiento de los alumnos más talentosos según las áreas de conocimiento, entre las que se encontraban en un lugar privilegiado las matemáticas y la física por razones que ya se han comentado. Esta reforma condujo a la aparición durante los siguientes años de las **escuelas especializadas en matemáticas**, tomando como referencia las ya existentes, y bien conocidas, especializadas en artes (Shabanowitz, 1978, pp. 77-78). Su necesidad quedaba expresaba por el Comité Central encargado de la reforma, de la siguiente manera:

> Las escuelas y organismos de la educación pública deben prestar mayor atención al desarrollo de las aptitudes y capacidades de los alumnos tanto en artes como en matemáticas, física, biología y el resto de ciencias. Los Círculos, seminarios y conferencias tienen que estar perfectamente organizados dentro de estas instituciones. Deberían fundarse sociedades de jóvenes matemáticos, físicos,

> químicos, naturalistas y de técnicos en general, y se debería localizar el talento entre los jóvenes para que pueda ser debidamente atendido. Se deberá pensar seriamente la forma para el establecimiento de escuelas especializadas para jóvenes con buenas aptitudes en matemáticas, física, química y biología.
> (Traducción del autor del texto citado en Shabanowitz, 1978, p. 78)

Consecuentemente se inició un proceso de creación de escuelas especializadas durante los primeros años de la década de los 60, ampliándose a partir de 1963 con la incorporación de 4 **escuelas internado especializadas en matemáticas y física** que se abrieron inicialmente de forma experimental para la formación de los jóvenes más dotados (después se abrieron más en otros lugares), y la fundación de la **escuela de matemáticas nacional por correspondencia** en 1964, un proyecto que trataba de agrupar las diferentes escuelas puestas en marcha por correspondencia, para ofrecer una respuesta de calidad a los estudiantes de las zonas más alejadas de los principales núcleos universitarios que tuvieran la intención de prepararse para la entrada en la Universidad. La aparición de las escuelas especializadas e internados, hizo que la participación en los Círculos Matemáticos decayese (Moser, 1967, p. 888), ya que en ellas además se ofrecía un programa bajo el título de ***"actividades de la facultad"*** de forma extra-curricular consistente en seminarios bajo temas que eran elegidos a preferencia del profesor. Al igual que los Círculos, su carácter era voluntario.

A partir de 1966, se empezaron a materializar las ideas lideradas por Kolmogorov consistentes en un proceso de modernización y que se conocería como ***"la reforma de Kolmogorov"***. La reforma dio lugar a la implantación generalizada de un sistema de educación matemática extra-ordinaria, que hasta ahora no había consistido en una solución con carácter global, sino concreta. De esta forma se hacía accesible la posibilidad de un mayor estudio de las matemáticas previamente a la llegada a la Universidad, cuestión hasta ahora desarrollada únicamente de forma extra-curricular y que según expresaba Kolmogorov en una carta a A. I. Markushevich (29-12-1964) no parecía serle suficiente: "dejar todas las apuestas en los Círculos Matemáticos y en las escuelas especializadas de matemáticas no me parece muy prometedor" (Traducción del autor del texto citado en Abramov, 2010, p. 106). A partir de entonces se establecieron los **cursos optativos de matemáticas** en todas las escuelas que permitían la elección de determinadas asignaturas de contenido matemático según los intereses de cada uno. A partir de la década de los 70, ya en lo que se conoce como el *periodo de estancamiento* con Brezhnev en el poder, las escuelas especializadas e internados, pasarían a ocupar un segundo plano frente a los nuevos cursos optativos, que de alguna forma tomaron el relevo de estas escuelas, a la vez que se produjo un decrecimiento en ellas por el interés científico. En su lugar, estas escuelas pasaron a servir de puente de acceso a la Universidad para la gente de provincias, y la admisión en ellas se vio influida por el estatus social al que se pertenecía (Karp, 2010b, p. 275). Muchas de estas escuelas llegaron a cerrarse o era imposible acceder a ellas (Karp, 2009, p. 23). En 1985 y en los sucesivos años, bajo el Gobierno de Gorbachov y una vez establecida la perestroika, se intentó recuperar el antiguo éxito de estas escuelas haciéndolo además extensible mediante el establecimiento de un **curso avanzado de matemáticas** en prácticamente todas las escuelas, con lo que se trataba también de impedir la marcha de los alumnos más talentosos a otros centros (Karp, 2010b, p. 278). Se explica a continuación en que consistieron estas escuelas y los nuevos programas de matemáticas orientados a la educación del talento matemático.

### 2.2.1. Las escuelas especializadas en matemáticas.

Con la **reforma de Kruschev** de 1958, y debido a la importancia de la recién aparecida computación electrónica coincidente con la pretendida respuesta a la necesidad de una educación más especializada demanda durante aquellos años, el Ministerio de Educación de la Unión Soviética aprobó en 1959 planes experimentales para establecer programas especiales para formar a estudiantes en programación informática. En sus primeros momentos el plan se implantó en algunas escuelas de Moscú en los cursos o grados 9 a 11, dando lugar a las primeras **escuelas especializadas en computación**. En un principio se fueron formando programadores, hasta que los educadores se dieron cuenta de que esta formación en programación informática podía servir de estímulo para aprender más conocimientos matemáticos en vez de constituir un fin en sí misma (Moser, 1967, p. 887), así que estas escuelas evolucionaron hacia **escuelas especializadas en matemáticas** definitivamente a partir de 1964 en los albores de la **reforma de Kolmogorov**. Se llegaron a establecer en torno a 100 de estas escuelas en toda la Unión Soviética, donde el número de horas semanales dedicadas a las matemáticas era de 10 o 12, frente a las 5 o 6 dedicadas en las escuelas ordinarias (Moser, 1967, p. 887). En la década de los 80, se ampliaría su programa para los cursos o grados 7 y 8. Para su establecimiento, se hizo necesaria una preparación especializada de sus profesores, la cual la recibirían en institutos especializados a través de un programa de 5 años, en vez del programa de 4 años establecido para las escuelas ordinarias, que incluía formación en programación informática. Este programa servía también para ser docente en las escuelas internado de matemáticas y física (Moser, 1967, p. 889).

La primera escuela especializada fue la **Escuela nº 425**, en Moscú, y comenzó a funcionar en Septiembre de 1959 bajo la supervisión de S. I. Shvartsburd, incrementándose considerablemente el número de ellas en 1960-61 (Karp, 2010b, p. 268). En Julio de 1961, el Ministerio de Educación de la RSFSR (Rusia) aprobó la primera versión de los programas de estas escuelas especializadas, que comenzaron a formarse, aunque cada una de una manera en cierto modo independiente y diferenciadora del resto, lo que hacía que pudieran adquirir cierta fama. Un ejemplo conocido de implantación de este programa, es el de la **Escuela nº 2**, localizada en el suroeste de la periferia de Moscú, que comenzó a funcionar en 1956. Su director, Vladimir Ovchinnikov consiguió que los conocidos matemáticos Israel Gelfand, en 1963, y Evgeny Dynkin, en 1964, se incorporasen a la escuela como profesores gracias a que había aceptado previamente la entrada en ellas de sus hijos, y a que habían acordado que desempeñarían algún tipo de participación laboral dentro del centro. Gelfand y Dynkin terminaron no sólo por impartir clases sino también por desarrollar el programa de la escuela, incorporando además a otros jóvenes matemáticos como Alexey Tolpygo, Boris Geidman, Sergei Smirnov y Valery Senderov, lo que hizo que aumentase el nivel de enseñanza y acabase convirtiéndose en una escuela especializada (Sossinsky A., 2010, pp. 200-201). Otro ejemplo conocido es el de la **Escuela nº 7** en Moscú, que fue organizada como escuela especializada en 1962, por el matemático Alexander Kronrod, obteniendo gran éxito debido al alto porcentaje de entrada de sus alumnos en la Facultad de Mecánica y Matemáticas de la Universidad de Moscú (Sossinsky A. B., 1993, pp. 199-200).

### 2.2.2. Las escuelas internado de matemáticas y física.

En respuesta a la necesidad de proporcionar una educación adecuada a los jóvenes más dotados y dispersos en todas las provincias de la Unión Soviética, y en especial para los más alejados de las ciudades principales, se creó el proyecto experimental de las **escuelas internado de matemáticas y física**. La idea fue propuesta por los matemáticos M. A. Lavrentiev y A. N. Kolmogorov, quienes se encargarían de organizar los dos primeros proyectos de este tipo. B. V. Gnedenko en algunas de las conversaciones mantenidas con Kolmogorov, mencionaba que:

> A. N. Kolmogorov expresaba repetidamente el pensamiento de que muchísimos estudiantes talentosos de matemáticas de pueblos y zonas rurales estaban lejos del alcance de la comunidad matemática, ya que era imposible organizar Círculos y grupos especiales para su formación es las escuelas de secundaria rurales, así como proporcionar a estas escuelas profesores cualificados que se ocuparan de desarrollar un programa matemático con un carácter científico.
> (Traducción del autor de Karp, 2010b, p. 273)

Las primeras 4 escuelas de este tipo se abrieron en 1963 bajo el auspicio de las principales Universidades, en Novosibirsk, Moscú, Leningrado y Kiev. En Junio de 1964, el Ministerio de Educación de la Unión Soviética aprobaría una resolución para la regulación de estas escuelas, a partir de la cual surgirían otras 9 en diferentes ciudades soviéticas, cada una con su propia identidad, y bajo el auspicio también de sus Universidades. Los métodos de selección en ellas, diferían de unas a otras, siendo una forma frecuente hacerlo a través de las Olimpiadas. En las escuelas internados de Moscú, Kiev y Leningrado, los alumnos elegían una asignatura principal que podía ser matemáticas, física, química o biología, en su entrada en el curso o grado 9, no así en Novosibirsk, cuya elección podía ser aplazada al curso o grado 10. El nivel de estas escuelas era elevado, pudiendo considerarse que en matemáticas y física el nivel alcanzado se correspondía con un tercer curso universitario (Citado en Dunstan, 1975, p. 563). Eran comunes las clases magistrales impartidas por renombrados científicos y profesores universitarios seguidas por sesiones prácticas en pequeños grupos conducidas a su vez por otros profesores, que podían ser graduados universitarios o a punto de hacerlo, o profesores de escuelas con gran experiencia. Las escuelas contaban con todas las condiciones para desarrollar una enseñanza de gran calidad a alumnos que habían sido seleccionados por sus habilidades o talento, dando lugar a cierto carácter elitista, tal es así, que todas ellas llegaron a adquirir cierta fama ya que de ellas surgiría una parte importante de la posterior comunidad científica de la Unión Soviética.

La primera de ellas se fundó en Akademgorodok (cerca de Novosibirsk), un centro científico siberiano fundado a principios de los 60. En 1962 empezó a funcionar como escuela matemática de correspondencia, y terminaría por organizar un concurso para toda Siberia de matemáticas y física. Algunos de sus ganadores serían admitidos en el futuro proyecto de una escuela internado, la cual se acabó fundando en Enero de 1963, impulsada por M. A. Lavrentiev con la ayuda de A. A. Lyapunov (Sossinsky A., 2010, p. 193).

La más conocida se abrió poco después en Moscú, con el nombre de Isaac Newton, pero conocida como la **Escuela de Kolmogorov nº 18**, ya que fue Kolmogorov (junto a Kikoin) quien se encargó de desarrollar todo el programa lectivo, además de su gran implicación personal durante largo tiempo. En ella había 360 alumnos en los grados 9 y 10, y algunos estudiantes llegaban a tener hasta 20 horas de clase a la semana dedicadas a matemáticas y física de un total de 36 horas semanales de clase, lo que significa un amplio porcentaje (Moser, 1967, p. 887). El currículo de matemáticas incluía análisis, álgebra lineal, matemática discreta y teoría de la probabilidad. Los métodos numéricos y programación tenían un peso menor que en las **escuelas especializadas de matemáticas** (inicialmente escuelas de computación) (Shabanowitz, 1978, pp. 79-80). La implicación de Kolmogorov era de tal importancia, que él personalmente se encargaba de la enseñanza dos días a la semana a un gran grupo. Matemáticas, física, química y bilogía se impartían de esta forma, a través de clases magistrales al gran grupo, para el resto del tiempo lectivo semanal dividirse en pequeños grupos de discusión (de 30 alumnos) y resolución de problemas (de 15 alumnos), conducidos por profesores graduados en la Universidad. El resto de asignaturas se impartían de manera convencional (Dunstan, 1975, p. 562). En total, había 27 profesores a tiempo completo y 40 a tiempo parcial, que realizaban tareas de apoyo en charlas especiales y seminarios sobre diferentes temas de forma extra-curricular (Moser, 1967, p. 887), las cuales no se limitaban únicamente a las matemáticas y la física, sino que se incluían otras actividades dedicadas a las humanidades, como las tardes literarias, lecturas colectivas de autores clásicos y modernos, viajes al campo, etc. (Karp, 2010b, p. 285). La intención de esta escuela no era tanto la de formar matemáticos, sino más bien la de formar científicos con una sólida formación matemática. Su distintivo como escuela, a parte de la gran importancia de establecer un intenso proceso racional en el aprendizaje, era el de ser capaz de fomentar las capacidades de discernimiento y creatividad, a través de la investigación y la resolución de problemas no comunes (Dunstan, 1975, pp. 563-564). La forma de acceso a la escuela se hacía a través de un proceso de selección en 3 rondas. La primera consistía en una prueba escrita de matemáticas y física, realizada durante los mismos días de las Olimpiadas a nivel regional. Para poder realizarla se necesitaba haber sido recomendado por el profesor. La segunda ronda consistía en un examen oral, unas semanas más tarde, para los que hubieran superado la primera. De esta prueba, los seleccionados serían invitados a un campamento de verano durante 20 días, donde serían seleccionados en función de los resultados del trabajo desarrollado.

### 2.2.3. Las escuelas de matemáticas por correspondencia.

Las escuelas matemáticas por correspondencia se concibieron para llegar a las zonas rurales donde no había posibilidad de establecer escuelas especializadas, y constituían un intento de encontrar al talento matemático oculto en esas zonas rurales más apartadas, así como de proporcionarles un aliciente matemático adecuado para su desarrollo. Entre las más conocidas se encontraba la auspiciada por la **Universidad Moscú**, aunque también eran importantes las que estaban bajo el control de las **Universidades de Leningrado y Kiev**. Para entrar en ellas se debía realizar un examen de ingreso una vez superado el curso o grado 8, previamente al comienzo del nuevo curso. Los exámenes los preparaban y evaluaban los profesores de la escuela por correspondencia, que eran propuestos por la Universidad correspondiente. Este tipo de proceso a distancia era también habitual utilizarlo para el transcurso de las Olimpiadas, en aquellos lugares (como Siberia) donde era muy difícil organizar los encuentros regionales debido a las largas distancias.

En 1964, el conocido matemático Israel Gelfand, antiguo alumno de Kolmogorov en la **Universidad de Moscú** e inicialmente supervisor de Círculos Matemáticos, creó la **Escuela nacional de matemáticas por correspondencia** (VZMSh) bajo el apoyo de I. G. Petrovsky, el entonces rector de la Universidad de Moscú. El proyecto pretendía que todas las escuelas de correspondencia funcionasen de manera unificada para proporcionar una enseñanza matemática y facilitar el acceso a la Universidad de los estudiantes de las zonas más alejadas de las ciudades principales de la Unión Soviética. Gelfand diseñó el programa de esta escuela junto con otros matemáticos, como A. A. Kirilov, y escribiría también libros de textos para ella. La escuela era muy eficiente gracias a la buena cooperación que se daba entre el equipo de profesores liderados por Gelfand y los profesores a nivel provincial (Sossinsky A., 2010, p. 197).

Para su entrada, había que superar una prueba, y el curso consistía en un programa de 2 años donde se debían cubrir 20 tareas, las cuales eran corregidas por estudiantes de matemáticas pertenecientes al departamento de matemáticas de la Universidad de Moscú. Cada uno de ellos tenía asignados 10 alumnos de la escuela a su cargo, y el trabajo de cada uno de ellos, era supervisado por un estudiante graduado también perteneciente al departamento. Se tasan en 200.000 los alumnos graduados en esta escuela hasta 2010 pertenecientes a diferentes lugares (Abramov, 2010, p. 109).

I. Gelfand fue primero profesor asociado (1935-1940) y después profesor de la **Universidad de Moscú** (MSU) hasta 1990 cuando emigró a Estados Unidos donde trabajó como profesor para la **Universidad de Rutgers**. Una vez en ella, exportaría la misma idea estableciendo **el programa de matemáticas por correspondencia Gelfand** (GCPM) cuya descripción del programa se puede encontrar en la siguiente URL: http://gcpm.rutgers.edu/former_description.html/ El programa aún continua bajo el nombre de EGCPM (Extended Gelfand Correspondence Program in Mathematics).

### 2.2.4. Las clases opcionales en escuelas ordinarias.

Las escuelas especializadas e internados fueron aceptadas tan solo parcialmente y en torno a ellas se establecía un debate sobre su idoneidad, por el hecho de que constituían núcleos con un marcado carácter elitista que obligaba además a los jóvenes a abandonar sus lugares de procedencia, lo que iba en contra del principio ideológico comunista de "uniformidad de la escuela". Esta discrepancia queda convenientemente expresada por el físico P. L. Kapitsa en una publicación de N. V. Alenksandrov de 1969:

> Actualmente, se han empezado a crear escuelas especializadas para los jóvenes dotados de la Unión Soviética y también de otros países, para prepararles para la labor científica. En el campo de las artes, esto puede estar justificado… Pero que existan escuelas que seleccionen a jóvenes talentosos en matemáticas, física, química y bilogía puede convertirse en algo negativo… Que un joven talentoso, salga de su escuela, equivale a sangrar a la institución, con el consecuente efecto negativo sobre toda la escuela. La razón es debida a que un joven talentoso comparte más su tiempo con el resto de sus compañeros que el propio profesor, y la ayuda y relación entre los compañeros se establece más fácilmente que entre el profesor y los alumnos. Los jóvenes talentosos suelen tener un papel más importante que los profesores en el aprendizaje de los alumnos.
> (Traducción del autor de Shabanowitz, 1978, p. 82)

Para Kolmogorov además, la propuesta experimental a tan pequeña escala de las **escuelas internado especializadas en matemáticas y física** resultaba ser una medida insuficiente, la cual quería ampliar, tal y como le sugiere en una carta A. I. Markushevich (29-12-1964) (Abramov, 2010, p. 107). En esa misma carta proponía que una forma de no vulnerar el principio de “uniformidad de la escuela” y favorecer la educación a partir del desarrollo de los intereses personales de los alumnos, podía ser la introducción de **clases opcionales**, a lo largo de todo el curso escolar, cubriendo varias áreas de conocimiento (dibujo técnico, tecnología de la radio, biología, lengua extrajera, matemáticas o física). Sugería, además, que se podían introducir para los últimos grados, 3 de estas clases para el grado 9 y 6 para el grado 10 (Citado en Abramov, 2010, p. 106). En respuesta A. I. Markushevich parecía estar de acuerdo, aunque le seguían pareciendo prioritarias las **escuelas especializadas en computación** (posteriormente en matemáticas), considerando absurdo el principio de “uniformidad de la escuela” (Abramov, 2010, p. 107).

Durante estos años, a partir de 1964, se había comenzado a fraguar la necesidad de una reforma (la posterior “reforma de Kolmogorov”), por lo que esta idea se incorporó a las ya proporcionadas por Kolmogorov en su propuesta de “modernización” de la enseñanza de las matemáticas. En 1966, el Comité Central del Partido Comunista, aprobó una resolución que permitía incorporar las **clases opcionales** en las escuelas, las cuales no tardaron en implementarse (y lo seguiría haciendo en la década de los 70) gracias a la amplia experiencia en la enseñanza extra-curricular y extra-ordinaria con la que contaban (Abramov, 2010, p. 207). A partir de los 80, y con la oposición de una contra-reforma, se fue deshaciendo gradualmente su implantación, y se cree que empezaron a dedicarse estas horas a la preparación de pruebas de competición (Abramov, 2010, p. 109). A principios de los 90 dejarían definitivamente de existir.

Se ofrecían dos tipos de cursos de matemáticas como **clases opcionales** para los cursos o grados 7 a 10: “Temas complementarios y Problemas de Matemáticas” y varios “Cursos especiales” (Maslova y Markushevitz, 1969, p. 237). El tiempo dedicado a ellos era de 2 h. semanales en los grados 7 y 8, y 6 h. semanales en los grados 9 y 10. El primero de ellos, tenía habitualmente la pretensión de constituir una ampliación de las clases ordinarias. Cada curso tenía una duración de 2 h. a la semana, y los temas eran propuestos por los profesores. En los cursos o grados 7 y 8, se trataban problemas de teoría de números, transformaciones geométricas y el concepto de función. Los cursos o grados 9 a 10, se centraban en teoría de conjuntos, teoria elemental de la probabilidad, métodos numéricos, resolución de ecuaciones, y ampliación del concepto de número; junto a gran parte del tiempo dedicado a la resolución de problemas de dificultad elevada (34 h. de las 70 h. del programa anual). Los llamados “Cursos especiales”, tenían menor relación con las clases ordinarias y normalmente eran elegidos por los alumnos de grados 9 y 10 según sus intereses personales. Se ofrecían 4 cursos de 70 h.: programación informática, matemática numérica, espacios vectoriales y problemas de programación lineal, y elementos de matemática discreta.

## 3. Los Círculos Matemáticos: objetivos y organización de experiencias recientes

Hoy en día se siguen llevando a cabo Círculos Matemáticos en Rusia y no hace mucho Estados Unidos ha recogido la tradición rusa adaptándola según las necesidades específicas de cada lugar a través de la creación de numerosos Círculos en diversas ciudades, constituyendo una práctica importante y ampliamente normalizada. En Europa no existen Círculos Matemáticos, y las experiencias más similares que podemos encontrar, al menos en cuanto al desarrollo de sus sesiones, se conciben bajo la fórmula de **Laboratorios Matemáticos**, dos de las cuales las incluimos aquí debido a su relevancia y repercusión, describiendo lo más característico sobre su funcionamiento.

### 3.1. Los Círculos Matemáticos en Rusia: "Lesser Mekh-Mat" en la Universidad de Moscú (MSU), Círculos de la Escuela nº 58 de Moscú y la Escuela nº 239 de San Petersbugo, y el proyecto Mathematical Etudes.

Actualmente en Rusia, la enseñanza extra-curricular a través de los Círculos Matématicos se continua en diferentes ciudades (Moscú, San Petersburgo, Yaroslavl, Krasnodar, Kirov, Chelyabinsk, Irkutsk, Omsk, etc.), aunque los principales núcleos de referencia siguen siendo Moscú y San Petersburgo. En el resto de ciudades, suelen ser tan sólo unos pocos los profesores involucrados en desarrollar estas actividades, aún así, es posible ver su repercusión en las Olimpiadas, donde se puede encontrar a sus participantes (Marushina y Pratusevich, 2010, p. 400). Cada Círculo Matemático tiene su propia idiosincrasia, y dentro de una misma gran ciudad suelen, hasta cierto punto, competir entre diferentes de ellos tratando de involucrar a cientos de jóvenes. En 2010 San Petersburgo, por ejemplo, contaba con 700 participantes en su red de Círculos, entre los que se encontraban los Círculos de la **Escuela nº 239** de Matemáticas y Física de San Petersburgo, del **proyecto del Palacio de San Petersburgo** para estimular la creatividad de los jóvenes y de la **Escuela matemática para jóvenes** (Marushina y Pratusevich, 2010, p. 400). En Moscú, se encuentran los Círculos del **MCCME** (Centro de Eduación de Matemáticas para todas las edades de Moscú), del **Lesser Mekh-mat** (Escuela de tarde del departamento de Mécanica y Matemáticas de la Universidad de Moscú) y los vinculados a las escuelas que cuentan con cursos avanzados de matemáticas.

En la actualidad la tradición de los Círculos Matemáticos en Moscú, se agrupa bajo el nombre de **Lesser Mekh-Mat**, que recoge la reunión semanal de Círculos para distintos cursos o grados (6 a 11), y la enseñanza a distancia o por correspondencia. El Círculo Matemático se reune en las aulas de la **Facultad de Mecánica y Matemáticas de la MSU** (Mekh-Mat) durante el sábado por lar tarde y mantiene las habituales características de ser gratuito, no realizar exámenes, y de estar abierto a todos los posibles alumnos. Cada Círculo lo conforma un número entre 15 y 30 participantes y entre 3 y 6 instructores (profesores-guía) que suelen ser estudiantes universitarios, salvo el instructor principal, un matemático con experiencia, responsable del diseño del programa del Círculo, aunque a los instructores también les está permitido desarrollar material. Los alumnos de curso o grado 5 pueden acceder al Circulo de curso o grado 6. En total, hay de 100 a 200 participantes en los Círculos de cursos o grados de 6 a 8 y algunos menos en los de 9 a 11 (Dorichenko, 2012, p. xii). El grupo de participantes varía a lo largo del año, aunque se suele mantener un núcleo fijo. La plantilla de instructores involucrados suele mantenerse también fija, salvo ausencias puntuales. Hay instructores que crean su propio curso individual, sobre todo en los cursos más altos, a los que se suelen apuntar normalmente un número reducido de participantes.

La **Escuela nº 58** de Moscú es conocida por haber sido convertida en escuela especializada en matemáticas en los años 70, a partir de que N. N. Konstantinov organizara un curso de matemáticas en ella, alcanzando el éxito y reconocimiento educativo en la segunda mitad de los años 80 (Sossinsky A., 2010, p. 214). En ella se desarrolla un Círculo

Matemático durante los jueves por la tarde, un día importante en matemáticas para la escuela, y con una organización similar al **Lesser Mekh-Mat**. Participan en torno a 100 alumnos por curso o grado (de 6 a 8) repartidos en grupos de 15 y con 3 o 4 instructores encargados de cada uno de ellos, y su organización es similar al descrito anteriormente (Dorichenko, 2012, p. xiv). En este Círculo, algunas sesiones se utilizan para realizar entrevistas a los alumnos con el propósito de seleccionarlos para las clases especiales de matemáticas que se imparten en los cursos o grados 8 y 9, a las que pueden acudir no sólo los alumnos que participan en los Círculos, sino cualquier otro alumno de la escuela.

La **Escuela nº 239** de Matemáticas y Física de San Petersburgo, cuenta con sesiones de Círculos Matemáticos para los cursos o grados 5 a 11, que se realizan dos veces por semana. Las sesiones en este Cículo mantienen un carácter variado, en torno a: exposiciones teóricas, resolución de problemas de forma individual, discusiones con los instructores sobre las soluciones de los problemas, análisis de las soluciones por parte del instructor, entrevistas y pruebas sobre teoría, seminarios, artículos e investigaciones de los alumnos, competiciones matemáticas. Un programa de estas características puede requerir entre 140-150 horas al año, con 80-90 de ellas dedicadas a la preparación de las Olimpiadas (Marushina y Pratusevich, 2010, p. 401).

**Mathematical Etudes** es un proyecto iniciado en 2002 y apoyado por el **Instituto Steklov de Matemáticas** para fomentar y ayudar al estudio de las matemáticas inspirado en la tradición rusa de los Círculos Matemáticos. La única diferencia radica en que, para lograrlo, han elaborado abundante material virtual, consistente en vídeos, animaciones, maquetas y apps, cuyo objetivo es favorecer la comprensión de los objetos matemáticos a través de la visualización y la interactuación (Andreev, Dolbilin, Konovalov, y Panyunin, 2014, p. 48). El proyecto también comprende la elaboración de una bibiloteca virtual a partir de la digitalización de revistas de contenido matemático, y un subproyecto basado en explicar el funcionamiento de mecanismos diseñados por Chebyshev, matemático del s. XIX, también modelados utilizando computación informática. Todo el material se puede encontrar en la URL: http://www.etudes.ru/ completamente disponible en ruso y parcialmente en inglés, francés e italiano. El material ha sido empleado, por Nikolai Andreev, uno de los desarrolladores del proyecto, en más de 500 conferencias dirigidas a alumnos y profesores de primaria y sencundaria así como a estudiantes universitarios de diferentes regiones de Rusia. Además, ha sido y es utilizado actualmente por profesores de matemáticas para sus clases y en actividades extra-curriculares como los propios Círculos Matemáticos (Andreev et al., 2014, p. 52).

### 3.2. Los Círculos Matemáticos en Estados Unidos: las experiencias de Boston y Berkeley.

Los Estados Unidos, han comenzado a desarrollar estos programas en la década de los 90, gracias a la herencia cultural de inmigrantes asentados en el país (Vandervelde, 2009, p. 6). Estas experiencias en Estados Unidos se recogen en la **Asociación Nacional de Circulos Matemáticos**, cuya URL es https://www.mathcircles.org/, donde actulamente existen 169 Círculos Matemáticos activos. La asociación también recoge otros 5 Círculos Matemáticos en Canada, 2 en China, 1 en Rusia (el acogido por el MCCME), y 1 online bajo el nombre de **Camp Euclid**. Podemos por tanto considerar a Estados Unidos, el principal país heredero de la tradición rusa de los Círculos Matemáticos para adaptarlos a sus propias circunstancias, aunque hay que tener en cuenta que cada uno de estos Círculos tiene sus propias características y pueden desarrollarse de maneras diferentes, lo cual podría ser objeto de una futura investigación. Además cuentan con su propio programa para la formación de profesores bajo el nombre de **Math Teacher's Circle Network**, inciado en 2006 bajo el apoyo del **AIM (American Institute of Mathematics)**, consistente en que los profesores trabajen en equipo, en Círculos de 15 a 20 participantes, en la resolución de problemas y la discusión entre ellos. Cuenta en la actualidad con unos 300 Círculos y publican un boletín informativo semestralmente. La URL de esta iniciativa es:

http://www.mathteacherscircle.org/ Un video, a modo de resumen, sobre estos programas, se pude visionar en la URL: https://www.youtube.com/watch?v=ObT8F4mZbEg

La primera de estas experiencias fue el **Círculo Matemático de Boston** fundado por Robert Kaplan y Ellen Kaplan con la ayuda de Tomás Guillermo en Boston (Masachussets) en Septiembre de 1994 ante la necesidad de encontrar “una manera más apropiada de enseñar matemáticas” (Traducción del autor de Kennedy, 2003, p. 27). Perseguían la idea de establecer un método mayeútico en su aprendizaje inspirados por el **método de Moore**, empleado a principios del s. XX por el matemático Robert Lee Moore, que consistía en enseñar a un grupo de estudiantes sin ningún conocimiento sobre el tema a tratar, a través de que realizaran sus propias demostraciones y resolución de problemas sin consultar ningún otro material, nada más que su propia asistencia a las clases y un trabajo de investigación constante. Los Kaplan habían aprendido sobre la existencia de los Círculos Matemáticos rusos a través del contacto con la comunidad de emigrantes rusos en Boston, por lo que finalmente y en consecución de su inciativa, decideron dejar sus trabajos como profesores de escuela secundaria y fundar el primer Círculo Matemático en Estados Unidos, alquilando para ello un espacio, los sábados por la mañana, en una iglesia local. Contaron en su incio (septiembre de 1994) con los primeros 29 participantes de diferentes edades, con los que formaron un único grupo. En el siguiente semestre, cuando ya contaban con 34 alumnos, el programa fue definitivamente acogido por las **Universidades de Northeastern** y **Harvard**, lo que les permitió incorporar a un matemático graduado por Harvard y dedicarse a la difusión y ampliación del Círculo. En su segundo año, con 38 estudiantes, impartieron las sesiones en domingo con una duración de 2 h. Se siguió incrementando la participación, y en 2001 contaban con 205 estudiantes, (The Math Circle) por lo que se vieron obligados a constituir 9 sesiones entre semana, estableciendo diferentes grupos para los grados 1 a 6, además de la sesión del domingo con 3 niveles: inicial (grados 4-6), medio (grados 7-9) y superior (grados 10-11). A partir de entonces diferentes matemáticos graduados, del **MIT** (Instituto tecnológico de Massachusetts) y de la **Universidad de Tufts** se han sumando ocasionalmente impartiendo cursos.

Las sesiones se establecen basándose en diferentes temas y manteniendo generalmente la conexión con el curriculo obligatorio de las escuelas, aunque focalizando en que los estudiantes formulen sus propias conjeturas y contraejemplos, dejando a un lado la necesidad de cubrir todo el temario propuesto para el desarrollo de las sesiones (The Math Circle).

Después, en 1998, otra experiencia que ha conseguido un gran éxito es la del **Círculo Matemático de Berkeley**, fundado por iniciativa de la actual directora Zvezdelina Stankova, matemática de origen húngaro, llegada a Estados Unidos en 1989 para estudiar matemáticas consiguiendo su doctorado por la **Universidad de Harvard** en 1997. En su adolescencia participó en los Círculos Matemáticos en Hungría, además de ser medallista en diferentes Olimpiadas Matemáticas, En 2011 fue premiada por la **Sociedad Estadounidense de Matemática (ASM)** con el **Premio Haimo** a la mejor profesora de matemáticas a nivel nacional, ya que también es profesora en el **Mills College** (Okland). (MAA Mathematical Association of America, 2015) Hasta entonces sólo existía el Circulo de Matemáticas de Boston, así que decidió fundar un nuevo Círculo en Berkeley, para que sirviera de ejemplo a otros profesores para crear sus propios Círculos en las escuelas. Finalmente, acabaría convirtiéndose en una experiencia importante y de largo alcance ya que inspiraría en la creación de más de 100 futuros Círculos (Weld, 2014). El Círculo debe su apoyo a la **Universidad de California** (Berkeley), el **MSRI** (Mathematical Sciences Research Institute) y los propios padres de los participantes. El programa, financiado a través de varias instituciones, intenta mantener el carácter económicamente accesible para todos ofreciendo un coste reducido de matrícula anual. Actualmente cuenta con unos 400 alumnos. El Círculo, fiel a sus orígenes, ofrece la posibilidad de entrar en contacto con las matemáticas, de forma más cercana a los propios procesos de investigación matemática, a través de matemáticos profesionales o graduados, que son los encargados de organizar y dirigir las sesiones. Además, prestan gran atención a los procesos colaborativos mediante la exploración matemática con los compañeros, así como a la exposición escrita y comunicación verbal matemática. Se puede ver una presentación sobre el Círculo en el siguiente enlace:
http://www.globalpres.com/mediasite/Viewer/?peid=d503dcf468b0475b9265bf2244ed7aac/

El objetivo de su programa se puede resumir en "aumentar, el número y la calidad de estudiantes que en un futuro se puedan convertir en matemáticos, dedicados a la educación o a la investigación, o de quienes simplemente gusten y empleen las matemáticas en sus estudios, trabajo y actividad diaria" (Traducción del autor de Berkeley Math Circle, 1998-2015). El programa se lleva a cabo a través de sesiones semanales y de concursos mensuales, y se desarrolla en torno a temas de matemática avanzada a través de la resolución de problemas a la vez que potencian sus habilidades de escritura en la exposición de sus soluciones. Se inició en 1998 con 50 participantes entre los cursos o grados 7 a 12, con una única sesión semanal de 2 h. El número de participantes se fue incrementando y la organización del Círculo se amplió a diferentes grupos para todos los niveles, hasta incluir en 2009 por primera vez un grupo para nivel elemental. A partir de entonces el Círculo se organiza en torno a la división en estos dos grandes niveles que a su vez se subdividen en otros.

El nivel elemental (**BMC-Elementary**), dirigido actualmente por Laura Givental, cuenta con los niveles elementales I y II (grados 1-2 y 3-4). La creación de este nivel constituyó una propuesta innovadora ya que los Círculos suelen formarse, normalmente a partir del grado 5, a lo sumo a partir de grado 4. Al principio se inció con un grupo de 12 a 15 estudiantes de grado 2, durante el primer semestre del curso 2009-2010, a cargo de la profesora Natalia Rozhkovskaya*,* al que pronto se sumaron más participantes hasta un total de 40 de grados 1-3, por lo que se incorporaron dos nuevos profesores (Rozhkovskaya, 2014, p. xii). El proyecto continuó así por 2 años más, hasta que se vio la necesidad de conformar dos niveles para incluir a los estudiantes de grado 4 en el primer semestre del curso 2011-2012. El objetivo de las sesiones a nivel elemental (BMC-Elementary) se dirige hacia "el fomento de una actitud positiva de los niños hacia las matemáticas introduciéndoles elementos de la cultura matemática" (Traducción del autor de Berkely Math Circle, 1998-2015), lo que se consigue fundamentalmente haciendo que la matemáticas sean accesibles y divertidas a través del juego (Rozhkovskaya, 2014, p. xv).

El nivel superior (**BMC-Upper**), está dirigido por Zvezdelina Stankova, y cuenta con los niveles intermedios I y II (grados 5-6 y 7-8) y un nivel avanzado (grados 9-12). Las sesiones se organizan en torno a temas al cargo de algún matemático experto, algunos de ellos son: combinatoria, teoría de grafos, álgebra lineal, transformaciones geométricas, recursividad, series, teoría de conjuntos, teoría de grupos, teoría de números, curvas elípticas, geometría algebraica y aplicaciones en informática, ciencias naturales y economía, etc. Los **concursos mensuales** consisten en una serie de problemas que hay que resolver individualmente, y a los que solo pueden presentarse los alumnos matriculados en el Círculo. Para ello, se proporciona un tiempo de 3 o 4 semanas, y se puede consultar cualquier libro o material, pero debe resolverse de forma individual, es decir, sin que nadie les haya proporcionado ninguna ayuda. Siguen la misma estructura que las **Olimpiadas Matemáticas del Área de la Bahía de San Francisco**, con 7 problemas, de los cuales, los problemas 1-4 son para los grados 4 a 8 y los problemas 3-7 corresponden a los grados 9 a 12. Cualquier alumno puede optar a presentarse también a un nivel superior, pero no al contrario. Entre los propósitos de la competición, se encuentra el de mejorar la calidad de la exposición escrita en la resolución de problemas, además de que los problemas seleccionados suelen dar lugar a la comprensión y desarrollo de una teoría importante, por tanto, las soluciones se evalúan desde ese punto de vista.

Adicionalmente, se imparte semanalmente en domingo, una conferencia a cargo de profesores universitarios sobre algún tema con una duración de 2 horas en el campus de la **Universidad de California** (Berkeley), dirigido a los estudiantes de secundaria (grados 6-

12) del Área de la Bahía de San Francisco. Durante su desarrollo suelen establecerse discusiones con los estudiantes, quienes pueden exponer sus problemas y soluciones. Algunas de estas conferencias se organizan también en torno a la resolución de problemas de competiciones matemáticas, o en temas aplicados de las matemáticas en relación con otras ciencias.

Otra experiencia importante que también es directa heredera de la tradición rusa en Estados Unidos es la Escuela de Matemáticas de Rusia (**RSM**, **Rusian School of Mathematics**) fundada en 1997 por Inessa Rifkin e Irina Khavinson, que cuenta con un importante programa de matemáticas extra-escolar para todos los niveles, cuyo “currículo y metodología están basados en la experiencia personal en escuelas especializadas de matemáticas en Rusia” (Traducción del autor de RSM, 1997-2015).

### 3.3. Algunos ejemplos de Laboratorios Matemáticos en Europa: MateMatita en Italia y ESTALMAT en España.

**MateMatita** es un Centro Interuniversitario de investigación para la comunicación y el aprendizaje informal de las matemáticas, fundado en 2005, en el que participan las Universidades de Milán, Milán-Bicocca, Pisa y Trento. Su objetivo principal es averiguar cuales son los contenidos y las metodologías más idoneas dentro del aprendizaje informal. Están interesados en el papel que juega el lenguaje, las posibilidades que ofrecen las tecnologías de la información y la comunicación (TIC), la relación entre las matemáticas y las disciplinas artísticas u otras ciencias, y la relación entre las matemáticas aplicadas y la tecnología. Para poder estudiarlo, elaboran diferentes tipos de materiales (exposiciones, libros, revistas  y material multimedia), y observan su impacto en contextos educativos y en el público en general. (MateMatita) Debe su nombre a la unión de “Mate” (matemática) y “Matita” (lápiz), y su intención es también la de encontrar fórmulas que hagan las matemáticas más atractivas y accesibles, mediante el empleo de todos estos materiales y a través del uso de metodologías no convencionales.

Debido a ello la parte más ambiciosa  del proyecto ha consistido en emplear estos materiales en las escuelas a través de **laboratorios matemáticos** en vez de las clases tradicionales donde el profesor explica y luego los alumnos realizan la tarea que se les propone en referencia a lo expuesto. Un laboratorio matemático se caracteriza entre otras cosas, porque los alumnos tienen un papel activo y dirigido hacia un objetivo, mientras que los profesores tienen únicamente un papel de guía que consiste en observar, escuchar y responder a las preguntas que surgen. Se suelen utilizar materiales manipulables y también es frecuente el uso de la tecnología (Dedò y Di Sieno, 2013, p. 324). Para poder desarrollar estos laboratorios se requiere de toda una planificación que comienza con la preparación adecuada de los profesores que harán de guías en el proceso de enseñanza-aprendizaje. Los laboratorios matemáticos se desarrollan utilizando el método científico, donde la observación, el establecimiento de planteamientos o conjeturas, el adecuado uso del rigor y el error, la discusión, la adopción de conclusiones y el repaso de los puntos fundametales en toda la elaboración son especialmente tenidos en cuenta. Es muy importante para la investigación, la correcta evaluación del desarrollo de estos laboratorios, lo que hace necesario la permanente discusión entre los profesores involucrados adoptando nuevas posturas que la favorezcan (Dedò y Di Sieno, 2013, p. 340). En la siguiente URL: http://www.matematita.it/realizzazioni/materiale_didattico.php/ se pueden encontrar algunos de los kit de laboratorio utilizados en las escuelas, y que se han empleado en 500 laboratorios involucrando a más de 12000 estudiantes y 1000 profesores. Otros laboratorios, más de 1000, se han lllevado a cabo en el **Departamento de Matemáticas de la Universidad de Milan**, involucrando a unos 25000 estudiantes y 2000 profesores (Dedò y Di Sieno, 2013, p. 341). Dentro de los materiales elaborados, los que han tenido una gran repercusión en las aulas de primaria y secundaria italianas son los juegos on-line que se pueden encontrar en la URL:  http://www.quadernoaquadretti.it/ La opinión de los numerosos profesores que han participado en estas actividades ha sido muy positiva, ya que según ellos los estudiantes que han contado con esta oportunidad han tenido muchas menos dificultades durante los exámenes nacionales italianos que los compañeros que no (Dedò y Di Sieno, 2013, p. 341).

El proyecto **ESTALMAT** en España, es una iniciativa del profesor Miguel de Guzmán Ozámiz apoyada por la **Real Academia de Ciencias**, pensada para detectar, estimular y orientar el talento matemático de los jóvenes españoles. El nombre del proyecto es un acrónimo obtenido de las primeras letras de “estímulo del talento matemático” y fue iniciado en 1998 en la **Comunidad Autónoma de Madrid**, aunque actualmente también está presente en otras Comunidades. El proyecto, en su momento innovador y pionero (Citado en ICMAT, 2013), se fundó con la intención de detectar y guiar, a jóvenes con especiales habilidades matemáticas, sin desarraigarlos de su entorno, y con el objetivo de que la sociedad se beneficiase de ello, ya que “la sociedad que facilita el desarrollo de sus individuos, avanza indudablemente más que la que permanece indiferente” (Traducción del autor de La Real Academia de Ciencias Exactas, Físicas y Naturales, Spain, 2002-2003, p.120). La participación en el proyecto no tiene ningún coste para los alumnos seleccionados ya que está financiado por diversas instituciones públicas y empresas privadas.

El proyecto selecciona anualmente un máximo de 25 jóvenes con edades de 12 y 13 años, en cada una de las ciudades donde está presente, presentándose en total, en España durante 2013, unos 2500 alumnos (ICMAT, 2013), los cuales para poder hacerlo tienen que haber sido recomendados por sus profesores de matemáticas de las escuelas. La selección se hace a través de una prueba (de 4 a 6 problemas originales) que intenta medir las habilidades matemáticas, más que los conocimientos que hayan podido adquirir, con la que se hace una primera selección, y posteriormente mediante entrevistas individuales (también a los padres) para ver el compromiso y disposición que tienen con el programa. Los seleccionados participarán durante 2 años en un programa de 1 sesión semanal de 3h., y posteriormente en 1 sesión mensual hasta finalizar el bachillerato. El programa es completamente independiente del de los centros escolares. Miguel de Guzmán organizó también de forma paralela, un curso de 75 h. en el primer cuatrimestre del año académico 1999-2000, para los alumnos de primer año titulado ***“Laboratorio de Matemáticas”*** dedicado al trabajo práctico y autónomo, bajo la guía de un profesor, en la resolución de problemas, (Madrid) y que pretendía “ayudar a los estudiantes de primer curso a hacerse con procedimientos prácticos básicos para afrontar las dificultades de adaptación al estudio de las matemáticas en la Universidad” (Madrid).

El programa se desarrolla bajo un formato similar al de los Círculos Matemáticos. Con los alumnos seleccionados se forman grupos de 5, cada uno de ellos conducido por 2 instructores, que son además, profesores universitarios o profesores de secundaria especialmente competentes. Las 3 h. de la sesión se dividen en dos partes de 1 h. y 20 min. cada una de ellas, con un descanso intermedio de 20 min. A cada parte de la sesión se le asigna un tema del que es responsable uno de los dos instructores, de tal forma que todo el programa anual está perfectamente planificado. En las sesiones se tratan diferentes temas como: pre-álgebra, grafos, algoritmos, teoría de números, juegos de estrategia, geometría interactiva, poliedros, mosaicos, fractales, etc. En todos los casos, las sesiones se desarrollan principalmente a través de la resolución de problemas donde el texto ***“Mathematical Circles (Russian Experience)”*** constituye una fuente principal de ellos. (Hernandez, 2006, p. 5) La bibliografía principal y parte del material utilizado puede ser consultado en la URL: http://www.uam.es/personal_pdi/ciencias/ehernan/Talento/Indice.htm/

Las sesiones constituyen **laboratorios matemáticos** donde “se intenta proporcionar un lugar donde poder experimentar y tener la libertad para jugar con conceptos, así como discutir sobre las ideas que surgen gracias a la habilidad de los alumnos en ofrecer respuestas a las preguntas que se les plantean” (Traducción del autor de Hernandez, 2006, p. 7). En ellas se trabaja de forma colectiva e individual. En el trabajo en grupos, después de la presentación de un tema por el instructor, los alumnos tienen que buscar ideas y soluciones a los problemas que se les plantean, para lo que se hace necesaria la discusión entre ellos. En el trabajo individual, se suele exponer un problema de cierta complejidad para su posterior solución. En ocasiones, no se llega a ninguna, y aunque no hay tarea para casa, las cuestiones quedan abiertas para quien libremente quiera seguir dedicándose a ellas. Entre los objetivos del programa, a diferencia de los Círculos, no se incluyen la preparación para competiciones, lo que no impide que los alumnos lo hagan por su cuenta y participen.

Sobre la evaluación de los resultados ofrecidos hasta ahora por el programa, se puede afirmar que la tendencia mayoritaria no es la de estudiar grado de matemáticas y dedicarse a la investigación, en algunos casos ni siquiera estudian grados de ciencia (Citado en ICMAT, 2013). Una amplia mayoría termina cursando grados en diferentes tipos de ingeniería (Hernandez, 2006, p. 8). En cualquier caso, los participantes de este programa están de acuerdo en que les ha proporcionado “una visión profunda sobre lo que significa trabajar con matemáticas” (Traducción del autor del texto citado en La Real Academia de Ciencias Exactas, Físicas y Naturales, Spain, 2002-2003, p. 122). Por el momento, parece difícil conocer el impacto (cultural, científico y técnologico) a largo plazo del programa (Hernandez, 2006, p. 8).

## 4. Los Círculos Matemáticos: cómo establecer el proceso didáctico (una cuestión práctica)

Un **Círculo Matemático** es una experiencia educativa informal, ya que no es obligatoria su asistencia y no se realizan exámenes ni se otorgan títulos, pero al mismo tiempo es una experiencia educativa extremadamente seria, pues constituye una práctica cuidadosamente pensada y planificada para la consecución de un aprendizaje que tiene una repercusión primordial en la formación futura de sus participantes. Otra característica importante es que se trata de una actividad gratuita y normalmente impulsada por personas con una gran admiración y dedicación hacia las matemáticas, que deciden invertir su tiempo, en su orígenes de forma altruista, para formalizarlos. Es de especial importancia en su carácter de red informal, la vinculación con ellos de importantes matemáticos dedicados a la investigación a través de la Universidad, lo que les proporciona una imagen profesional, atractiva y confiable. Por tanto, la creación de un Círculo Matemático implica esfuerzos importantes, a través de los cuales, se formaliza también una práctica educativa que no es cubierta por la educación obligatoria u ordinaria, y de la que se desprenden resultados importantes, como fuera el caso de la formación inicial de abundantes matemáticos y científicos durante su época con mayor desarrollo en la Unión Soviética. La experiencia se extiende, tal y como hemos visto, hasta la actualidad y se vislumbra un importante futuro para ella, siempre y cuando se cuente con una comunidad de profesores y matemáticos dispuesta a impulsarla.

En Noviembre de 2006 en un viaje organizado, algunos integrantes del **Círculo Matemático de San José**, tuvieron la oportunidad de ver y en ocasiones también de participar en algunos Círculos Matemáticos de San Petersburgo (durante 3 días) y Moscú (durante 1 día). La profesora Tatiana Shubin, perteneciente al Departamento de Matemáticas de la Universidad de San José y directora del Círculo Matematico de San José, describe la experiencia:

> La organización era bastante diferente a lo que estábamos acostumbrados en nuestros propios Círculos. No se trataba de una clase, no había nigún profesor en frente de la clase dirigiendo las actividades. En vez de ello, se repartían a los alumnos hojas con problemas para que trabajasen en ellos de forma individual, y cuando un alumno quería hablar de su solución, un profesor se sentaba junto a él y le escuchaba. Con el tiempo el silencio de la clase se fue llenando con animadas conversaciones, y cada alumno tenía un matemático a quién dirigirse y para que le ayudara en su solución con los problemas, o mejor aún, descrubriendo otras cuestiones a raíz del problema o dilucidando implicaciones importantes.
> (Traducción del autor de Shubin, 2012, p. 215)

A partir de entonces decidieron organizar dos sesiones de Círculos al “estilo ruso” anualmente, lo que conllevaba la participación de muchos voluntarios como profesores asistentes. En Marzo de 2007, tuvieron el privilegio de contar con la participación del matemático ruso Vladimir I. Arnold en una de estas sesiones. La forma de trabajar consistió, en un principio, en proponer una semana antes una tanda de problemas, para que posteriormente, los alumnos, pudieran terminar de solucionarlos y comentarlos durante la sesión en el Círculo ayundandoles en la búsqueda de la mejor opción para su solución y en como exponerla. Esta tanda previa se daba tan sólo a modo de entrenamiento y, a continuación, se hacía entrega de una nueva tanda durante el inicio de la sesión para que la solucionaran durante los primeros 40 min. Las sesiones se desarrollaban durante 2 h. y los participanes eran de los cursos o grados 7 y 8 (Shubin, 2012, pp. 216-217).

Con lo anterior queremos enfatizar la idea de que un Círculo Matemático puede organizarse de una manera flexible, adaptándose a los recursos y posilidades existentes, pero siembre tienen como objetivo primordial el hecho de que se desarrollen en torno al aprendizaje mediante la resolución de problemas que no requieran necesariamente, al menos en sus etapas iniciales, de profundos conocimientos matemáticos. Lo que se pretende es, desde los primeros cursos o grados (5 y 6 o incluso antes) desarrollar estas habilidades tan importantes en matemáticas, que son esencialmente las mismas que las que utiliza un futuro matemático en sus procesos de investigación, pues como dice el matemático español Lluis A. Santaló: “en la historia de las matemáticas, la curiosidad por la resolución de problemas de ingenio ha sido un factor que ha contribuido a la creación matemática tanto más que sus posibles aplicaciones prácticas” (Traducción del autor del texto citado en Fomin et al., 2012, p. 5).

### 4.1. Consideraciones para la formalización de un Círculo Matemático.

En un Círculo Matemático la asistencia es totalmente libre, por lo que a la hora de atraer posibles participantes, es bueno tener en cuenta cuál es su motivación y más aún, que es lo que hará que se vuelvan asiduos al mismo. Lo recomendable es que el Círculo permanezca abierto, y que, aunque haya una cierta variación de participantes, se pueda mantener un núcleo fijo. Por lo general, muchos de ellos, bien por su propia iniciativa o por la de sus padres, acuden para obtener una respuesta a la pregunta de si están verdaderamente interesados en las matemáticas y además se les dan bien, por lo que suele contarse en su mayor parte con un público afín (Vandervelde, 2009, p. 10). El Círculo Matemático deberá ofrecer una respuesta que sea capaz de enganchar a los estudiantes que acudan y haga que su asistencia se mantenga. Dado el carácter extra-curricular y por tanto de no obligatoriedad, se hace necesario que la actividad resulte amena e incluso divertida, por lo que una de las intenciones principales, además de desarrollar el pensamiento matemático, es el de mantener una comunicación viva y espontánea entre los estudiantes y profesores, así como tratar de fomentar la participación de todos sin miedo a equivocarse todas las veces que sea necesario (Fomin, Genkin, y Itenberg, 1996, p. 8). El profesor o instructor deberá ser un ejemplo continuo demostrando esta actitud, por lo que se hace imprescindible el estar bien preparado para poder liderar un Círculo Matemático. Dentro de las características de un Círculo Matemático, está la de su viculación con el departamento de Matemáticas de alguna Universidad, siendo impulsado por matemáticos importantes que desarrollen su contenido y programas de referencia que hagan de él algo característico además de asegurar un elevada calidad en su enseñanza. De ser así, se establece también que la comunidad universitaria interesada, alumnos y graduados, pueda participar y por tanto, dotar siempre al Círculo con instructores y asistentes suficientes según la asistencia y las actividades que se quieran desarrollar. En ellos también es bienvenida y deseable la presencia de profesores de matemáticas de secundaria, pues su experiencia e interés son especialmente valiosos. El carácter altruista, de todas las personas que hacen posible un Círculo Matemático es siempre el más deseable, y aunque es entendible que esto no pueda ser siempre así, hay que intentar que en la medida de lo posible, lo que lo impulse sea precisamente este gran interés hacia las matemáticas y un carácter de libre participación en él, pues la meta consiste precisamente en organizar una comunidad matemática, vículada a una ciudad o región, que ofrezca posibilidades más ambiciosas a los niños y adolescentes que estén interesados en descubrir el apasionante mundo matemático.

La organización más deseable es la que se realiza por niveles según cursos o grados, pero si esto no es posible, al menos en sus inicios, entonces habrá que establecer grupos más amplios, abarcando eso sí, los cursos más próximos. Lo importante será establecer un grupo cuya dimensión no sea ni demasiado pequeña ni demasiado grande, un número entre 5 y 15 participantes se considera el más deseable, ya que permite trabajar de forma individual pero también colectiva. Será importante, la elección de las instalaciones dónde desarrollar el Círculo. Una opción es, si se cuenta con el apoyo de la Universidad hacerlo en su Campus, pues las posibilidades de contar ampliamente con diferentes aulas y salones de actos, facilitará mucho la organización del Círculo, además de que la Universidad constituye un lugar de referencia con valor añadido que hará que mucha gente pueda sentirse atraída a participar, proporcionando un primer acercamiento a ella a todos los alumnos de escuelas e institutos que pueden verlo como algo lejano e inaccesible. Los

Círculos también pueden organizarse a nivel de la propia escuela, o de varias de ellas, pero lo que sí que es conveniente es que la Universidad esté implicada en el desarrollo del mismo, pues se trata de ofrecer una programación muy cuidada que estimule el pensamiento matemático, tal y como lo realiza un matemático profesional dedicado a la investigación, por encima de la adquisición de contenidos matemáticos determinados, aunque estos también se consiguen pues son multitud las áreas y temas matemáticos que se pueden ofrecer como programa, y que no tienen porque ser obligatoriamente una extensión de aquellos que se estudian en la educación ordinaria, sino que serán elegidos según los intereses que se crean más convenientes. Todo ello, son decisiones importantes y características en la formalización del Círculo Matemático. La elección del horario, es otra de ellas, que se podrá establecer según el área geográfica a la que el Círculo quiera atender así como en función de otras actividades que se puedan desarrollar de forma extra-escolar y contra las que el Círculo tendrá inevitablemente que competir. Lo que sí que es deseable es que sea un día o dos a la semana en horario fijo con la intención de constituir una rutina. La carga horaria de la sesión, es aconsejable que se realice durante 1 hora de forma continuada, o utlizando 2 horas con un descanso en medio, pues con ello se favorece la atención y concentración, además de ser tiempo suficiente para que los participantes indaguen en el tema propuesto. Lo más importante será la organización de la sesión, de la que es responsable el instructor principal o profesor guía, ya que la forma en que se realiza es muy diferente a la de una clase ordinaria, y es posible que exista la tendencia a que se asemeje a una de ellas por cuestiones de herencia cultural en nuestra propia formación, por lo que tan importante como presentar las matemáticas de “otra manera” mucho más atractiva, es que la sesión se conduzca bajo la forma de un Círculo Matemático (Vandervelde, 2009, p. 17). Es por ello, que el mejor entrenamiento, si no existe una formación específica para liderar uno, será poder asistir a algunos de ellos para tener una idea certera sobre cuales son las posibilidades para que su conducción sea la idónea.

Para el establecimiento de los sesiones, se considera erroneo desarrollarlas sobre un único tema para los participantes más jóvenes, y en su lugar, es preferible que haya cambios suficientes que la dinamicen. También es importante retornar sobre temas pasados, que se discutan suficientemente las soluciones y demostraciones, y que se realicen sesiones lúdicas y juegos matemáticos (Fomin et al., 1996, p. 8). Lo más característico de un Círculo Matemático es que las sesiones se organicen en torno a temas, al menos para los grupos a partir de un grado medio, en cualquier caso, tal y como dice Mark Saul, el objetivo de un Círculo Matemático nunca debería ser el de cubrir un temario determinado sino que se utilicen determinados temas por su propio interés (Citado en Vandervelde, 2009, p. 68). Tom Davis, propone la siguiente regla en referencia a como se ha de dirigir una sesión: “puntuate a ti mismo, 1 punto por cada minuto que estés hablando, 5 puntos por cada minuto que un estudiante hable, 10 puntos por cada minuto que un estudiante discuta con otro, y 50 puntos por cada minuto que los estudiantes discutan entre ellos” (Traducción del autor de Davis, 2014), todo ello con el objetivo de que sean los propios estudiantes quienes se interesen por los temas y realicen sus pequeñas investigaciones que tendrán que demostrar y debatir con los otros, tal y como lo haría un grupo de investigadores. Por tanto el objetivo de una sesión es bien distinto al de una clase magistral, donde el interés, la iniciativa propia y la interactuación son muy reducidas, primando el cumplimiento del temario y la forma que para ello propone el profesor. Paul Zeitz, propone conducir la sesión de forma que se presente de forma atractiva un tema, se proponga una serie de problemas que resulten atractivos para los participantes, se trabaje como grupo para llegar a un entendimiento sobre los mismos, y por último y como algo especialmente deseable de que acontezca, que la sesión desemboque un resultado satisfactorio o conclusión.

Es fundamental, para la configuración de una sesión, la elección de los problemas que se van a tratar, pudiendo abarcar diferentes niveles de dificultad. El instructor o profesor-guía será el encargado de elegirlos, así como el tema a tratar, y se establecen en torno a su criterio personal en base a los objetivos que quiera cubrir en la sesión. Durante la sesión el instructor y asistentes trabajarán con los estudiantes, escuchando sus propuestas, dándoles pistas y ofreciéndoles el mejor diálogo para que ellos mismos encuentren su propio camino en lo que tratan de conseguir. Hay que evitar en cualquier caso ofrecerles la solución de forma inmediata, si el problema es interesante, es preferible dejar el tiempo necesario para que se piense sobre él, no tiene porqué resolverse inmediatamente en la sesión en la que se ha propuesto, sino que puede abarcar varias sesiones. El instructor o profesor guía, puede hacer las introducciones o aclaraciones que sean necesarias, pero es importante que no consuman mucho tiempo, y no supongan un desarrollo teórico complejo, ya que todas las formulaciones es preferible dejarlas para el final, cuando un tema haya sido lo suficientemente investigado, de forma que resulte parte natural del proceso.

La forma de promover la participación en el Círculo Matemático se hará a través de su difusión por medios publicitarios o directamente por invitación dándolo a conocer a través de los profesores de matemáticas de las escuelas, quienes pueden a su vez recomendárselo a determinados alumnos. Será importante organizar una adecuada presentación del mismo, si es posible coincidiendo con el inicio del año lectivo, y en un auditorio donde se pueda llevar una presentación diseñada para la ocasión. Este acto será muy importante ya que se trata de atraer la atención de los posibles participantes, sobre todo la de aquellos que sean los más idóneos, que no necesariamente han de ser lo que hayan obtenido mejores resultados en la asignatura obligatoria de matemáticas, sino únicamente aquellos que se sientan especialmente atraídos hacia ellas (Vandervelde, 2009, p. 22). Debido a ello sería muy apropiado mostrar como se va a desarrollar el Círculo Matemático preparando alguna actividad significativa. Si es posible, es muy conveniente realizar presentaciones previas locales en las escuelas e institutos, a cargo de los futuros instructores y asistentes de los Círculos Matemáticos, que muestren cómo funcionan y distribuyan la información más relevante sobre el mismo. También es de gran utilidad contar con aquellos padres que se muestren más decididos a apoyar la actividad, ya que pueden ser un buen canal de comunicación de sus relaciones sociales, y de hecho, su implicación se ha demostrado de gran importancia para el desarrollo y mantenimiento de algún que otro Círculo Matemático (Vandervelde, 2009, pp. 23-34). La configuración de una página de internet que se mantenga actualizada de tal forma que sea representantiva de la actividad del Círculo, se considera una opción indispensable hoy en día. Es importante que el diseño de la página sea atractivo y esté muy bien organizado su contenido, ya que los tipos de información que puede contener seguramente serán variados y amplios (Vandervelde, 2009, p. 29). Por ejemplo, colgando tan sólo todo el material de las actividades para las sesiones, conferencias, competiciones, etc., que suele ser el material más característico de un Círculo Matemático, ya se requiere de un espacio y organización importante. Dotar a la página de material adicional a modo de biblioteca virtual, puede resultar una opción muy apropiada, ya que hay que tener en cuenta que parte de la finalidad de un Círculo Matemático consiste precisamente en el fomento y divulgación de las matemáticas, por lo tanto, atraer la mayor cantidad de visitas diarias será un dato muy significativo. Otra función que puede albergar la página es el registro de todos los participantes, de tal forma que tengan acceso a información personal y privada como la corrección de los problemas, foros temáticos, tablón de anuncios, etc., ofreciendo además la posibilidad de ser notificados a través de sus correos electrónicos.

Dado el carácter extra-curricular y de libre asistencia es muy importante cuidar el aspecto social, facilitando en un principio que todo el mundo se conozca y se sienta agusto, de esta forma los participantes estarán más predispuestos a continuar. Para que así sea, será importante contar con un coordinador en cada uno de los grupos, de modo que el grupo lo tenga como referencia y ayude en la presentación de los diferentes instructores o asistentes. En las primeras sesiones del Círculo Matemático se pueden proponer problemas de diferente difícultad para resolver en grupo, para posteriormente debatir las soluciones, o también se pueden incluir juegos matemáticos compitiendo entre pequeños equipos. Adicionalmente, puede dedicarse alguna sesión con algún carácter más festivo, incluyendo también por ejemplo una conferencia breve y divertida que reuna a todos los integrantes del

Círculo Matemático. Repetir estos encuentros especiales a lo largo del año, será siempre bien recibido y ayudará a la cohesión.

La organización de las conferencias puede ir dirigida a varios grupos y niveles, con una periodicidad mensual, invitando a diferentes matemáticos dedicados a la investigación. Es una baza fundamental para que todos los participantes en los Círculos tengan un contacto directo con matemáticos que están desarrollando su propia investigación en un campo matemático determinado. Diseñar este programa a lo largo del año, constituye una opción muy atractiva, a la que además pueden sumarse un determinado número de invitados que no estén dentro del Círculo, de modo que favorezca su difusión y conocimiento a futuros participantes. Estas sesiones, ayudan a dotar al Círculo de un cierto ritmo, rompiendo con las sesiones en pequeños grupos semanales, además de facilitar la integración y socialización entre todos ellos.

Otra baza con la que pueden contar los Círculos Matemático, son las competiciones, bien estableciendo el propio Círculo la suya, o bien dedicando partes de las sesiones a la preparación de una de ellas, tal y como se hacía en la tradición rusa con las Olimpiadas. El hecho de dedicar parte de las sesiones a su preparación, puede de alguna forma determinar el carácter del Círculo y favorecer un tipo de participación en él, pues puede haber participantes que estén especialmente interesados en prepararse para estas competiciones, y el hecho de que exista un Círculo Matemático donde se cuide este aspecto, puede ser un factor decisivo para apuntarse en él. Por otra parte también puede favorecer el aumento en las competiciones. Todas estas cuestiones son importanes a la hora de querer dotar al Círculo Matemático de un determinado carácter específico.

### 4.2. Los contenidos de las sesiones de un Círculo Matemático.

Es característico en un Círculo Matemático que sus sesiones se organicen en torno a temas. A la hora de elegir uno de estos temas, y los problemas en los que se va a desarrollar la sesión, lo más conveniente, por encima del dominio que sobre él tenga el instructor, ya que puede estar más o menos familiarizado con ellos, será que realmente le resulte atractivo y motivador, es decir, debe sentirse fascinado con aquello sobre lo que se va a tratar y trabajar, para poder trasladar estas sensaciones a sus estudiantes, en los que tratará de despertar esta misma visión. El instructor deberá conocer el material que va a utilizar en profundidad, pero no para revelar las cuestiones inmediatamente, sino todo lo contrario, para ser capaz de generar el diálogo, las preguntas, las pistas, que sean necesarias que ayuden en las investigaciones. Que en el grupo se acaben estableciendo, discusiones sobre el tema a tratar, será, tal y como hemos visto, algo muy deseable, por encima de la consabida "explicación del instructor" que es en cualquier caso, lo que hay que tratar de evitar.

No hay limitación para la elección de temas a utilizar en un Círculo Matemático. La única premisa, será la de que esté adaptado al nivel de los estudiantes. Existen multitud de publicaciones, desde los orígenes de los Círculos Matemáticos hasta la actualidad, que recogen material elaborado especialmente para las sesiones, la más reciente a través de la la serie **MSRI Mathematical Circles Library**, que cuenta hasta la fecha con 17 volúmenes, y que son una fuente más que suficiente para poder elaborar las sesiones. (American Mathematical Society, 2014) Según Vandervelde (2009): "existen maneras para adaptar prácticamente cualquier tema para que resulte exitoso dentro de una sesión" (p. 73) [traducción del autor]. Otra de las recomendaciones a la hora de elaborar una sesión, es la de no tener la pretensión de querer introducir mucho material, pues los estudiantes tendrán tendencia, en ese caso, a aburrirse o sentirse perdidos (Vandervelde, 2009, p. 73). En vez de ello, será mucho mejor centrarse sobre unas pocas cuestiones, bien acotadas, que se adapten al tiempo de la sesión. Será conveniente tener preparado material de 2 a 3 veces al que se tenga pensado usar para la sesión (McCullough y Davis, 2013, p. 9), y estar dispuesto a usar tan sólo aquel que se estime necesario, que puede ser una parte o todo, según la situación que se presente en la sesión. Como el ritmo que se establezca dependera de varios condicionantes, y dado que no hay necesidad alguna de cumplimentar el tema, lo más importante será dirigir la sesión sin prisas pero avanzando de forma continua.

La elección de temas puede ser muy amplia, pero conviene establecer cierta continuidad e incremento gradual de la dificultad, para ir profundizando en ellos (Stankova y Rike, 2008, p. 288). De entre los posibles temas a elegir, se expone aquí una lista, a modo de ejemplo, utilizada en el Círculo Matemático de Berkeley: Inversiones en el plano, combinatoria, el cubo de Rubik, teoría de números, pruebas, inducción, geometría de masas, números complejos, juegos con invariantes y pisadas, monovariantes, reconstrucciones geométricas, teoría de nudos, funciones multiplicativas, introducción a teoría de grupos e introducción a inecuaciones (Obtenida de Stankova y Rike, 2008, 2015).

Para una lista más amplia sobre posibles temas se puede consultar la siguiente dirección: https://www.mathcircles.org/Wiki_OrganizingTheAcademics_TranscendantTopics/ Bajo un tema importante, se pueden desarrollar una variedad de subtemas que sean de nuestro interés, por ejemplo, bajo el tema Geometría y Transformaciones, Tom Davis ha desarrollado los siguientes subtemas de su interés:

Problemas de geometría en general (mayormente de áreas)
Inversión en una circunferencia
Cálculo de áreas de polígonos
Introducción a Matrices
Transformaciones geométricas con matrices
Construcciones geométricas clásicas
Geometría de la circunferencia
4 puntos sobre una circunferencia
Geometría para concursos
Geometría y geografía
Polígonos
Teorema de Pick
Teorema de Euler
Geometría proyectiva
Coordenadas homogéneas para gráficos informáticos
Coordenadas homogéneas
Curvas Spline
(Traducción del autor de Davis, 2014)

Dentro de la elección de problemas para el tema propuesto, es conveniente tener en cuenta las posibles orientaciones que se derivan al utilizar unos u otros, entre las cuales establecemos la siguiente lista:

Problemas interactivos
Problemas abiertos (algo para cada uno)
Problemas de explicación sencilla
Problemas interesantes
Problemas con elementos manipulables
Problemas para fortalecer el trabajo en equipo
Problemas orientados a la experimentación
Problemas con trucos para desarrollar métodos
(Traducción del autor de McCullough y Davis, 2013, p. 11)

Es por todo ello que la configuración de las sesiones de un Círculo Matemático presentan su propia idiosincrasia, y aunque tienen una filosofía y estructura común, pueden diferir ampliamente en algunos aspectos, principalmente en el de los contenidos. Lo importante será, en cualquier caso, el tipo de actividad que se desarrolle, el interés hacia las actividades propuestas de los estudiantes y su involucración personal en las investigaciones matemáticas que se les planteen, de tal forma que las sesiones semanales de los Círculos resulten ser algo fascinante a las que estén deseando acudir para conocer quiénes serán los encargados de proponerles nuevas e interesantes actividades que pongan a prueba una vez más su pensamiento y creatividad matemática.

### 4.3. Acciones frente a comportamientos ocasionales no idóneos de los participantes en una sesión de un Círculo Matemático.

Durante el desarrollo de las sesiones es posible que se den determinados comportamientos no idóneos, y que haya que reconducirlo para que la sesión del Círculo se siga desarrollando de una forma adecuada. Tom Davis propone dos reglas para asistir a una sesión: (1) que el participante tiene que divertirse haciendo matemáticas y (2) que desea estar allí (McCullough y Davis, 2013, p. 10). Teniendo en cuenta que la asistencia es totalmente libre, en principio no se deberían dar situaciones en las que algún participante estuviese allí de forma obligada o tuviese un comportamiento habitualmente disruptivo. Nos ocuparemos de los casos en los que participantes que sí desean estar allí puedan estar aburriéndose, sintiéndose perdidos, o tener un comportamiento ocasionalmente disruptivo. En esos casos es conveniente tener una serie de pautas adquiridas por parte del instructor para que no se vea alterado el clima de la sesión. Algunas de ellas son:

Para los participantes que se sientan perdidos:

- Ponerlo en pareja con otro y explicarles que tienen que trabajar en equipo.
- Tener problemas más sencillos preparados para ellos o estrategias para que puedan resolver algún problema más difícil (hacer un seguimiento de lo que realizan).
- Tener buenas pistas preparadas. Preguntar por ellas al resto del grupo.
- Formar un grupo con todos ellos para que intenten resolver el problema. Se les puede ofrecer una versión más simple del mismo.

Para los participantes que se aburran:

- Intentar involucrarles. Puede ser que simplemente no entiendan el problema, o que les resulte muy difícil o muy fácil. Si es muy fácil, hay que asegurarse de que entiendan todos los pasos y sepan explicarlos.
- Darles más problemas y algunos que sean más difíciles.
- Preguntarles qué es lo que saben o han descubierto. Puede que simplemente necesiten atención. Hay que intentar motivarlos con algún reto.

Para los participantes con un comportamiento disruptivo:

- Normalmente suelen estar aburridos, por lo que conviene tener más problemas para ellos.
- Darles tareas a realizar como ayudar a otros compañeros.
- Darles la oportunidad de que expresen sus ideas al grupo, pero siempre bajo el control del instructor.
- Hacerles saber que están teniendo un coportamiento disruptivo y que no es admisible, por lo que tendrán que salir fuera si continuan así.
- Si es debido a que se aburren porque el nivel es demasiado bajo para ellos, cambiarlos a otro grupo con un nivel más alto.

(Traducción del autor de McCullough y Davis, 2013, pp. 10-11)

La forma en que el instructor se dirige a cada uno de los participantes y el trato que les dispensa es esencial para mantener el formato del Círculo Matemático. Por tanto, se deberá aprovechar cada oportunidad que se presente para que así sea, y tener siempre algo que ofrecer en respuesta a las demandas que se presenten, que hagan que el ritmo no decaiga para ninguno de los participantes, los cuales pueden tener necesidades diferentes. El objetivo es aprender a la vez que todo el grupo se divierte, para lo que es necesario matener una atención y participación con un grado elevado.

### 4.4. Un Círculo Matemático según la tradición rusa.

En un Círculo Matemático, según la tradición rusa, suele haber un número elevado de asistentes o instructores, en proporción al número de alumnos (por ejemplo, puede haber 4 instructores para 15 alumnos), además del instructor principal, que suele ser un matemático con experiencia encargado de dirigir todo el proceso, estableciendo contenidos, pautas y criterios a seguir. Se expondrá aquí, cómo se desarrolla uno de estos Círculos según Dorichenko (2012).

La primera clase del año, suele servir para organizarse ya que se desconoce el número de asistentes. En ella se reparte un hoja de problemas sencillos con el objetivo de averiguar que es lo que saben, y lo que no, los futuros integrantes de los Círculos. Además, se les advierte del propósito para que puedan realizarlo con total tranquilidad. Se les menciona que no se trata de ninguna evaluación y de que no pasa absolutamente nada, incluso si no resuelven ningún problema. Si por el contrario algún alumno resuelve los problemas rápidamente, se le podrán dar problemas adicionales. En esta primera sesión, también se les pide que rellenen un cuestionario, para poder asignarles un grupo en la siguiente sesión. Si los estudiantes no están conformes podrán cambiarse de grupo. Durante las siguientes sesiones, todos los nuevos asistentes, son asignados a un nuevo grupo, donde se les ofrecerá una primera sesión, para después integrarlos dentro de los grupos que ya se hayan conformado. Otra opción, es ocupar un auditorio con todos los alumnos nuevos que vayan llegando, y según sean capaces de resolver ciertos problemas se van integrando en los grupos. Los alumnos suelen tener diferentes niveles, por lo que los grupos se pueden conformar en base a ellos, lo que conlleva un mayor trabajo, pues hay que elaborar diferentes hojas de problemas. Otra opción, es elaborar una única hoja de problemas con diferentes niveles de dificultad.

Cada sesión tiene una duración de 2 horas y se desarolla en dos partes con un descanso en el medio. En ellas las explicaciones del instructor tienen que reducirse al mínimo, ocupando no más de 20 min., dedicando el resto del tiempo a la resolución de problemas y a discusiones personales con los instructores. Siempre será el instructor principal el encargado de conducir la sesión, ofreciendo diferentes posibilidades según estime conveniente, como dar alguna pista, proponer algún problema divertido, pedir a alguien que salga a exponer su solución o jugar a algún juego.

Durante una sesión se entregará una hoja de problemas para que los estudiantes intenten resolverla. Se les aconseja, que no tienen porqué hacerlo en orden y que se pasen de uno a otro, si no son capaces de resolverlo o se han quedado atascados en un cierto punto. Cuando un estudiante, tiene algo que preguntar o discutir, puede levantar la mano y un asistente acudirá. En estas conversaciones el instructor tratará de orientar al estudiante, que puede encontrarse en diferentes tipos de situaciones frente a la resolución. Por ejemplo, puede suceder que sepa cómo encontrar la solución pero debido a su falta de experiencia no sepa como explicarla correctamente. Justo antes del descanso se entregan las soluciones de los problemas de la sesión previa y se discuten con el grupo. De esta forma, se ha dado tiempo a que hayan podido terminar de trabajar sobre ellos. A su vez se pueden hacer comentarios sobre la nueva tanda de problemas.

A los estudiantes que hayan terminado con los problemas se les entregarán problemas adicionales. También en el caso de que les quede uno en el que se encuentren atascados. Esta hoja de problemas adicionales se puede entregar al final de la sesión a todo aquel que la solicite. Los problemas que queden sin resolver, son libres de resolverlos hasta la próxima sesión. Para los problemas más difíciles se suele dar más tiempo hasta que se exponen en grupo, para que todo el mundo tenga el tiempo suficiente de pensar en ellos. Lo más importante, es que durante las sesiones, todos hayan podido discutir con los instructores suficientemente sobre los problemas. Al final de la sesión el progreso de cada estudiante se registra, anotando los problemas que han sido capaces de resolver, de esta forma se puede evaluar cómo se ha desarrollado la sesión, conociendo que problemas han sido más dificultosos y si es necesario hacer algún cambio.

En la resolución de problemas, a diferencia de los ejercicios, no existe un camino previo evidente en forma de algoritmo que podamos aplicar. Las habilidades necesarias para resolver cualquier problema son de otra índole, y requieren de gran creatividad y pensamiento lateral. Por ello, las sesiones en los Círculos son diferentes a las de una clase ordinaria, la cual tiene otros objetivos. Si se explicase como se resuelven un determinado tipo de problemas y luego se ofreciesen problemas similares para resolverlos, se estaría dejando poco espacio para desarrollar esta capacidad. Por ello, el método, depende precisamente de la elección de los problemas, que desde un inicio, y sin apenas conocimiento previo, los estudiantes deben tratar de resolver. Un hoja característica de problemas suele contener un primer problema sencillo que se puede resolver de diferentes maneras. Una buena idea es proponer algún problema que intuitivamente tenga una respuesta evidente, pero que para nuestro asombro, sea falsa. La hoja de problemas suele estar planteada en torno a un tema, proponiendo además un problema de nivel bajo-medio para introducir el nuevo tema de la próxima sesión, que puede ser presentado cuando, cerca del descanso, se expongan los problemas de la sesión anterior. Si durante varias sesiones se está tratando un mismo tema, en todas ellas se pueden ir ofreciendo un problema introductorio sobre el siguiente.

La hoja de problemas puede contener problemas repetitivos, es decir, el mismo problema presentado de diferentes maneras, de tal forma, que en un principio los estudiantes los resuelvan como si se tratase, cada uno de ellos, de un problema novedoso, teniendo que llegar a la conclusión final de que en realidad no es así. También se les puede dar una cadena de problemas a lo largo de diferentes sesiones, donde un problema ayuda a resolver el siguiente. Otra opción consiste en entregarles un problema sobre un nuevo tema que sea difícil y atractivo, dándoles el tiempo suficiente, a través de diferentes sesiones para que lo resuelvan. En cada sesión se puede discutir sobre el problema, dando lugar a una lluvia de ideas por parte de los estudiantes, que ayude en la resolución. En lo grupos con más estudiantes de menor grado, es una buena idea ofrecer algún que otro juego matemático entre los problemas, de forma que puedan jugar con el instructor, con algún compañero o resolverlo por ellos mismos. Se pueden aprovechar estos juegos para incorporar algunas ideas, para ello conviene adoptar las estrategias más adecuadas, como aquella que es perdedora para demostrarle algo importante. Otra idea, es incorporar en cada sesión un problema geométrico.

El objetivo de un Círculo Matemático, no es aprender una determinada cantidad de matemáticas, sino más bien hacer que los estudiantes se sientan atraídos hacia ellas descubriendo todo lo que tienen de interesante a la vez que se les enseña a razonar y distinguir cuando se ha encontrado una solución. Lo principal es que los alumnos se sientan interesados, para así poder trabajar con ellos, pues aprender a resolver problemas, requiere de un trabajo constante hasta lograr un cierto grado de confianza a la hora de enfrentarse a ellos. Por tanto, el Círculo Matemático debe ser capaz de capacitar en este aspecto a los estudiantes.

Es muy importante toda la preparación previa de las sesiones. Para ello, todos los problemas deben haber sido resueltos previamente por los instructores y asistentes, y si no fueran capaces, consultando y aprendiendo la solución. El instructor principal debería ser capaz de resolver todos los problemas, contar con profundos conocimientos y experiencia. Durante las sesiones los instructores no deben tener miedo de expresar su desconocimiento sobre algo, pues parte del Círculo Matemático, es que todo el mundo se sienta libre de expresar sus conclusiones, consituyendo los errores una fuente de aprendizaje continua. El instructor puede ofrecer respuestas a los estudiantes durante las sucesivas sesiones.

Si las sesiones discurren adecuadamente, los estudiantes suelen tener una actitud animosa entre ellos, pero puede suceder que no sea así, y que una tanda de problemas no resulte ser tan divertida como se esperaba, bien porque sea demasiado difícil o porque sencillamente no les resulte interesante. En ese caso, conviene tener algunos problemas preparados que puedan ayudar a incentivar la dinámica del grupo. Otra posiblidad es formar pequeños grupos, e ir proporcionándoles pistas para que se animen. Es importante tener paciencia con los estudiantes, pues inicialmente pueden ser inexpertos y tener pocos conocimientos matemáticos, pero pueden poseer una gran capacidad de aprendizaje, y desarrollarla será el objetivo de los instructores. Que un estudiante sea capaz de adquirir un determinado conocimiento por ellos mismos a través de un determinado problema que les resulte revelador debe ser lo más gratificante para un instructor, así que el objetivo, será que esto suceda el mayor número de veces posible. El ambiente informal, es propicio para que aprendan a saber exponer sus conclusiones frente a los demás, además de escuchar y ser capaces de mantener una discusión con el instructor o sus compañeros.

## Conclusiones

Un primer acercamiento a la historia de la educación moderna y contemporánea de Rusia sirve (además de para conocer en detalle el contexto en el que se desarrolla la posterior investigación) para entender la problemática que conlleva la formalización de la educación matemática que históricamente se va gestando a base de reformas y contrarreformas en la búsqueda de dotarla del sentido que la sociedad del momento considera más oportuno y rescatarla así de los “errores” del pasado. En su formalización se van estableciendo debates en torno a si en la educación matemática debe primar la asimilación de contenidos muy específicos o la ejercitación de determinados procesos inherentes a la práctica de la investigación matemática, si debe ser un estudio en sí misma o si debe estar orientada para poder ser aplicada en otras disciplinas o hacia la vida (o incluso estudiarse en su relación con otros conocimientos en el desarrollo de proyectos concretos), si tiene interés porque desarrolla las capacidades de razonamiento y abstracción del individuo o sólo en la medida en que tiene un sentido utilitario para la sociedad, o si su estudio debe realizarse a través de una aproximación intuitiva o mediante una exposición rigurosamente formalizada desde sus inicios. Todos estos debates, y algunos más, constituyen una problemática, que como vemos, no es ajena a las cuestiones actuales y la respuesta a ellos es lo que origina el continuo cambio de los currículos. Este acercamiento sirve también para entender como una vez que la educación matemática queda plasmada en el currículo, se atienden algunos aspectos pero a su vez otros, quizás no menos importantes quedan desatendidos. De esta forma tiene pleno sentido pensar que la educación matemática pueda desarrollarse en campos más amplios y ambiciosos que la educación meramente ordinaria y de carácter obligatorio.

Entender bajo que fórmulas y posibilidades puede llevarse a cabo este segundo aspecto, es el objetivo del segundo acercamiento. Gracias a la perspectiva histórica que ofrece el desarrollo de la educación matemática extra-curricular y extra-ordinaria en la Unión Soviética entre 1930 y 1985, y cuyas influencias alcanzan a la Rusia actual, es posible indagar sobre otras posibilidades, algunas de ellas con plena vigencia hoy en día, no sólo en Rusia, sino también en otros países europeos y en Estados Unidos. Estas propuestas constituyen popularmente algo novedoso, aunque como se ha visto, no es así. En cualquier caso, estudiar estos programas, sus aciertos y sus fracasos, puede ser una baza importante para mejorar el actual transcurso de la educación matemática. Si no directamente, es decir, desarrollando programas más ambiciosos y con un carácter optativo, lo cual, aunque muy deseable requeriría de un esfuerzo social y económico importante, sí al menos indirectamente, lo que se puede conseguir aprendiendo de algunas de sus prácticas y formas de proceder para poder aplicarlo en nuestras aulas. Nos referimos al objeto final de la investigación, el de los Círculos Matemáticos, los cuales constituyen, tal y como se ha visto una experiencia con un carácter extra-curricular que cuenta ya con un siglo de historia. Experiencia, que por su importancia y calidad en el proceso didáctico, ha cosechado inmejorables resultados, y de la cual se puede y se debería, aprender mucho.

Para aprender sobre la experiencia de los Círculos Matemáticos, lo mejor es remitirse a los ejemplos sobre los mismos en la actualidad. Lo deseable sería poder participar directamente en alguno de ellos, para lo cual primeramente se hace necesario conocer qué Círculos Matemáticos existen en la actualidad y que tipo de organización, funcionamiento y propósitos los impulsan. Ese y no otro es el objetivo de un tercer acercamiento, conocer como se organizan estas prácticas desde un primer momento y que posibilidades reales ofrecen, aportando datos relevantes que lo acompañen. Puesta la visión sobre ellos, se puede concluir que se siguen desarrollando gracias a iniciativas individuales que recuperan estas prácticas para ponerlas de rabiosa actualidad. Con la expectativa de la mejora educativa matemática existente, se empeñan en crear estos Círculos Matemáticos, como algo puntual y en vistas de desarrollo. Afortunadamente la demanda parece ser lo suficientemente importante para que estos Círculos, una vez instaurados, puedan seguir creciendo hasta que finalmente, instituciones públicas, empresas privadas y fundamentalmente universidades, acaben por auspiciarlas, dotando así al programa con la suficiente consistencia y garantías para que en el futuro su programa siga ampliándose y llegando cada vez a más alumnos. Por todo ello, se podría considerar incorporar esta práctica en colaboración entre los institutos y las universidades, ayudando no sólo a mejorar en general la educación matemática de los jóvenes, si no, ofreciendo una posibilidad de gran valor para todos aquellos que se sientan particularmente atraídos hacia las matemáticas, y que de no ser por estos programas, no podrían encontrar una respuesta de calidad.

Teniendo claro cuáles pueden ser los objetivos y los beneficios de la formalización de los Círculos Matemáticos, un cuarto y último acercamiento se dirige hacia como debería formalizarse su proceso didáctico. Para ello, es importante conocer en profundidad las posibilidades en que puede organizarse una sesión, que fundamentalmente ha de hacerse con un constante carácter investigador que ayude a desarrollar en los participantes sus habilidades matemáticas frente a la resolución de problemas. Tomando como referencia a personas relevantes por su implicación en la formación de Círculos Matemáticos, se exponen cuestiones que van desde, la promoción y formalización de un Círculo, la recepción inicial de los participantes, la actitud y funciones del instructor o profesor guía, la elección y preparación de los contenidos, el trato frente a comportamientos no idóneos de los participantes, el desarrollo ideal de las sesiones, y todas las posibles actividades e ideas que pueden tener lugar a lo largo del año para configurar una experiencia, además de matemática, social y atractiva. Todo ello debe ayudar al docente en la comprensión de cómo se debe establecer una sesión tipo Círculo, para así planificar sesiones conforme a este estilo y ponerlas en práctica de una forma efectiva.

Esta última idea, es en definitiva el objetivo final. Por ello se propone, a partir de todo el conocimiento teórico proporcionado diseñar posibles futuras sesiones al estilo de un Círculo Matemático, lo cual constituye la conclusión del trabajo, y con el que se quiere acercar al lector, aún más todavía si cabe, el significado de lo que comprende una sesión tipo Círculo. El anexo, además, constituye una propuesta práctica en sí misma, ofreciendo dos posibilidades, una orientada a la propia incorporación de lo aprendido sobre los Círculos Matemáticos en el aula y otra a la formalización de un Círculo Matemático en un Instituto. Creemos que las posibilidades que ofrece una sesión al estilo de un Círculo Matemático deberían de ser tenidas muy en cuenta por el profesorado de secundaria, tanto para aprender de ellas como para ser capaz de incorporarlas con cierta periodicidad en el aula, ya que como se ha visto están orientadas a cubrir una parte muy específica de la formación matemática, de carácter fundamental y que no puede dejar de eludirse. Además, muchos alumnos se sentirán especialmente gratificados de encontrarse, al menos de vez en cuando, con este tipo de sesiones. Desde el punto de vista del profesor, debe ser además del mayor interés, ya que está en su mano el mejorar en todos los aspectos posibles la calidad de la educación matemática que ofrece, cuestión que va más allá de la impartición de los contenidos obligatorios que el currículo vigente establece.

**Daniel Grilo Bartolomé**[1]

[1] Estudiante de la Universidad Autónoma de Madrid del Máster en formación de profesorado de educación secundaria y bachillerato 2014/2015, daniel.grilo@alumni.uam.es / dagribar@gmail.com

## Bibliografía

## Otro material consultado

# Anexo

## Propuesta práctica consistente en sesiones al estilo de un Círculo Matemático.

Hand out a set of problems a week before your session. Not too many, but seductive. Include an easy one and a challenging one.

Be encouraging, even about wrong answers. Find something positive in any attempt, but don't be satisfied until there is a rigorous solution. Wrap up each problem by reviewing the key steps and techniques used.

If the kids don't answer your question immediately, don't just tell them the answer -- let them think. If they're still stuck, give hints, not solutions.

(Davis, 2014)

*Reparte una hoja de problemas una semana antes de la sesión. No demasiados, pero atractivos. Incluye uno fácil y otro que constituya un buen reto.*

*Se alentador, a pesar de las respuestas erróneas. Encuentra algo positivo en cada propuesta, pero no te des por satisfecho hasta que la solución sea rigurosa. Aprovecha cada problema para volver sobre los puntos clave y los métodos empleados.*

*Si los chicos y chicas no responden la cuestión inmediatamente, no les digas la respuesta sin más, déjales que puedan pensar. Si aún así, siguen atascados, proporcionales pistas, no soluciones.*

(Traducción del autor de Davis, 2014)

# Índice

## A.1. Introducción.

En esta parte se proponen una serie de sesiones, a modo de ejemplo y guía, para formalizar una propuesta práctica en referencia a la exposición teórica desarrollada hasta ahora, diseñada conforme al estilo de un Círculo Matemático, con el objetivo de lograr su proceso didáctico característico. Para ello, se ha pensado en dos posibilidades, aunque podría caber alguna más:

- **Clase-sesión tipo Círculo**: opción para ofrecer una de entre las 4 o 5 sesiones de 50 min. de las que habitualmente se dispone semanalmente para un curso de educación secundaria, en un estilo próximo al de un Círculo Matemático. La sesión, además, se adaptará al Currículo obligatorio, por lo que se utilizará alguno de sus temas para planificar las sesiones. Aunque esta sesión semanal se apoye en contenidos ya vistos durante las sesiones ordinarias, la idea es que, entre ellas mismas se establezca una continuidad, de forma que sea claramente diferenciadora del resto de las sesiones ordinarias, donde las explicaciones del profesor y la realización de ejercicios tienen un papel predominante. Para llevarla a cabo se recomienda contar con al menos 2 profesores para un grupo de 25-30 alumnos. Las sesiones consideradas para desarrollar esta propuesta serán las denominadas “sesiones de investigación tipo 1” y “sesiones de problemas tipo 1”.

- **Círculo Matemático:** opción para la formalización de un Círculo Matemático en una Escuela o Instituto de ESO formando un grupo de entre 5 y 15 alumnos provenientes de 1º y 2º de ESO. La asistencia será libre, por lo que se entiende que deberían participar aquellos alumnos que se sientan especialmente atraídos hacia las matemáticas. Lo interesante sería asegurar su continuidad, por lo que planteamos formar un nuevo grupo en el año sucesivo para alumnos de 3º y 4º de ESO, y finalmente, en el tercer año, para alumnos de 1º y 2º de Bachillerato. Se aconseja involucrar a todos aquellos profesores del área de ciencias que puedan estar interesados y buscar apoyos de instituciones públicas. En los Círculos Matemáticos las sesiones se desarrollaran con total independencia del currículo obligatorio, lo que no significa que no puedan tratarse algunos de sus temas pero con otra perspectiva más ambiciosa. Se llevará a cabo en una única sesión semanal de 90 min. con un descanso intermedio de 10 min, por lo que cada parte de la sesión, que pasamos a denominar **semisesión** tendrá una duración de 40 min. Las sesiones consideradas para desarrollar esta propuesta serán las denominadas “sesiones de investigación tipo 2” y “sesiones de problemas tipo 2”.

Los materiales que se proponen para conducir las sesiones, se han diseñado bajo dos estructuras características, “sesiones de investigación tipos 1 y 2” y “sesiones de problemas tipos 1 y 2”, ambas orientadas hacia el objetivo principal de desarrollar el razonamiento matemático y la capacidad para plantear y resolver problemas, y como objetivo secundario el de facilitar la adquisición y empleo del lenguaje simbólico, formal y técnico. A continuación se exponen lo característico de cada una de ellas.

- **Sesiones de investigación**: diseñadas para establecer un proceso guiado de trabajo, que facilite la exploración o investigación sobre algún tema, de modo que se fomenten las estrategias de análisis, de resolución y comprobaciones en torno a algún objeto matemático empleando para ello métodos de demostración y procedimientos característicos de la actividad matemática.

- **Tipo 1 (SI1)**: en conexión con un tema del currículo obligatorio, de una duración de 50 min. y con una propuesta final de investigación a realizar de forma voluntaria durante una o dos de las semanas siguientes en el tiempo libre del alumno. Preferentemente se trabajará de forma individual o por parejas.
- **Tipo 2 (SI2)**: con total libertad en la elección del tema, divididas en dos semisesiones de 40 min. cada una de ellas. La primera de estas semisesiones estará orientada a la exposición de algún tema involucrando en pequeñas investigaciones guiadas a los participantes, para que después, en la siguiente semisesión, bien perteneciente a ese mismo día o al siguiente, puedan dedicarse a realizar sus propias investigaciones sobre las cuestiones que se les propongan, o bien a la exposición y discusión sobre cuestiones realizadas durante semanas anteriores. Preferentemente se trabajará de forma individual o en grupos.

- **Sesiones de problemas**: el objetivo fundamental será fomentar las estrategias frente a la resolución de problemas de forma autónoma por parte del alumno, y donde podrán poner a prueba la competencia matemática que hayan ido adquiriendo mediante las sesiones de investigación que se habrán desarrollado previamente. Por ello, se hará fundamental plantear estas sesiones intercaladas entre las llamadas de investigación proporcionando el ritmo que se estime conveniente. Se hará entrega de la hoja de problemas una semana antes de la fecha en la que se desarrolle la sesión y se reservará la hoja de problemas adicionales.

  - **Tipo 1 (SP1)**: en conexión con un tema del currículo obligatorio, de 50 min. de duración, con una primera parte de trabajo autónomo de 35 min. y una segunda parte de discusión de 15 min. Si se considera necesario se podrá organizar algún grupo. La primera parte consistirá en la resolución de 3 problemas de diferente dificultad y dentro de la Unidad Didáctica que se esté impartiendo en ese momento. Se propondrá 1 problema adicional basado en algún tema correspondiente a Unidades Didácticas vistas anteriormente, a modo de refuerzo, o de la siguiente Unidad Didáctica prevista, con carácter introductorio de la siguiente hoja de problemas. Se contará además, con otros 4 problemas alternativos, también de diferente dificultad y que puedan ser utilizados de diversos modos: para sustituir o facilitar el camino hacia la resolución de alguno de los propuestos durante la sesión, para quienes hayan terminado, y como trabajo adicional en el tiempo libre para quién los solicite. Los problemas no se recogerán, sino que se podrá seguir trabajando en ellos hasta que sean comentados finalmente en una sesión posterior. La segunda parte se dedicará a la discusión de algunos problemas vistos durante las clases-sesiones tipo Círculo previas. Para ello se estima conveniente que sean los propios alumnos quienes traten de exponer los razonamientos y estrategias que hayan utilizado de forma justificada. El profesor, será el encargado de establecer la dinámica adecuada con la clase a través del diálogo y mediante el uso de preguntas y pistas.
  - **Tipo 2 (SP2)**: con total libertad en la elección del tema, y para desarrollar en una semisesión de 40 min. Consistente en la resolución de 3 problemas de diferente dificultad. Para estas sesiones se recomienda contar con un ratio de 1 profesor-guía por cada 5 alumnos, de forma que puedan establecerse todas las discusiones que sean necesarias durante esos 40 min. Se diseñan con total libertad, por lo que cabe la posibilidad de presentar problemas que abarquen diferentes temas. Al igual que las sesiones de problemas tipo I (SP1), se contará con una hoja adicional de otros 3 problemas, que puedan utilizarse también en los modos que ya se han mencionado. 10 min. de la semisesión podrán ser utilizados en la exposición de alguno de los problemas tratados en semisesiones anteriores, en la aclaración sobre alguna cuestión importante, o en la propuesta de algún juego matemático entre grupos.

## A.2. Sesiones de investigación.

### A.2.1. Patrones, inducción matemática: suma de los ángulos interiores de un polígono simple (SI1)[1].

**Desarrollo**

La sesión está diseñada para 1º de ESO, y se requiere haber tratado antes los siguientes contenidos:

- Definición de polígono simple y sus elementos característicos.
- Clasificación de los polígonos según el número de lados o ángulos.
- Polígonos regulares e irregulares.
- Polígonos cóncavos y convexos.
- Demostración de la suma de los ángulos interiores de un triángulo.

<u>Problema 2.1.1</u>

Calcular la suma de los ángulos interiores de los polígonos mostrados mediante la descomposición del polígono en triángulos a través del trazo de diagonales. ¿Existe alguna relación entre el número de triángulos y algún elemento del polígono que se mantenga constante? De ser así, ¿podrías encontrar una expresión para calcular de forma inmediata la suma de los ángulos interiores de cualquier polígono conocidos su número de lados?

Recomendaciones para el desarrollo de la sesión:

En conexión con el estudio de los polígonos, se diseña una sesión en la que se tiene que investigar sobre la suma de los ángulos interiores de un polígono simple. Para ello se proporciona una hoja de trabajo donde hay dibujados polígonos que aumentan en su número de lados progresivamente, ofreciendo una versión cóncava y otra convexa. Se recomienda seguir los siguientes pasos:

- Previamente, deben conocer que los ángulos interiores de un triángulo suman 180º. Primeramente se les pedirá, que conocida esa propiedad, busquen algún procedimiento que al utilizarla, sirva para calcular el área de los polígonos proporcionados. Se les preguntará por qué creen que se les ofrece una versión cóncava y otra convexa, además de un número de polígonos cuyo número de lados va sucesivamente en aumento.
- Una vez que hayan conseguido averiguar, ofreciendo el apoyo y las pistas que se estimen convenientes, el procedimiento de triangulación de los polígonos mediante el trazo de diagonales, se les pedirá que determinen la suma de los ángulos interiores de los polígonos que se les había presentado, en relación al número de triángulos hallados.
- Realizado el cálculo, se les pedirá que busquen si existe alguna relación entre el número de triángulos hallados respecto de alguno de los elementos del polígono.
- Seguidamente se les pedirá conjeturar una fórmula general y comprobarla para otros polígonos.

Sol.:

$$S_{ang\ int} = 180(n-2)\ \ siendo\ n = nº\ de\ lados\ del\ polígono$$

[1] Esta sesión tuve la oportunidad de llevarla a cabo con un grupo de 1º de ESO durante las prácticas del Master correspondientes al módulo específico.

- Dentro de los objetivos de la sesión se encuentra además el de llegar a conclusiones importantes, ya que no vamos a realizar una demostración por inducción, a través de razonamientos más elaborados que constituyan una prueba en sí misma. Como en la Fig. 1, donde en un polígono convexo, triangulando apoyándose en un mismo vértice (el 1), se hace evidente que el número de triángulos está en directa relación con el número de vértices, ya que podemos emparejar cada nuevo triángulo con un vértice, a través de la diagonal trazada, a excepción del último triángulo que cierra el polígono, y que siempre se puede emparejar con el último vértice (el n) de forma que tan sólo el vértice 1 y el vértice 2  quedan libres, y por tanto siempre sobrarán 2 vértices cualquiera que sea el polígono convexo que tracemos.

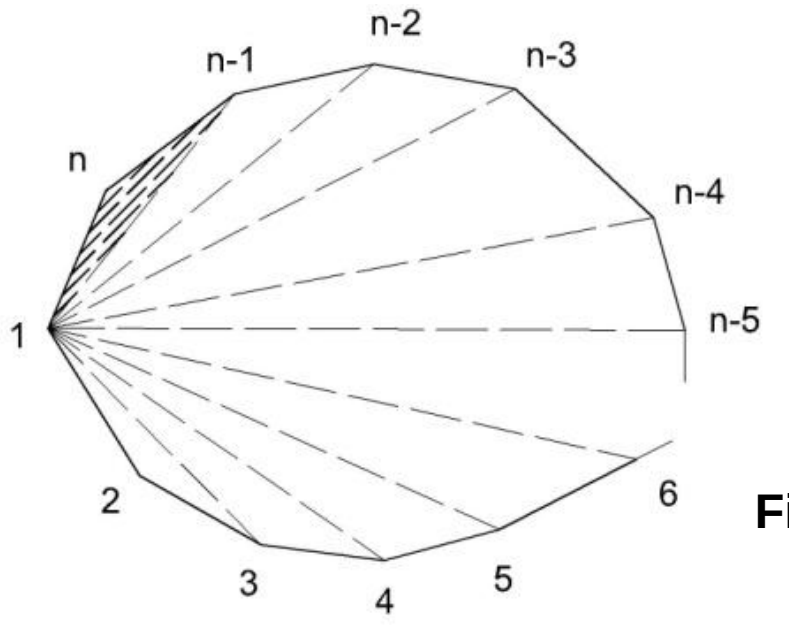

**Figura 1**

- A continuación cabría preguntar, si es posible realizar algo semejante con un polígono de n lados pero esta vez para el caso cóncavo. Podemos plantear un caso arbitrario (Fig. 2), y estudiar qué es lo que sucede. Tratamos de aplicar el mismo método, para lo cual se hace necesario dividir el polígono en un número de subcasos determinados (2 en este ejemplo). Si renombramos los vértices, podemos fácilmente aplicar el método a cada uno de ellos. En uno de ellos, el de la Fig. 3, obtenemos 7-2=5 triángulos, y en el otro, el de la Fig. 4, obtenemos n-5-2=n-7 triángulos. Si sumamos los triángulos de ambos subcasos obtendremos el resultado deseado: 5+n-7=n-2 triángulos.

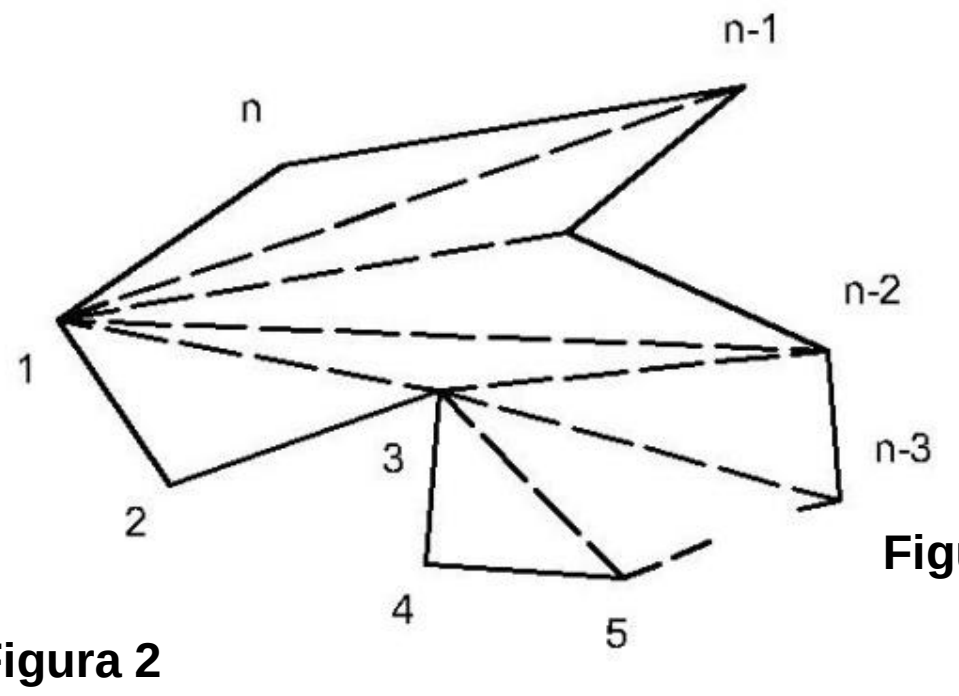

**Figura 2**

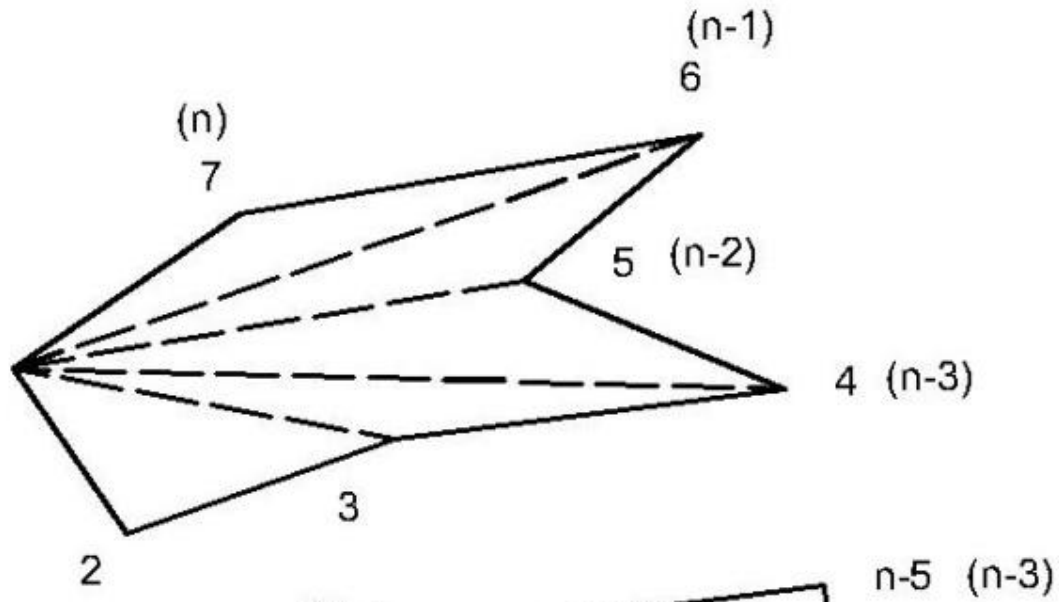

**Figura 3**

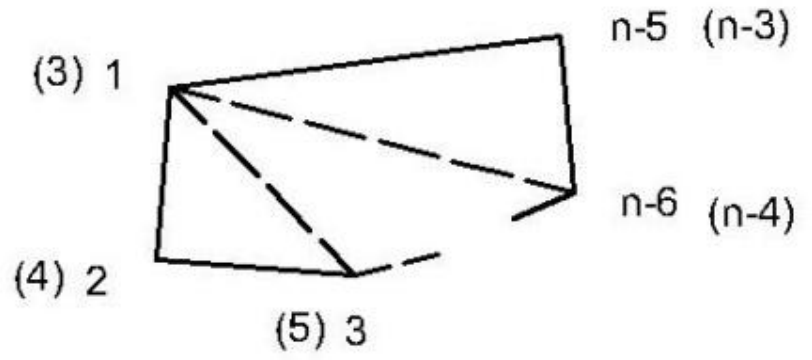

**Figura 4**

- Llegados a este punto parece adecuado pensar, que cualquiera que sea el polígono, podremos reducirlo a otros subcasos conocidos, y por tanto que siempre se obtendrán n-2 triángulos. Sin embargo, habrá que indicar, que aún así es necesario demostrarlo. ¿Es posible generar y demostrar toda una serie de subcasos de tal forma que cualquiera que sea el polígono podamos reducirlos a ellos? La respuesta a esta pregunta puede conducirnos a una nueva investigación más ambiciosa.

- Por último, se preguntará sobre la posibilidad de formular algún otro problema que pueda tener previsión de resolverse en relación a la investigación realizada. La idea es conducirlos hasta un problema de similares características que también habrán de intentar de resolver, y que son una consecuencia del resultado ya obtenido.

Problema 2.1.2.

A partir de la expresión general hallada para calcular la suma de los ángulos interiores de un polígono, determinar la **suma de los ángulos exteriores de un polígono** y el **ángulo interior de un polígono regular.**

Sol.:

$$S_{ang\ ext} = 180n - 180(n-2) \ = 360º \ \ siendo\ n\ el\ nº\ de\ lados\ del\ polígono$$

$$\hat{A}_{int} = \frac{180(n-2)}{n} \ \ siendo\ n\ el\ nº\ de\ lados\ de\ un\ polígno\ regular$$

- Al final se propondrá una investigación avanzada, optativa y a realizar en el tiempo libre, para lo cual contarán con 2 semanas. La investigación se centrará sobre el siguiente problema:

Problema 2.1.3.

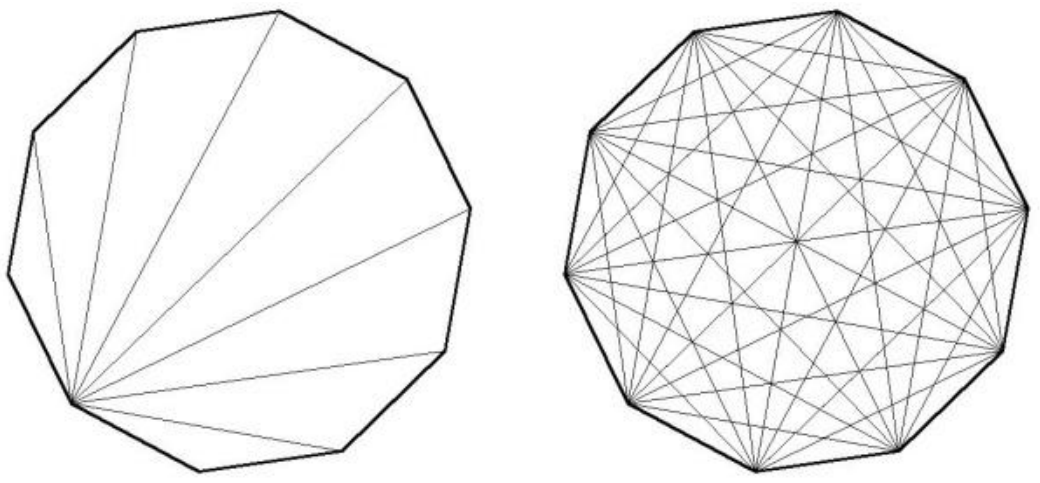

Figura 5

Encontrar una expresión general para calcular el número de diagonales que se pueden trazar desde el vértice de un polígono. Encontrar otra expresión general para el número total de diagonales (Fig. 5). Se puede consultar la siguiente dirección como material de apoyo:

http://www.ceibal.edu.uy/UserFiles/P0001/ODEA/ORIGINAL/110601_diagonales.elp/polgonos_convexos_y_no_convexos.html

Sol.:

$$d_{un\ vértice} = n - 3\ siendo\ n\ el\ nº\ de\ lados\ del\ polígono$$

$$d_{totales} = \frac{n(n-3)}{2}\ siendo\ n\ el\ nº\ de\ lados\ del\ polígono$$

**Objetivos**

- Descomponer en triángulos cualquier superficie poligonal mediante un proceso de triangulación a través del trazo de diagonales sucesivas, donde cada una de ellas comparte al menos un vértice con otra, de forma que se originan triángulos adyacentes.
- Reconocer patrones
- Conjeturar una fórmula general a partir de la observación de casos sucesivos (introducción a la inducción matemática)
- Determinar la suma de los ángulos interiores y exteriores de un polígono, el ángulo interior de un polígono regular, el número de diagonales que se pueden trazar desde un vértice en un polígono, el número total de diagonales que se pueden trazar en un polígono.
- Plantear problemas a partir de otros problemas.
- Expresar de forma algebraica las conclusiones geométricas halladas.
- Realizar demostraciones directas.

**Hoja de trabajo**

(4 lados) Cuadrilátero cóncavo

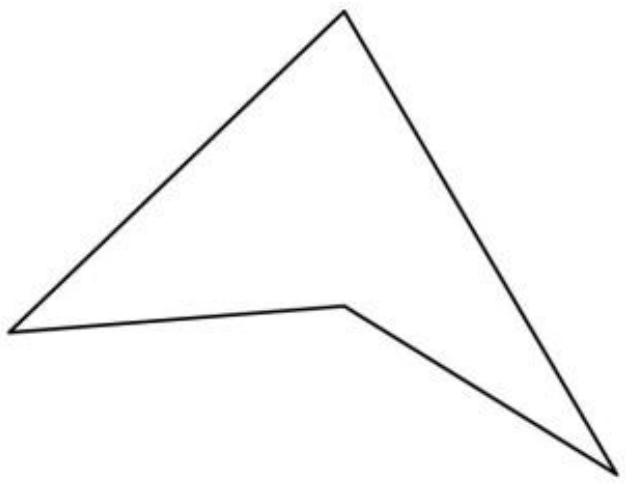

(4 lados) Cuadrilátero convexo

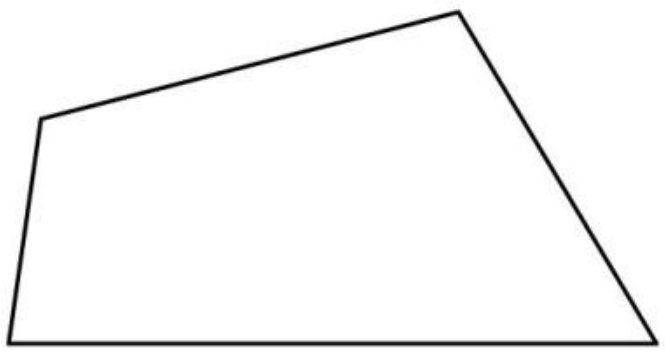

(5 lados) Pentágono cóncavo

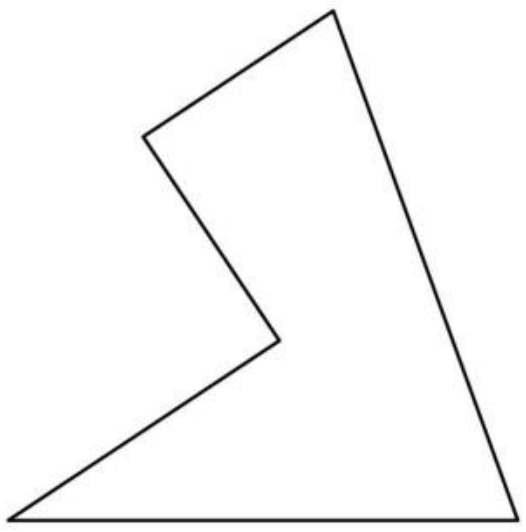

(5 lados) Pentágono convexo

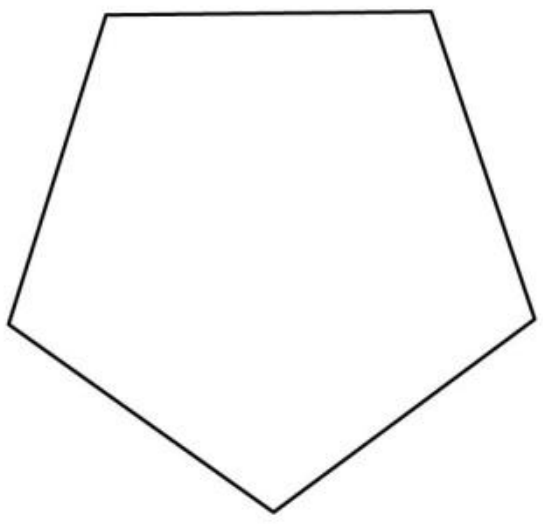

(6 lados) Hexágono cóncavo

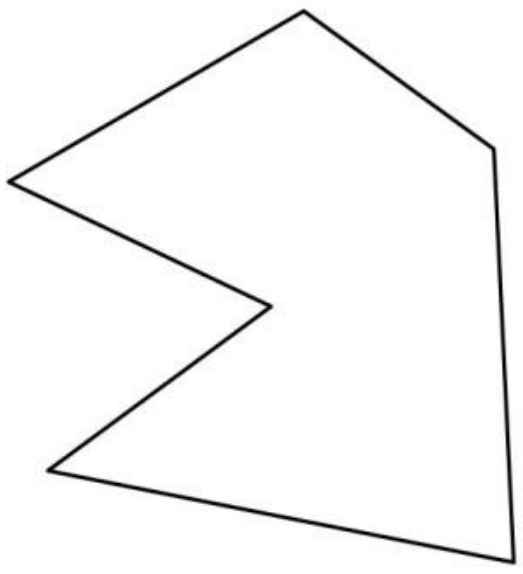

(6 lados) Hexágono convexo

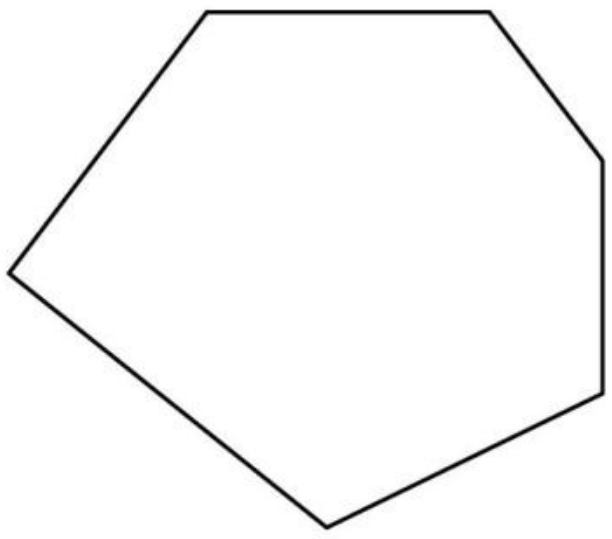

…

(9 lados) Eneágono cóncavo

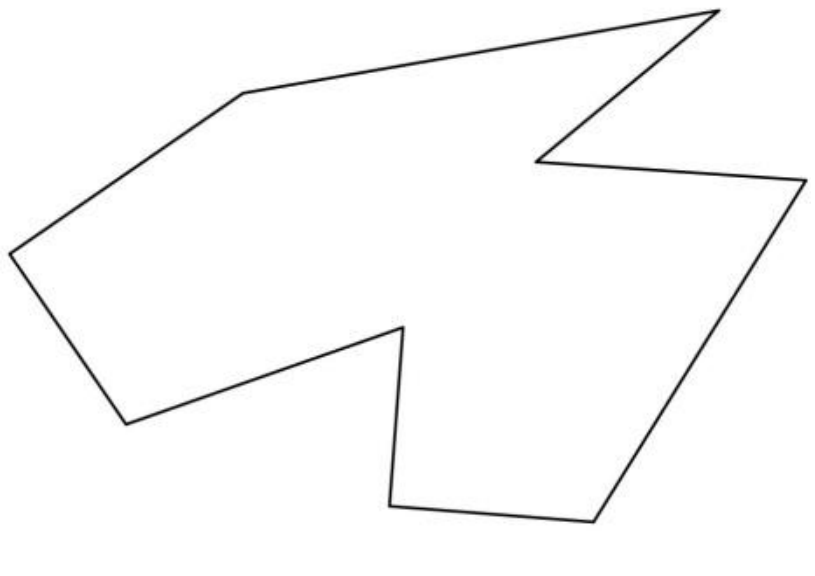

…

(9 lados) Eneágono convexo

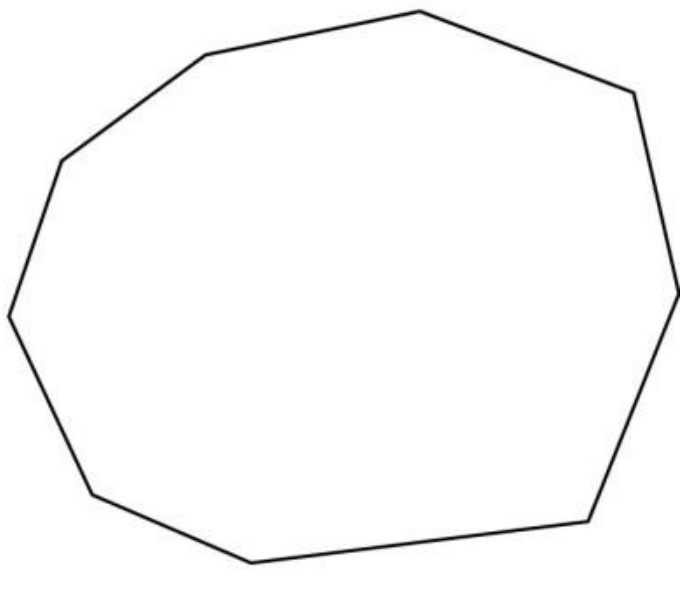

…

(12 lados) Dodecágono cóncavo

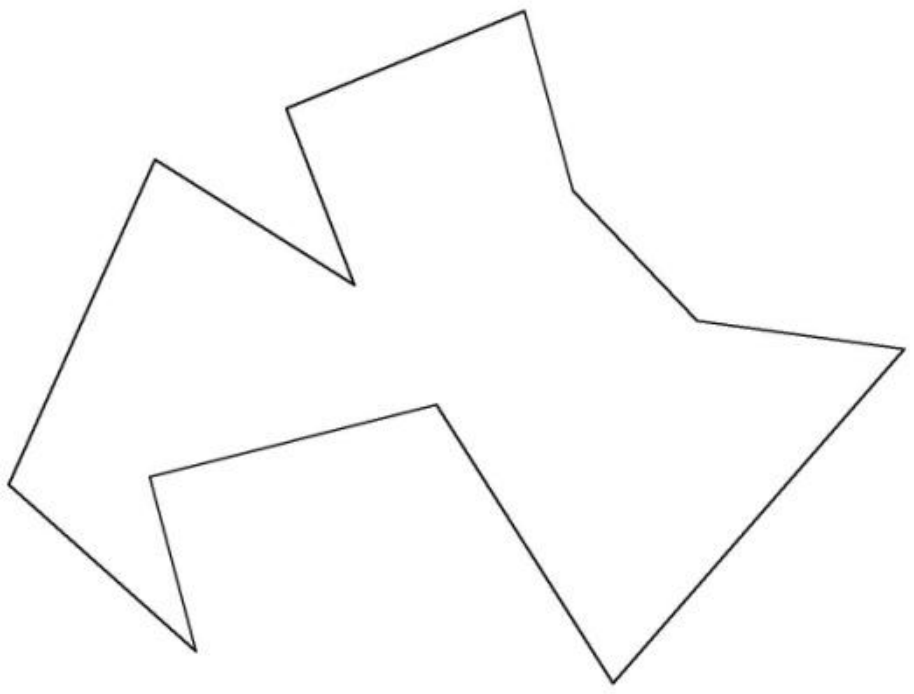

**…**

(12 lados) Dodecágono convexo

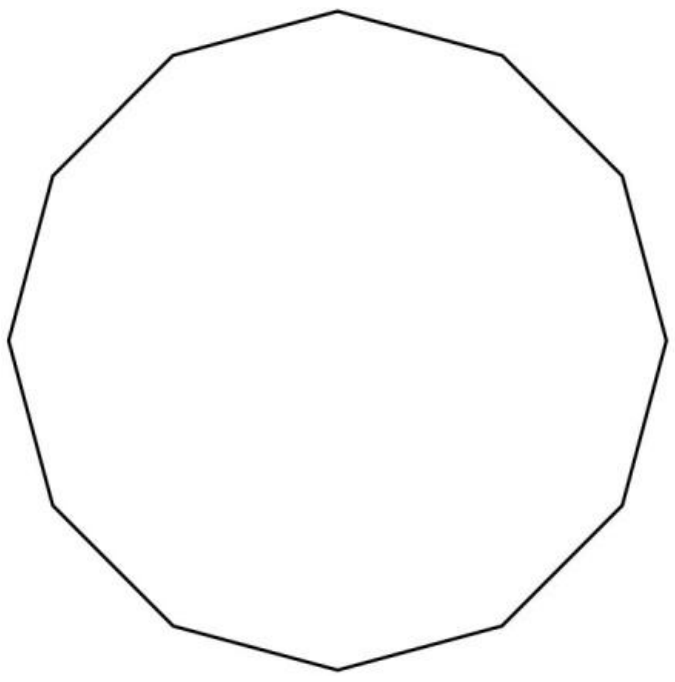

**Figura 6**

## A.2.2. Resolución de problemas mediante grafos (SI2).

### Desarrollo

La sesión está diseñada para 1º de ESO. El propósito es el de introducirles en el empleo de grafos en la resolución de problemas. Para ello se proponen dos problemas, uno sencillo de forma inicial para introducir de forma intuitiva el objeto matemático grafo, con sus elementos característicos; y un segundo problema, de una dificultad mayor, donde se siguen explorando las posibilidades que ofrece un grafo.

Problema 2.2.1

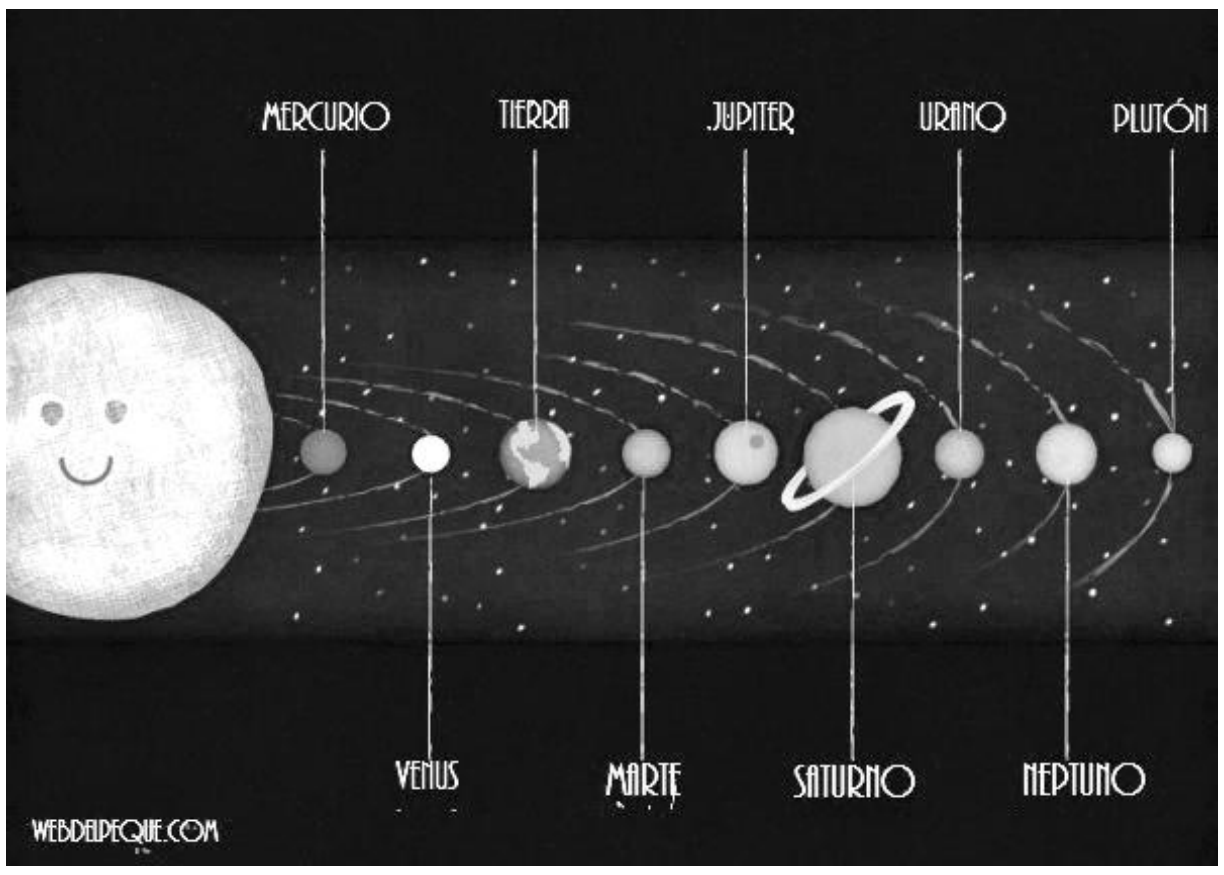

**Figura 7**

Supongamos que en un futuro es posible viajar, a través de cohetes espaciales, dentro del Sistema Solar (Fig. 7). Los planetas que están comunicados, en viaje tanto de ida como de vuelta, son: Tierra-Mercurio, Plutón-Venus, Tierra-Plutón, Plutón-Mercurio, Mercurio-Venus, Urano-Neptuno, Neptuno-Saturno, Saturno-Júpiter, Júpiter-Marte y Marte-Urano. ¿Puede un viajero que parte desde la Tierra llegar a Marte? (Fomin, Genkin, y Itenberg, 1996, p. 39).

Dibujando un diagrama, donde los planetas sean representados por puntos, y la comunicación entre ellos se establezca mediante aristas, es fácil comprobar que la respuesta es no. A continuación, se pueden introducir el objeto 'grafo' nombrando los puntos, como 'vértices' y las uniones entre ellos como 'aristas' que indican las conexiones establecidas entre los vértices.

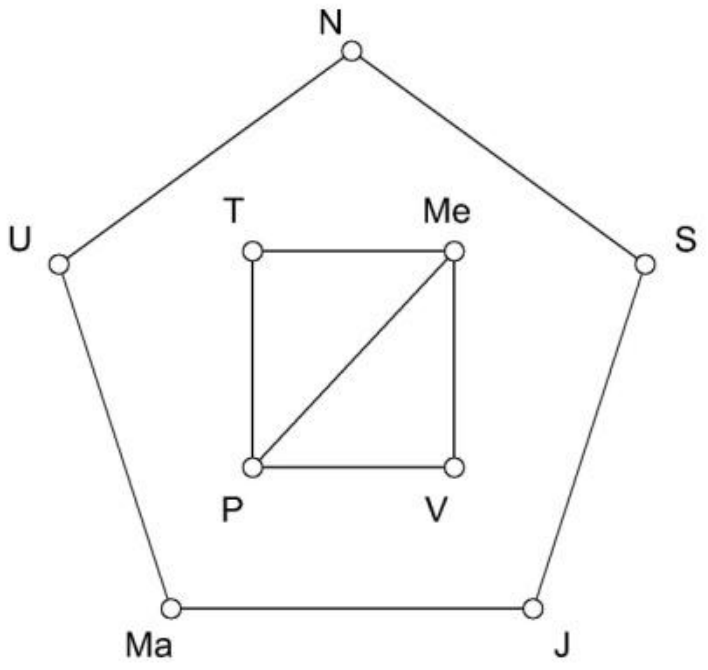

**Figura 8**

Preguntas a realizar:

- ¿cuál es la información más significativa que nos aporta el grafo?
- ¿de entre los viajes posibles, cuales son directos (no hay que viajar a destinos intermedios)?
- ¿se pueden clasificar los vértices según el número de conexiones que tiene establecidas? ¿puedes configurar una lista? (una vez realizado se mencionará que dicho número constituye el 'grado' del vértice, que puede ser par o impar.

A continuación, se propondrá un nuevo problema, esta vez construyendo un grafo dirigido, para que realicen una investigación por su cuenta. Se considera oportuno realizar la actividad en grupos de 3, para que puedan establecer discusiones entre ellos.

Problema 2.2.2.

En un gallinero, cada una de las gallinas es etiquetada como $C_1$, $C_2$, $C_3$, etc. Si consideramos cada pareja de gallinas posibles, siempre sucede que una de ellas pica a la otra (lo cual se indicará con una flecha): por ejemplo $C_1 \rightarrow C_2$ significa que $C_1$ pica a $C_2$. Un posible orden de picadas dentro de un gallinero con 5 gallinas podría ser:

$C_1 \rightarrow C_2$; $C_3 \rightarrow C_1$; $C_4 \rightarrow C_1$; $C_5 \rightarrow C_1$
$C_2 \rightarrow C_3$; $C_2 \rightarrow C_4$; $C_2 \rightarrow C_5$
$C_3 \rightarrow C_4$; $C_5 \rightarrow C_3$
$C_4 \rightarrow C_5$

Se pide:

- Dibujar el grafo (sol.: Fig. 9)
- Si consideramos una gallina ‘dominante’ aquella que bien directamente o a través de otra, que llamaremos gallina ‘soldado’, por ejemplo $C_1 \rightarrow C_2 \rightarrow C_4$, consigue que todas las gallinas hayan sido picoteadas, indicar cuales lo son.
- De entre todas las gallinas, que gallina puede ser considerada la mejor gallina ‘soldado’ y cual la peor.
- Encontrar un gallinero con 6 gallinas, donde cada una de ellas sea ‘dominante’. (sol.: Fig. 10)

(Vandervelde, 2009, pp. 95-102)

Sol.:

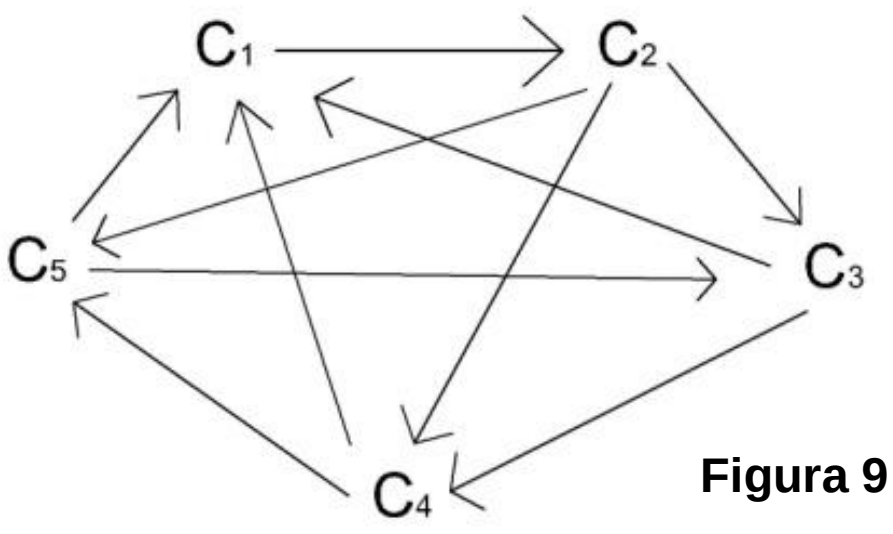

Figura 9

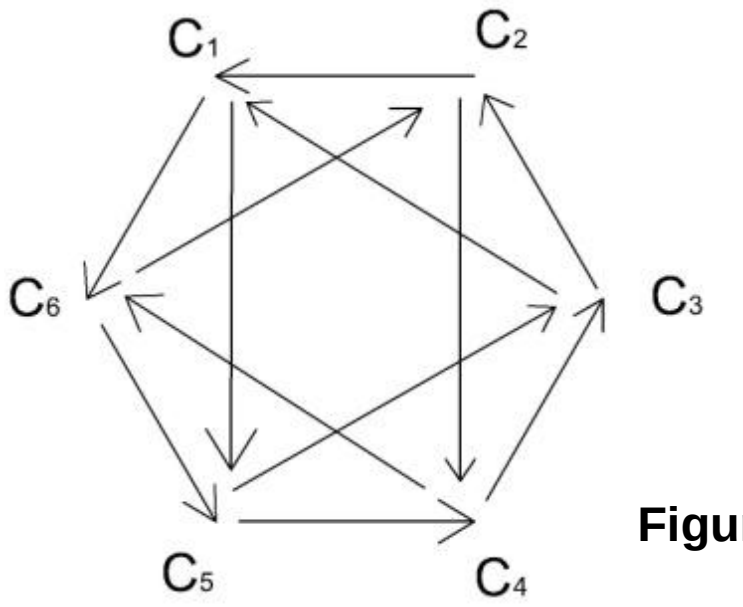

Figura 10

**Propuesta de investigación**

Consistirá en la resolución de dos nuevos problemas, y podrá realizarse en parejas.

<u>Problema 2.2.3.</u>

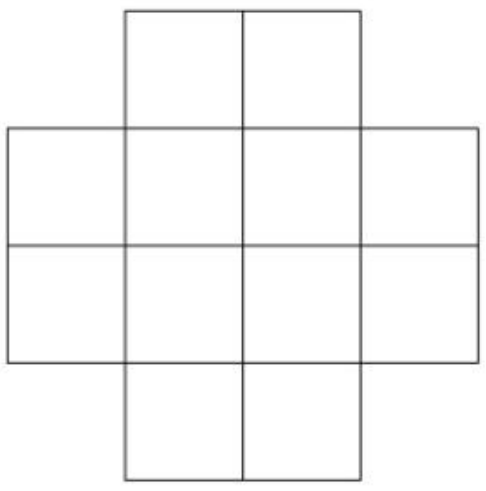

Recortando las esquinas de un tablero de ajedrez de dimensiones 4x4 obtenemos un tablero con forma de cruz (Fig. 11). ¿Puede recorrer este tablero un caballo pasando por cada cuadrado exáctamente una vez y acabar en el mismo cuadrado en el que empezó? (Fomin, Genkin, y Itenberg, 1996, p. 41).

Figura 11

Sol.:

Construyendo un grafo, vemos que ese recorrido sí es posible (Fig. 12):
1-7-5-11-10-4-12-6-8-2-3-9-1

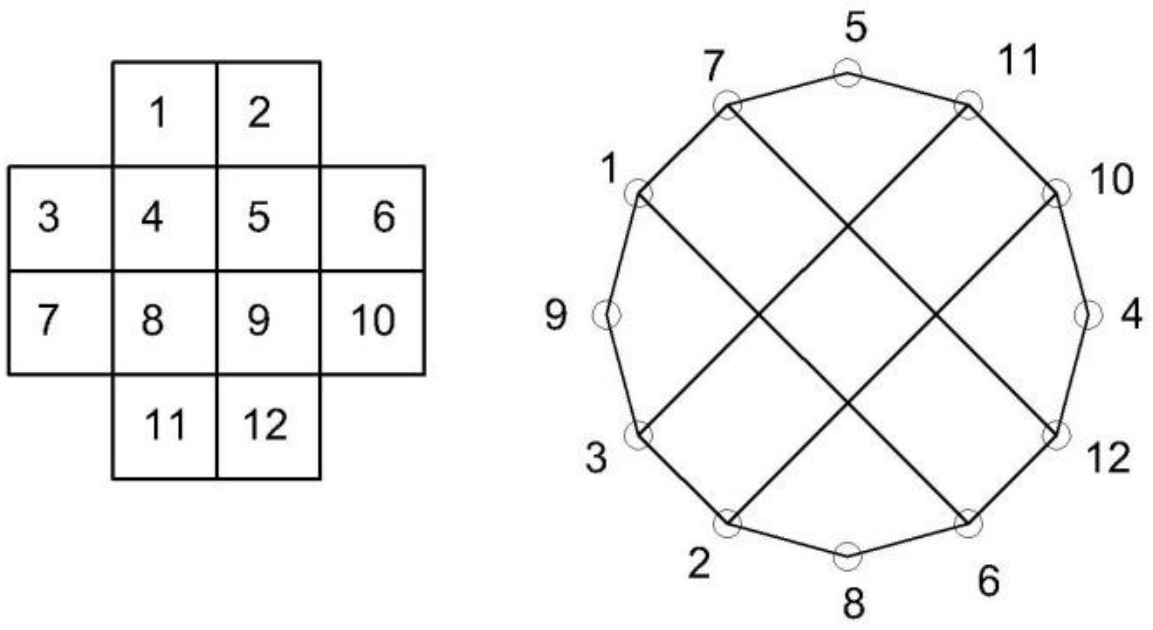

Figura 12

Problema 2.2.4.

Se dice que dos grafos son ‘isomorfos’ si tienen el mismo número de vértices (por ejemplo, n) y es posible numerar los vértices de cada grafo de 1 a n de forma que los vértices del primer grafo estén conectados con una arista si y solo si los vértices del segundo grafo, etiquetados con los mismos números, están conectados por una arista (Fomin, Genkin, y Itenberg, 1996, p. 135).

Teniendo en cuenta esta definición, averigua si los siguientes pares de grafos (Fig. 13, 14 y 15) son isomorfos o no:

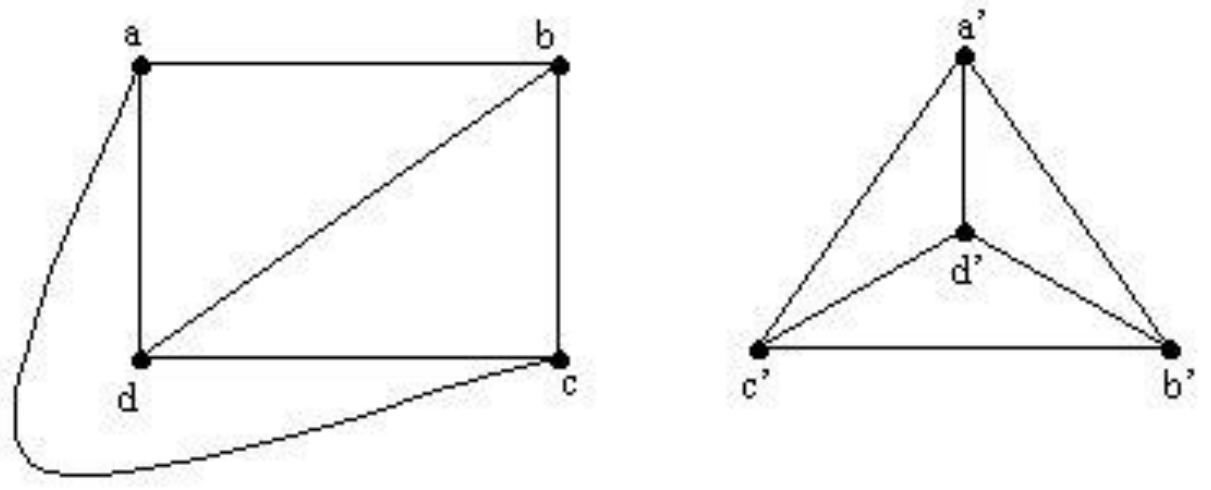

**Figura 13**

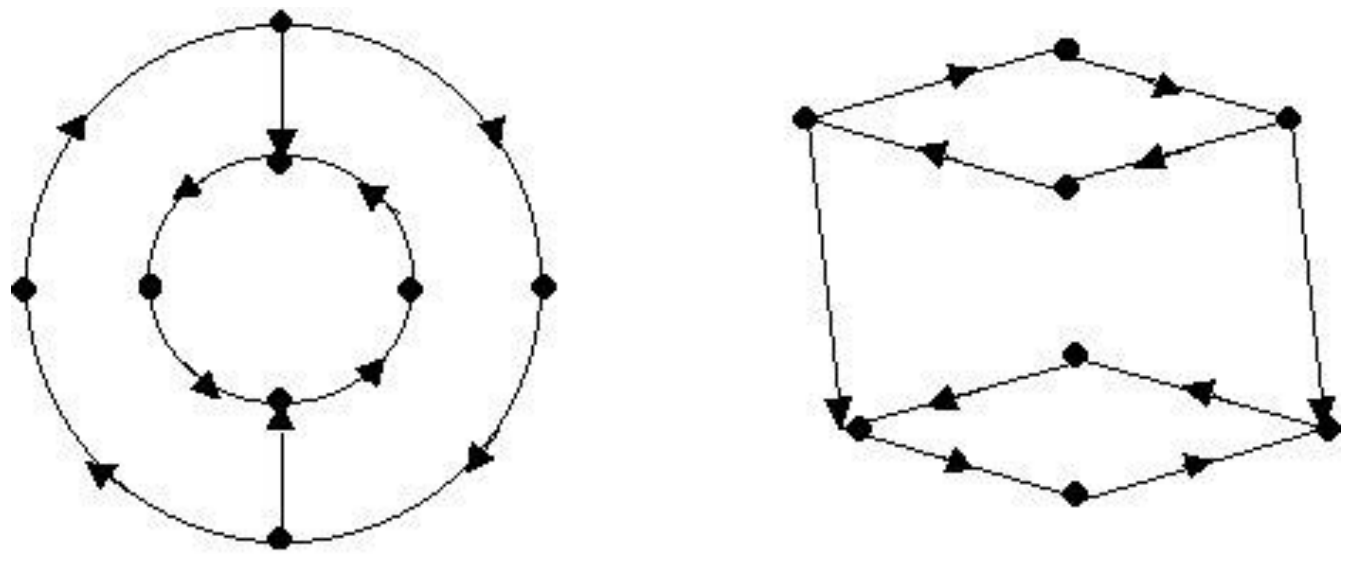

**Figura 14**

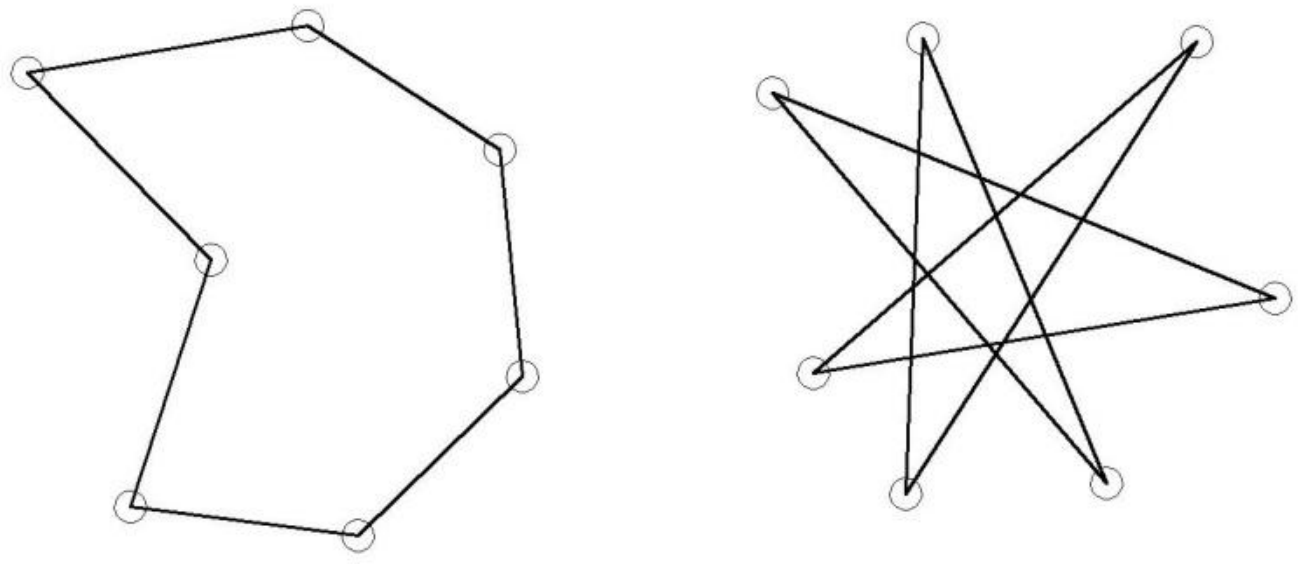

**Figura 15**

**Objetivos**

- Construir el grafo característico de un problema.
- Identificar la información que se nos pide en un grafo.
- Distinguir entre grafos conexos y grafos inconexos.
- Distinguir entre grafos dirigidos (dígrafos) y no dirigidos.
- Identificar y plantear grafos isomorfos.
- Resolver un problema a través de la construcción de un grafo
- Reconocer cuando se cumple la propiedad de transitividad

## A.2.3. Patrones, recursividad, inducción: la torre de Hanoi (SI2).

### Desarrollo

La sesión está diseñada para el grupo de 3º y 4º de ESO y se requiere haber introducido antes los siguientes contenidos de una forma básica:

- Sucesiones numéricas. Progresiones aritméticas y geométricas.
- Principio de inducción

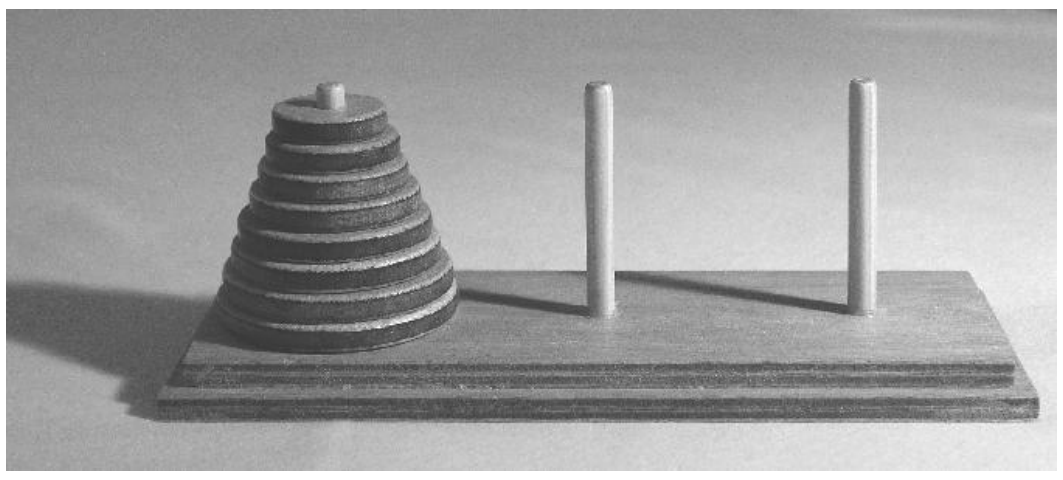

Figura 16

La torre de Hanoi (Fig. 16) es un juego que consiste en pasar unos discos superpuestos, sucesivamente ordenados de mayor a menor tamaño, ubicados en una varilla, a alguna de las otras 2 varillas. Para mover los discos se tienen que seguir las siguientes reglas:

- Sólo se puede mover un disco cada vez.
- Un disco de mayor tamaño que otro no puede superponerse sobre este último.
- Sólo se puede mover el disco que se encuentre en la posición superior de la varilla.

A modo de introducción, se mostrará a los participantes como funciona el juego proponiendo el siguiente problema:

<u>Problema 2.3.1.</u>

Utilizando tan sólo 3 discos, ¿cuál es el menor número de movimientos necesarios para trasladarlos a una de las 2 varillas libres?

Sol. Será fácil comprobar que el número mínimo es el de 7 movimientos, tal y como muestra la Fig. 17

La siguiente tarea, llevada a cabo en parejas, será dar respuesta al problema:

<u>Problema 2.3.2.</u>

Utilizando $n \geq 1$ discos, ¿cuál es el menor número de movimientos necesarios para trasladarlos a una de las 2 varillas libres? (Citado en Benito González y Madonna, p. 241).

Para ello tendrán que realizar su propia investigación averiguando cual es ese número, cuando en vez de tan sólo 3 discos, se trata de 1,2,3,4,5,etc. Durante el proceso, el objetivo es que los participantes puedan llegar a algún tipo de conclusión relevante, como por ejemplo, que para poder mover el último disco, es necesario haber movido previamente todos los discos superpuestos a él, a otra varilla.

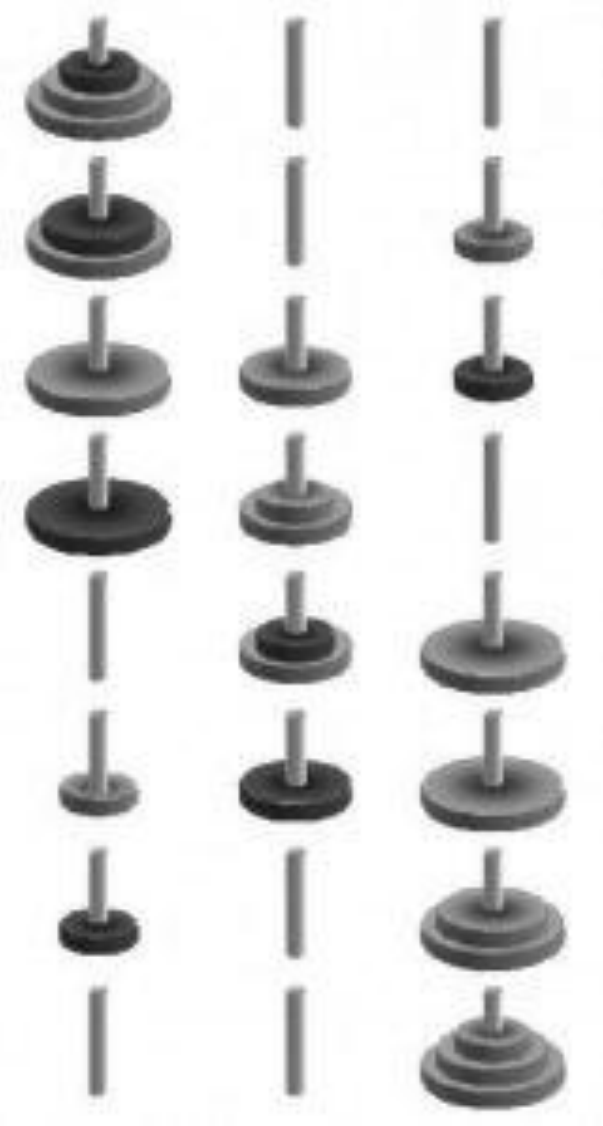

Figura 17

Llegados a esta conclusión, y con los resultados obtenidos para los otros números de discos, se preguntará si a partir de ella, se puede plantear alguna fórmula que sea generalizadora del proceso.

La respuesta es afirmativa, ya que partiendo de *n* discos situados en una varilla, siempre se tendrán que pasar *n−1* a otra de las dos varillas, para poder pasar el disco mayor a la varilla que nos quede libre. Después de esto, nuevamente, tendremos que pasar los *n–1* discos sobre el disco mayor, para que finalmente todos los discos se hayan colocado en una varilla diferente de la incial. Es importante comprender que al vernos obligados a realizar estos movimientos, estamos demostrando que no puede haber una forma que implique mover menos discos. De esta forma, se observa que los movimientos de los discos vienen **definidos por recursividad**, lo que da lugar a una sucesión recurrente:

- Para mover 1 disco se precisa de 1 movimiento: $a_1=1$
- Para mover *n>1* discos se precisa mover 2 veces *n−1* discos y 1 vez 1 disco:
  $a_n = a_{n-1} + 1 + a_{n-1} = 2\cdot a_{n-1} + 1$

Lo que equivale a realizar los siguientes movimientos según el número de discos:

$a_1 = 0+1+0 = 1$
$a_2 = 1+1+1 = 3$
$a_3 = 3+1+3 = 7$
$a_4 = 7+1+7 = 15$
.....
$a_n = a_{n-1} + 1 + a_{n-1} = 1 + 2\cdot a_{n-1}$

(Benito González y Madonna, p. 243)

Mediante simple iteración se obtiene que:

$a_n = 1 + 2^1 + \ldots + 2^{n-2} + 2^{n-1} = \sum_{j=0}^{n-1} 2^j$

Para realizar la iteración, se propondrá que se realice en pasos sucesivos, es decir, primero para n=2 donde se obtendría $a_2=1+2$, después para n=3, donde se obtendría, $a_3=1+2+2^2$, depués para n=4, donde se obtendría: $a_4=1+2+2^2+2^3$, etc., de forma que resulte evidente formularla para el caso enésimo.

A continuación, se comentará, que al aplicar el método de iteración, la expresión que resulta es la **suma parcial de una serie geométrica** de razón común r=2 y con $a = 1$ como primer término, que viene dada por la expresión: $U_n = a + ar + ar^2 + ar^3 + \cdots + ar^{n-1}$

Si multiplicamos a ambos lados de la expresión por la razón r, obtenemos la expresión:
$r \cdot U_n = ar + ar^2 + ar^3 + \cdots + ar^n$

Restando a la primera expresión la segunda, queda:

$$U_n - r \cdot U_n = a - ar^n \Leftrightarrow U_n(1-r) = a(1-r^n) \Leftrightarrow U_n = a\frac{1-r^n}{1-r}$$

(Rey Pastor, Pi Calleja, y Trejo, 1959, pp. 296-297)

Luego fácilmente se puede obtener una expresión explícita:

$$\sum_{j=0}^{n-1} 2^j = \frac{1-2^n}{1-2} = \frac{1-2^n}{-1} = 2^n - 1$$

(Fernández Gallardo y Fernández Pérez, 2003)

Quizás más sencillo, sea proponer hallar una conjetura para la expresión explícita a partir de los resultados que se iban mostrando, para después tratar de demostrarla utilizando el **principio de inducción** (véase Delgado Pineda y Muñoz Bouzo, 2010, pp. 39, 65, 66).

$a_1 = 0+1+0 = 1$
$a_2 = 1+1+1 = 3 = 4-1 = 2^2-1$
$a_3 = 3+1+3 = 7 = 8-1 = 2^3-1$
$a_4 = 7+1+7 = 15 = 16-1 = 2^4-1$
.....
$a_n = a_{n-1} + 1 + a_{n-1} = 1 + 2\cdot a_{n-1}$ conjetura: ¿$a_n = 2^n - 1$?

i) Lo comprobamos para n=1: $a_1 = 2^1 - 1 = 1$; luego cierto (definido por recurrencia)
ii) Para n, suponemos cierto que: $a_n = 2^n - 1$ (hipótesis de inducción)
iii) Lo tratamos de demostrar para n+1: ¿$a_{n+1} = 2^{n+1} - 1$?

Conocemos que $a_{n+1} = 2a_n + 1$ pues está definido por recurrencia, luego:

$a_{n+1} = 2a_n + 1 = 2(2^n - 1) + 1 = 2 \cdot 2^n - 2 + 1 = 2^{n+1} - 1$; luego cierto

Por tanto: $a_n = 2^n - 1 \quad \forall \quad n \geq 1$ Q.E.D.

Por último podemos preguntar si se les ocurre, alguna otra forma de ampliar el problema propuesto. Lo que se puede hacer considerando un número de varillas sucesivamente mayor, dando lugar a la completa generalización del problema, tanto en número de discos como en número de varillas (Benito González y Madonna, p. 244).

Problema 2.3.3.

Teniendo en cuenta las mismas reglas, si en vez de 3 varillas, hubiera 4 y utilizando 4 discos (Fig. 18), ¿cuál es el menor número de movimientos necesarios para trasladarlos a una de las 3 varillas libres? ¿Y si hay k≥1 varillas en lugar de 3 o 4 y n discos?

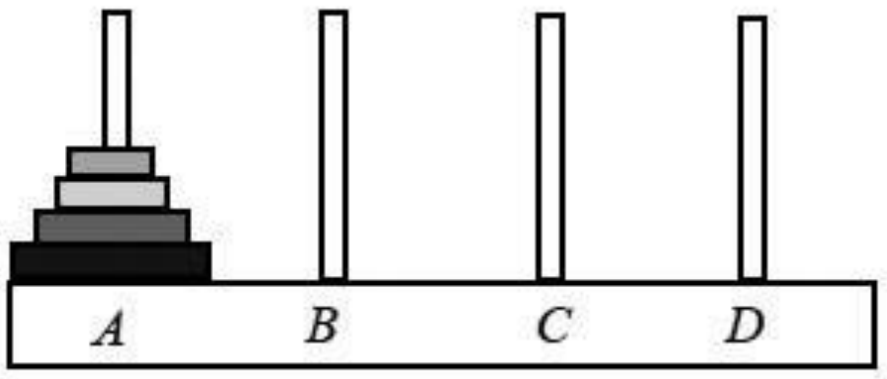

**Figura 18**

**Propuesta de investigación:**

Consistente en dos problemas que darán lugar a sendas investigaciones con un carácter similar. Para realizarla se agruparán en parejas.

Problema 2.3.4.

En relación al concepto de recursividad, se aprovechará para introducirles en la geometría fractal (que puede ser tratado en una sesión sucesiva). Se propone generar el **Triángulo y la Alfombra de Sierpinski** (Fig. 19), Unas parejas puede trabajar sobre el triángulo de Sierpinski y otras sobre la alfombra. Se pide:

- Proponer una definición por recurrencia para el número de triángulos, respectivamente cuadrados, negros según los casos sucesivos $n \geq 1$, siendo n=1 el caso inicial de 3 triángulos negros y 1 blanco, repectivamente 8 cuadrados negros y 1 blanco.
- Conjeturar una expresión explícita para cada caso y demostrarla por inducción.

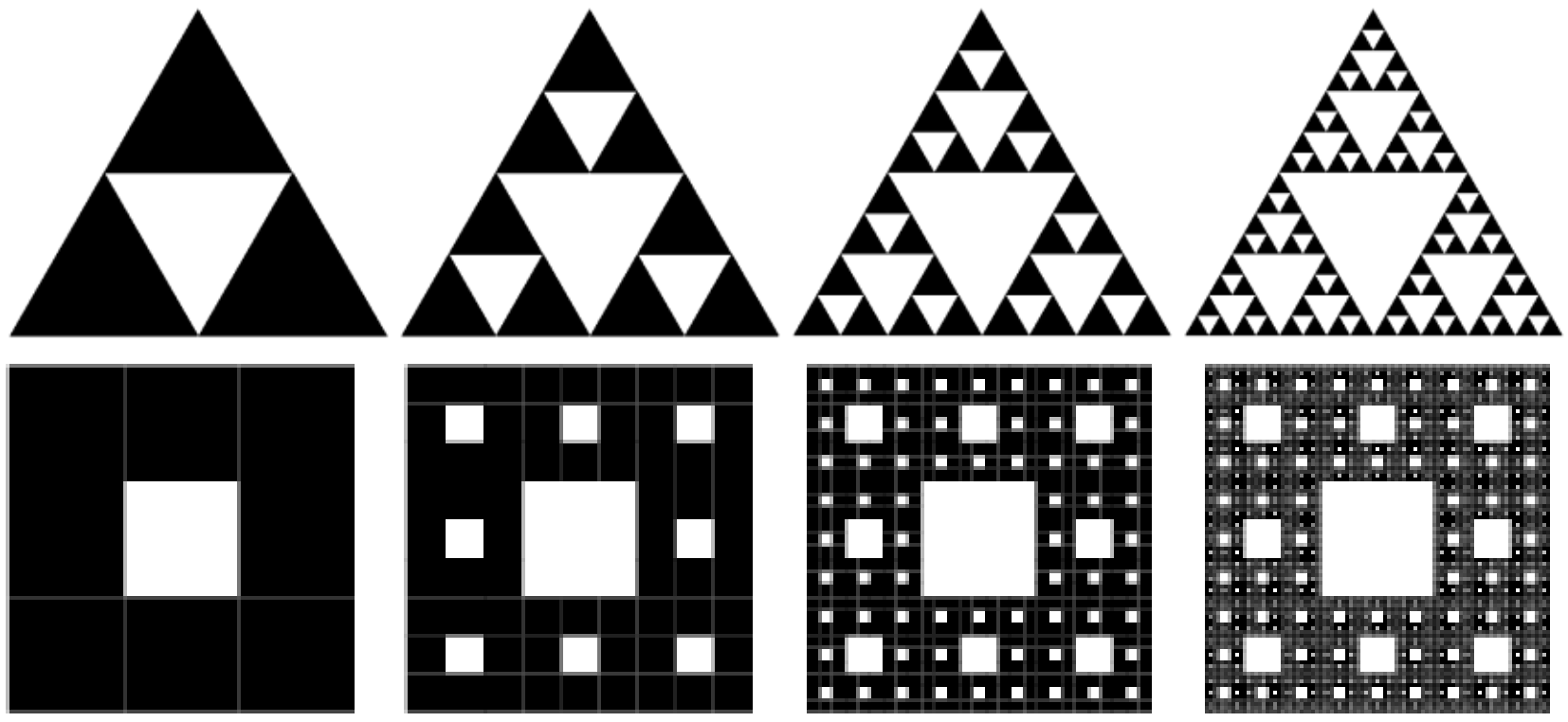

**Figura 19**

Sol.

Triángulo de Sierpinski (número de triángulos negros):

La sucesión recurrente queda definida por: $\begin{cases} c_1 = 3 \\ c_n = 3c_{n-1} \end{cases}$

Una expresión explícita para ella viene dada por: $c_n = 3^n$

Alfombra de Sierpinski (número de triángulos negros):

La sucesión recurrente queda definida por: $\begin{cases} c_1 = 8 \\ c_n = 8c_{n-1} \end{cases}$

Una expresión explícita para ella viene dada por: $c_n = 8^n$

Problema 2.3.5.

Dada la siguiente sucesión de figuras que se muestra en la Fig. 20, y sea $c_n$ el número de triángulos que aparecen en la n-ésima figura. Se pide:

- Averiguar si $c_n$ está definido por recurrencia
- Conjeturar una expresión explícita para $c_n$
- Demostrarla haciendo uso del principio de inducción

**Figura 20**

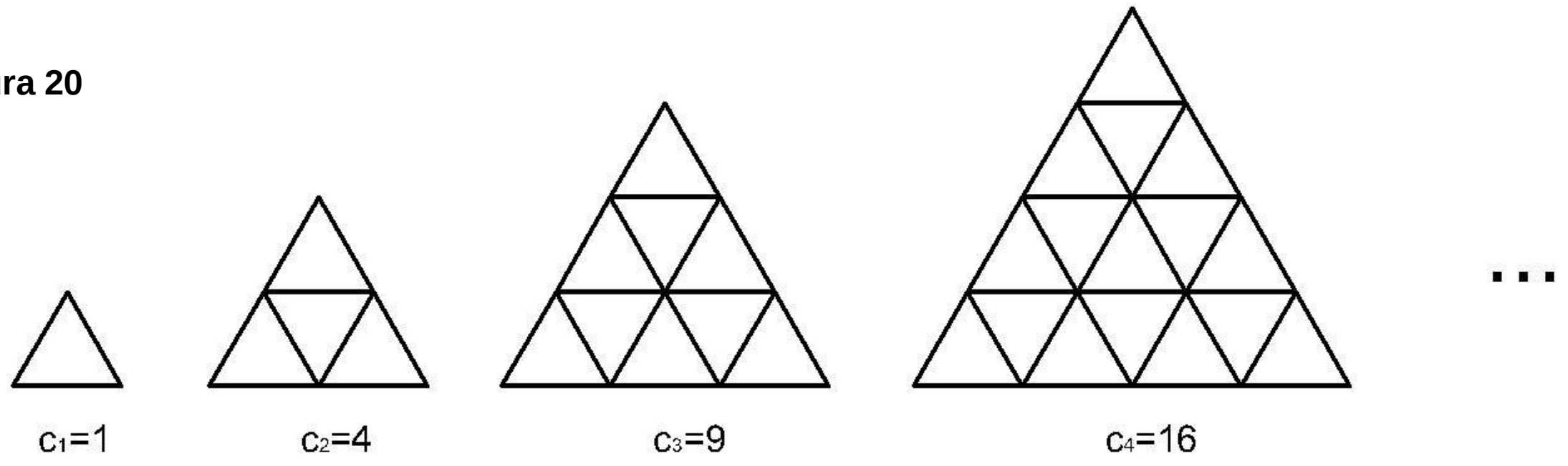

Se puede ofrecer como pista sombrear algunos triángulos, como muestra la Fig. 21:

**Figura 21**

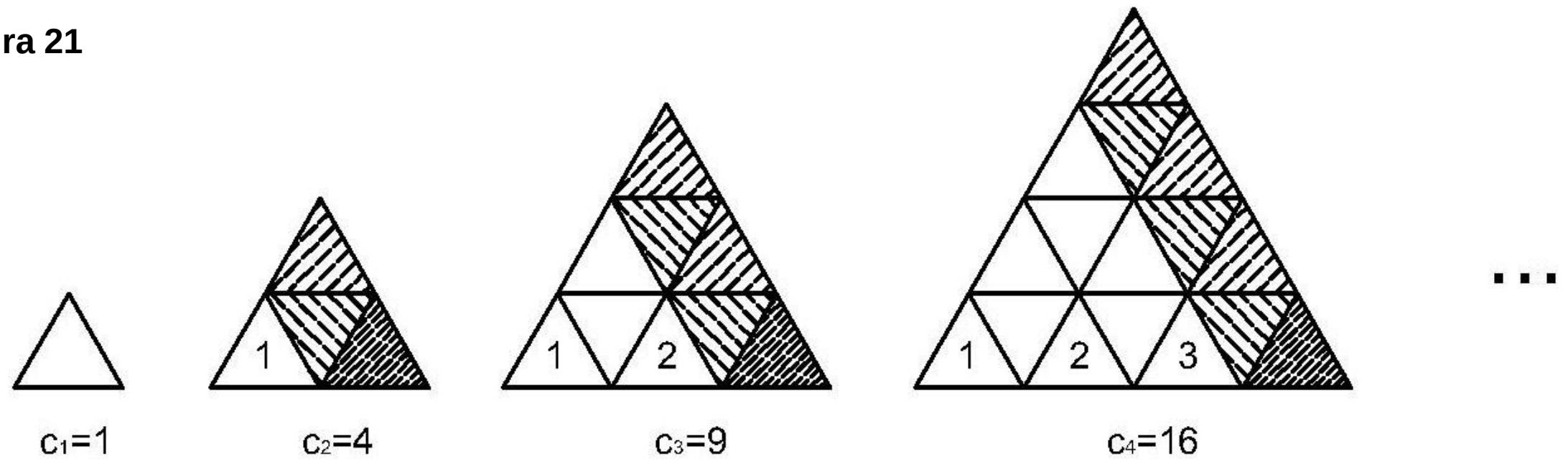

Sol.

La sucesión recurrente queda definida por: $\begin{cases} c_1 = 1 \\ c_n = c_{n-1} + 2n - 1 \end{cases}$

Una expresión explícita para ella viene dada por: $c_n = n^2$

**Objetivos:**

- Reconocer patrones
- Recursividad: definición por recurrencia, sucesión recurrente
- Conjeturar una expresión explícita a partir de la observación de casos sucesivos (inducción matemática) y se capaz de demostrarla
- Plantear problemas a partir de otros problemas
- Optimización: hallar un mínimo.
- Introducción a las geometrías fractales

## A.2.4. Demostración directa: áreas de polígonos (SI1).

### Desarrollo

La sesión está diseñada para el grupo de 1º de ESO y se requiere haber introducido antes los siguientes contenidos:

- Rectas y ángulos
- Paralelismo y perpendicularidad
- Clasificación de triángulos y cuadriláteros
- Polígonos simples y sus elementos
- Definición de circunferencia y círculo
- Polígonos regulares inscritos y circunscritos a una circunferencia

Se partirá del establecimiento de un axioma considerando el área de un cuadrado como la longitud de su lado elevado al cuadrado, $a^2$. Mediante demostraciones directas se irá estableciendo el área de otros polígonos: rectángulo, triángulo, paralelogramo y polígonos regulares con número de lados n≥5.

Para ello se comienza la sesión, focalizando en el método de demostración directa:

> "Este método utiliza las leyes transitivas, o silogismo hipotético. Para demostrar que el antecedente es condición suficiente para asegurar la verdad del consecuente, se busca una condición intermedia tal que el antecedente sea condición suficiente de ésta, y que ésta sea condición suficiente del consecuente." Se basa pues en la implicación $(p \Rightarrow q \Rightarrow r)$ o en la doble implicación $(p \Leftrightarrow q \Leftrightarrow r)$ (Delgado Pineda y Muñoz Bouzo, 2010, p. 28).

Se comentará el significado de axioma, como aquella 'evidencia' mínima de la que se parte, y que se establecen como verdadera, a partir de la cual se pueden ir construyendo otras proposiciones, haciendo uso del método de demostración directa. "Un axioma es una sentencia bien formada que se considera verdadera, es decir, una tautología primaria no deducible" (Delgado Pineda y Muñoz Bouzo, 2010, p. 30).

Problema 2.4.1.

Sea un rectángulo, de lados b y h, hallar y demostrar una fórmula para su área (R), conocida la del área del cuadrado, $A_{cuadrado} = a^2$, que ha sido establecida axiomáticamente (Sally y Sally, 2011).

Sol.:

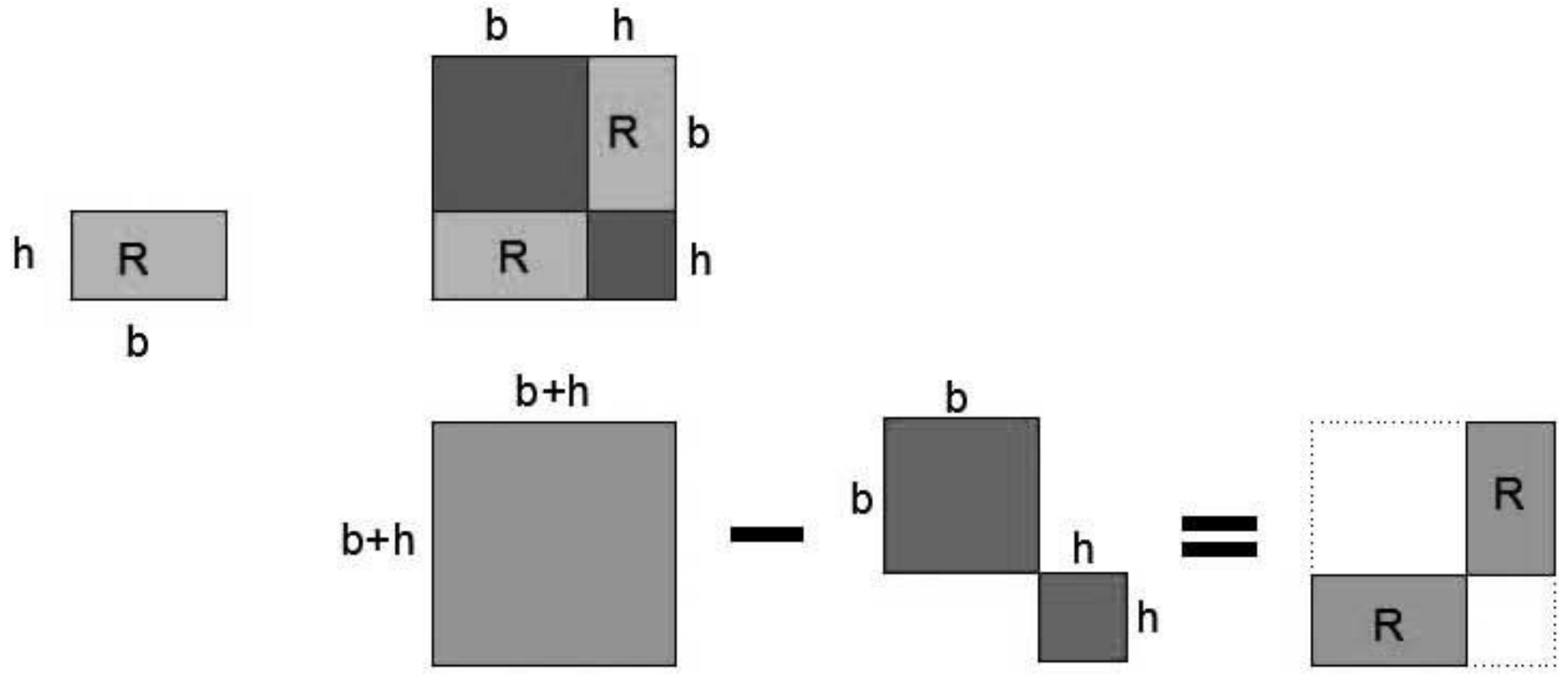

$$(b+h)^2 - (b^2 + h^2) = 2 \cdot R$$

$$b^2 + 2bh + h^2 - b^2 - h^2 = 2 \cdot R$$

$$2bh = 2 \cdot R$$

$$bh = R$$

**Figura 22**

Se plantearán otras opciones alternativas, de modo que se entienda que la demostración no tiene necesariamente porque ser única, como se muestra en la Fig. 23.

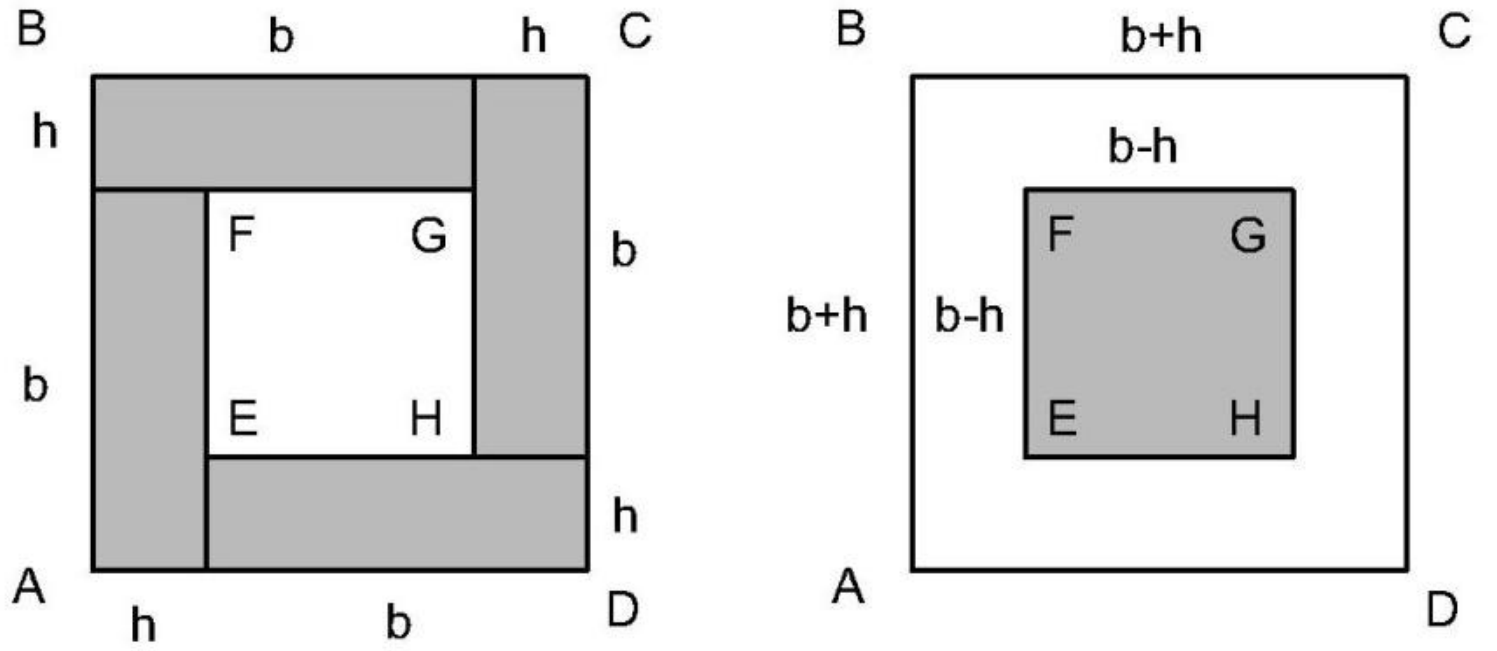

Sol.:

$A_{rectángulo} = \frac{1}{4}[(b+h)^2 - (b-h)^2] = \frac{1}{4}(b^2 + h^2 + 2bh - b^2 - h^2 + 2bh) = \frac{1}{4}4bh = bh$

**Figura 23**

Después de esta introducción se pedirá a los participantes realizar otras demostraciones, para calcular el área de otros polígonos característicos.

<u>Problema 2.4.2.</u>

Hallar y demostrar una fórmula para el área de un triángulo de base b y altura h, a partir de las ya conocidas y demostradas. Para ella habrá que tener en cuenta los posibles subcasos: (1) el de un triángulo rectángulo, (2) el de un triángulo acutángulo, (3) el de un triángulo obtusángulo.

Sol.:

Es sencillo obtener el área de un triángulo rectángulo, una vez se ha demostrado como hallar área del rectángulo, pues su area resulta ser exactamente su mitad. No lo es tanto en el caso de que se trate de un triángulo acutángulo u obtusángulo, donde habrá que buscar alguna estrategia algo más imaginativa, invitando a los alumnos a que la descubran. Una vez hayan conseguido dar con una demostración, se pueden realizar las siguiente pregunta: ¿Puede haber diferentes formas de demostrarlo? Daremos lugar a que esto suceda, por ejemplo, como se muestra en la Fig. 24 y en la Fig. 25

$$A_{triángulo} = \frac{1}{2}b_1h - \frac{1}{2}b_2h - = \frac{1}{2}h(b_1 - b_2) = \frac{1}{2}bh$$

$$A_{triángulo} = \frac{1}{2}x_1y + \frac{1}{2}x_2y - = \frac{1}{2}y(x_1 + x_2) = \frac{1}{2}xy$$

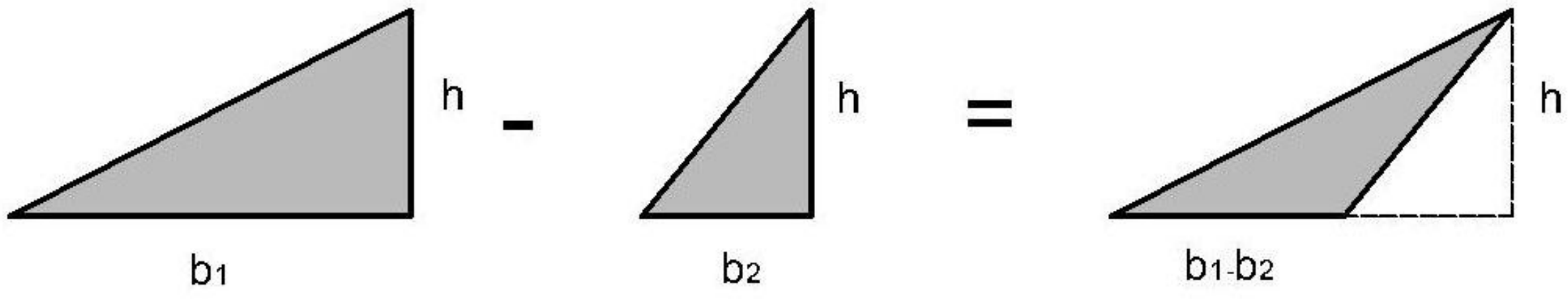

**Figura 24**

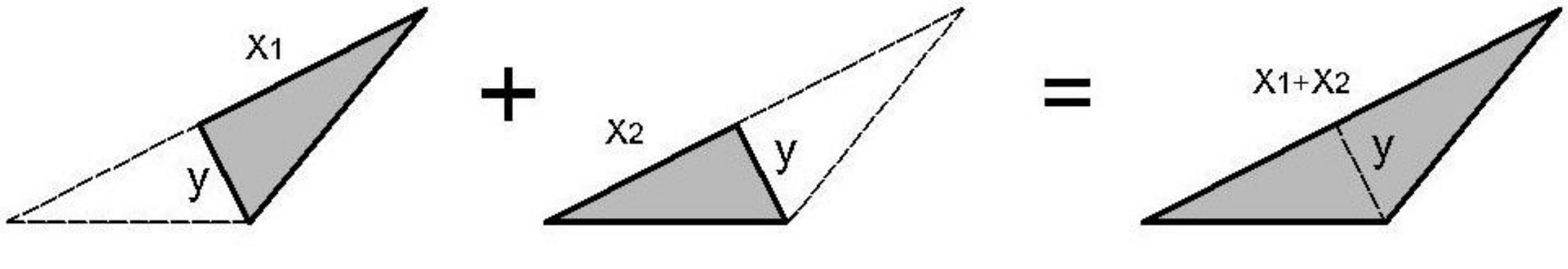

**Figura 25**

Problema 2.4.3.

Hallar y demostrar una fórmula general para el área de un paralelogramo romboide de lado mayor b y altura h, a partir de las ya conocidas y demostradas.

Se puede preguntar: ¿y si demostramos primeramente el área de un paralelogramo para posteriormente demostrar la de un triángulo? ¿Qué es mejor? Todas estas cuestiones, se irán formulando en la medida en que vayan avanzando en sus investigaciones, estableciendo un ritmo apropiado.

Sol.:

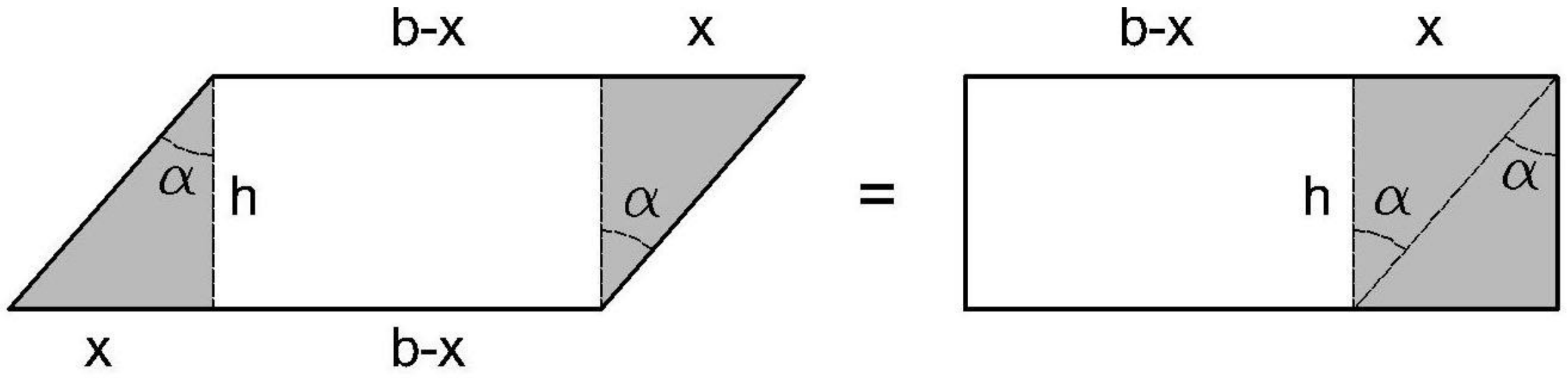

**Figura 26**

$A_{paralelogramo} = (b - x)h + xh = bh$

Problema 2.4.4.

Hallar y demostrar una fórmula general para el área de los poligónos regulares, de lado l y apotema a, a partir de las ya conocidas y demostradas.

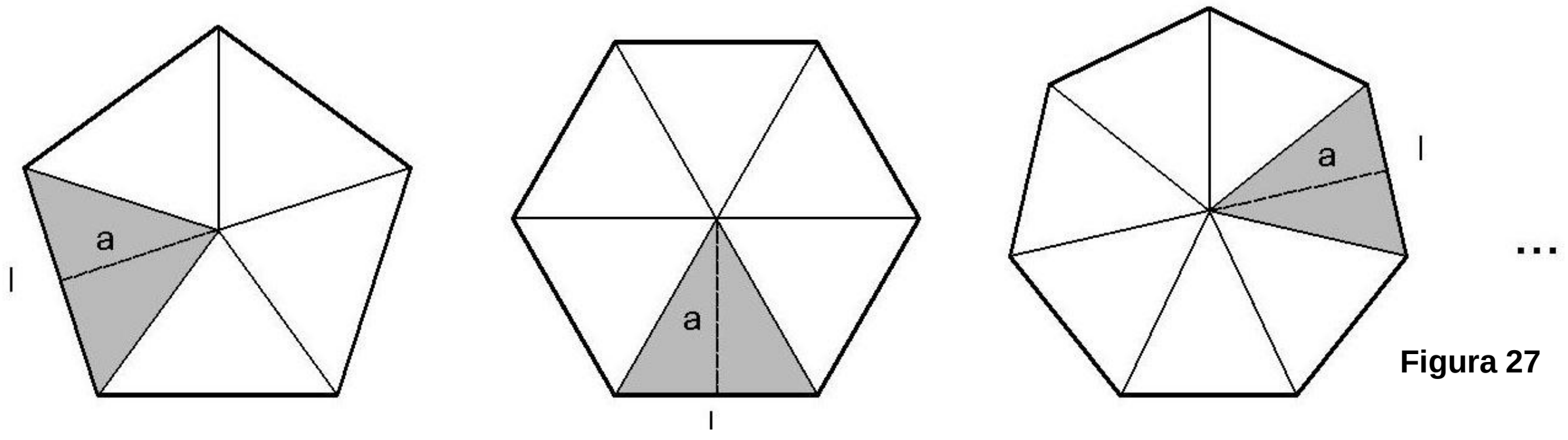

**Figura 27**

Sol.:

$A_{polígono\ regular} = n\frac{1}{2}l \cdot a; siendo\ n = nº\ de\ lados\ del\ polígono\ regular$

**Propuesta de investigación**

Se proponen dos problemas, uno sencillo, basado en la investigación desarrollada en la clase-sesión tipo Círculo, y otra más compleja como problema avanzado y con apoyo del profesor. El tiempo estimado es de 2 semanas.

Problema 2.4.6.

Dado un trapecio escaleno de altura h, base menor b y base mayor B, hallar y demostrar una fórmula general para su área, a partir de las ya conocidas y demostradas.

Sol.:

$A_{trapecio} = \frac{1}{2}(B + b)h$

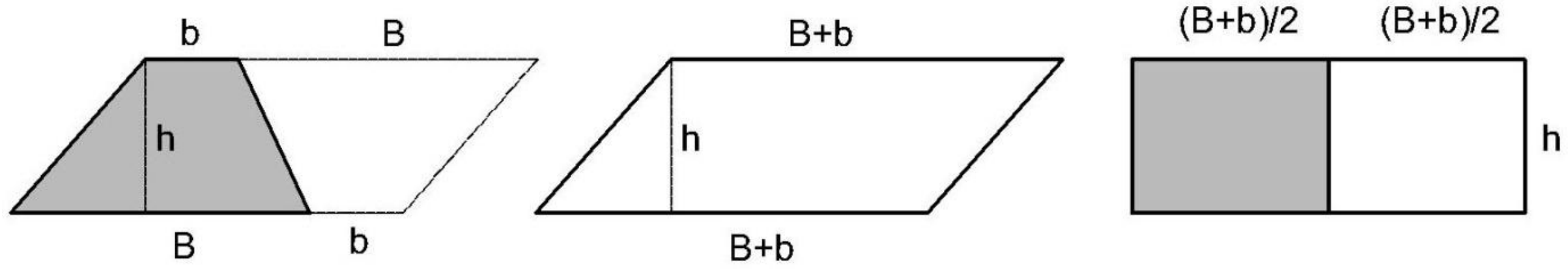

**Figura 28**

Problema 2.4.7.

Calcular el área de un círculo de radio r, a partir de observar lo que sucede con el aumento progresivo de los lados de un polígono regular inscrito en una circunferencia (Fig. 29).

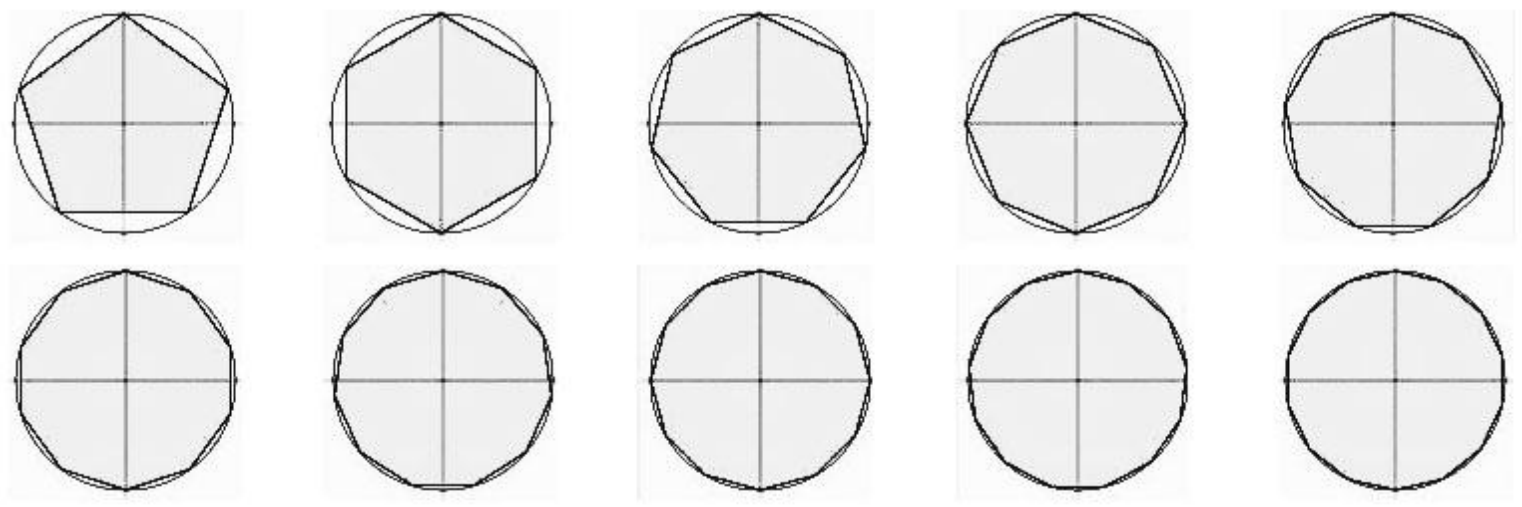

**Figura 29**

Se puede ofrecer como pista, que se calcule el área de un polígono regular de 10 lados inscrito en una circunferencia de radio 10 und., y después el de un polígono regular de 30 lados inscrito en esa misma circunferencia. Se utilizará geogebra, para facilitar la representación y los cálculos. Una vez se haya realizado se preguntará: ¿qué sucede a medica que aumentamos el número de lados?

> "Suele presentarse el círculo como un no polígono (porque no tiene lados rectos). Esto no es sino consecuencia de una visión tradicional y estática de la geometría. Desde una perspectiva dinámica, es fácil ver y comprobar cómo un polígono regular de 30 ó 40 lados, inscrito en un círculo, apenas puede diferenciarse del mismo. ¿Y si aumentamos el número de lados a 200 ó 1000? ¿Qué tendencia muestra su apotema? ¿Y la longitud de sus lados? No resulta chocante, pues, aceptar que un círculo es un polígono regular de infinitos lados rectos infinitamente pequeños. En el caso límite (al aumentar progresivamente el número de lados) la apotema se confunde con el radio del círculo, la longitud del lado del polígono regular tiende a cero y el perímetro tiende, sin sobrepasarlo, al valor de la longitud de la circunferencia."
> (García Moreno, 2012)

De esta forma se está introduciendo al mismo tiempo el concepto de **infinito**, que será el número de lados del polígono-círculo, de **infinitesimal**, que será la longitud de cada uno de los lados de dicho polígono-círculo. El área del círculo se entiende como el **límite** de un polígono regular inscrito en un circunferencia, cuando el número de lados tiende a infinito.

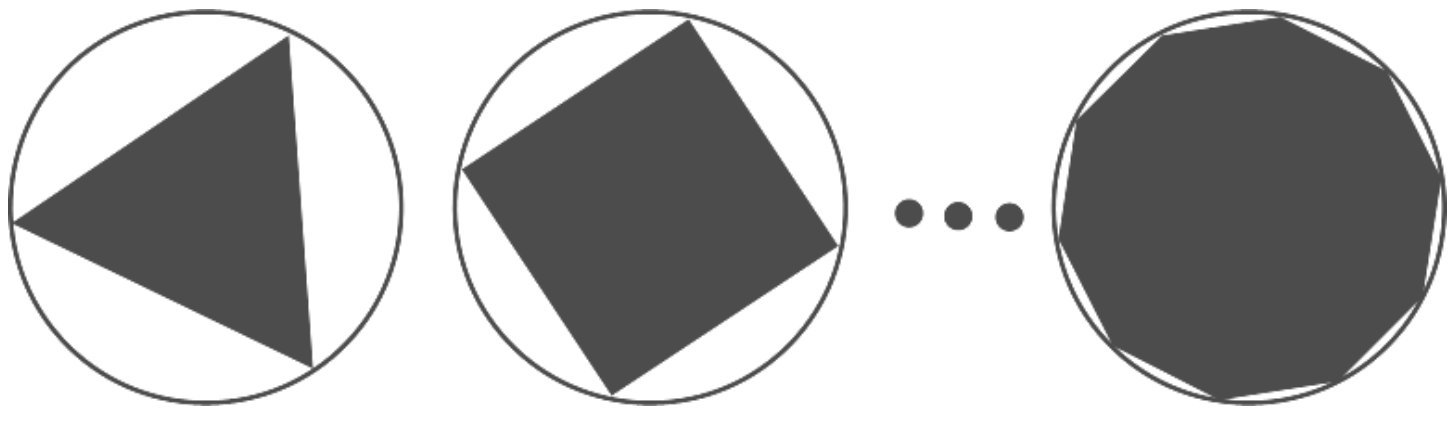

**Figura 30**

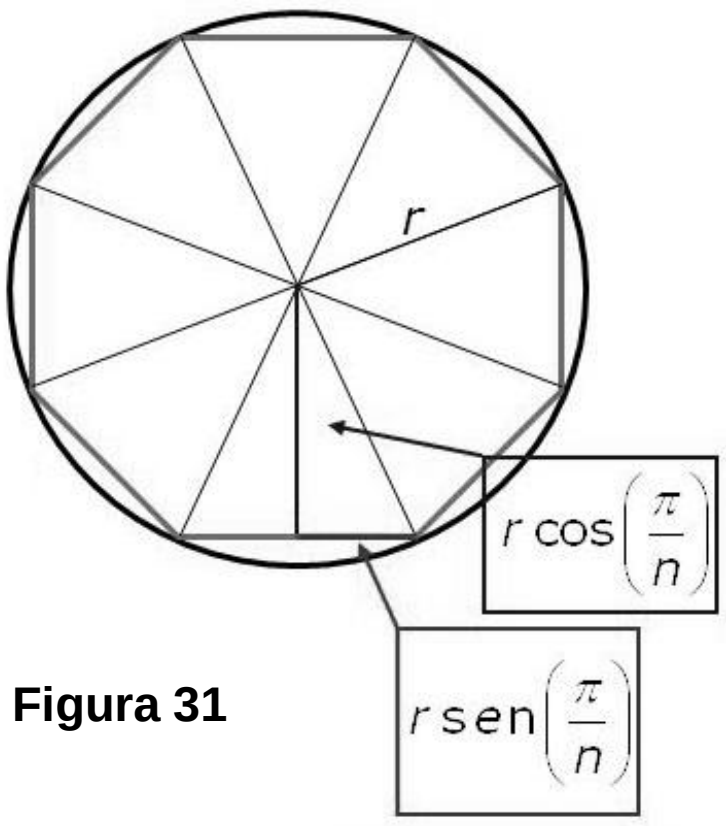

Figura 31

Es importante hacer una investigación para entender lo que sucede. Como en este nivel no se han introducido todavía conocimientos sobre trigonometría, bastará con observar que a medida que el número de lados aumenta, la longitud de la apotema del polígono se irá aproximando a la del radio. Construyendo un modelo (Fig. 31) en el que el área del polígono quede expresada en función del radio de la circunferencia en la que se encuentra inscrito, el área del círculo quedaría formulada con la siguiente expresión:

$\lim_{n\to\infty} nr^2 \sin\left(\frac{\pi}{n}\right) cos\left(\frac{\pi}{n}\right)$ ;

con cuyo cálculo se obtiene la conocida expresión del área del círculo[2] $A_c = \pi r^2$

| P(n) | n | Q(n) |
|---|---|---|
| $2\ r^2$ | 4 | $4\ r^2$ |
| $2.8284\ r^2$ | 8 | $3.3137\ r^2$ |
| $3.0614\ r^2$ | 16 | $3.1826\ r^2$ |
| $3.1214\ r^2$ | 32 | $3.1517\ r^2$ |
| $3.1363\ r^2$ | 64 | $3.1441\ r^2$ |
| $3.1405\ r^2$ | 128 | $3.1422\ r^2$ |
| $3.1401\ r^2$ | 256 | $3.1418\ r^2$ |
| $3.14157\ r^2$ | 1000 | $3.1416029\ r^2$ |
| $3.141592447\ r^2$ | 10000 | $3.141592757\ r^2$ |

Figura 32

Para realizar este proceso se requieren conocimientos avanzados, por eso se recomienda tan sólo proponerlo intuitivamente, y utilizar un método heurístico para el cálculo del área del círculo. Para ello se puede emplear el **método clásico de exhaución o agotamiento de Eudoxo (370 a.C.)**, y utilizar el programa informático **geogebra** como medio auxiliar para realizar los cálculos. El método de exhaución de Eudoxo, consiste en obtener dos sucesiones (Fig. 32) que se aproximan al área del círculo. Una sucesión viene determinada por las áreas sucesivas y crecientes, del polígono de *n* lados **inscrito** en una circunferencia, a medida que aumentamos el número de lados, y la otra por las áreas sucesivas y decrecientes para el polígono de *n* lados **circunscrito** a esa misma circunferencia a medida que aumenta también el número de lados (Fig. 33). De esta forma podemos establecer un valor muy preciso para el área de un círculo, que estará comprendido entre dos cotas, una superior y otra inferior.

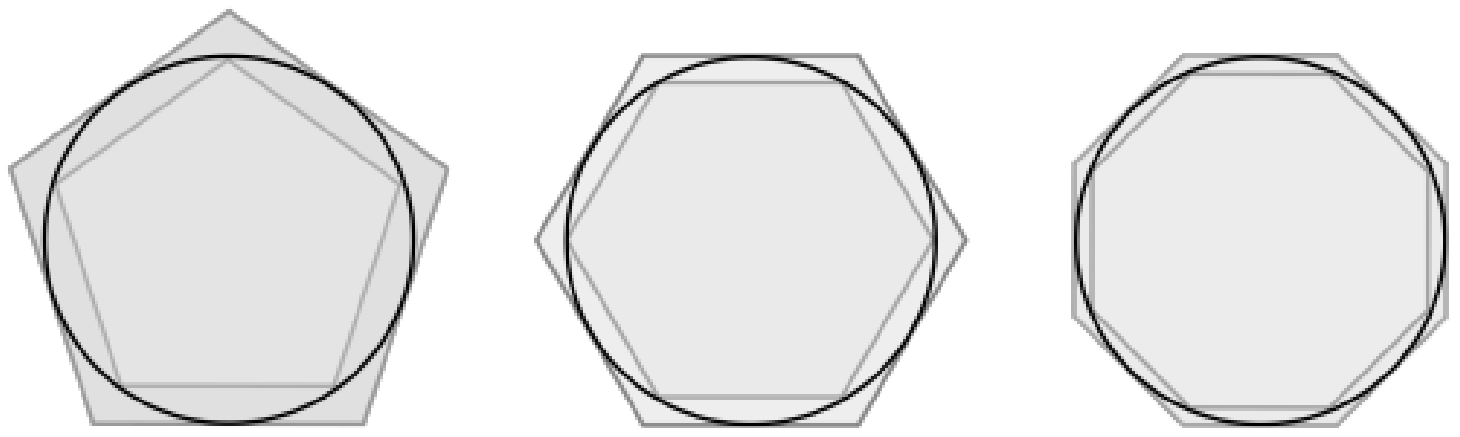

Figura 33

**Objetivos**

- Empleo del método de demostración directa
- Empleo del método heurístico
- Cálculo de áreas de polígonos: cuadrado, triángulo, rectángulo, paralelogramo, trapecio, polígonos regulares.
- Cálculo del área del círculo.
- Introducción de conceptos de cálculo infinitesimal: infinito, infinitesimal, límite

[2] Puede plantearse como problema avanzado para algunos alumnos

## A.2.5. Demostración directa: Teorema de pitágoras (SI1).[3]

### Desarrollo

La sesión está diseñada para el grupo de 2º de ESO y se requiere haber introducido antes los siguientes contenidos:

- Rectas y ángulos: intersección, paralelismo y perpendicularidad
- Polígonos simples y sus elementos
- Clasificación de triángulos y cuadriláteros
- Cálculo de áreas de polígonos simples

El objetivo es presentar el Teorema de Pitágoras de forma que constituya en sí mismo una investigación que vaya profundizando en los aspectos que vayan apareciendo. Se plantea como una exposición teórica fundamentada a la vez que se trabaja desarrollando un proceso de investigación. Se desarrollará según los siguientes puntos.

- Definición de un triángulo rectángulo y sus elementos característicos: catetos e hipotenusa.

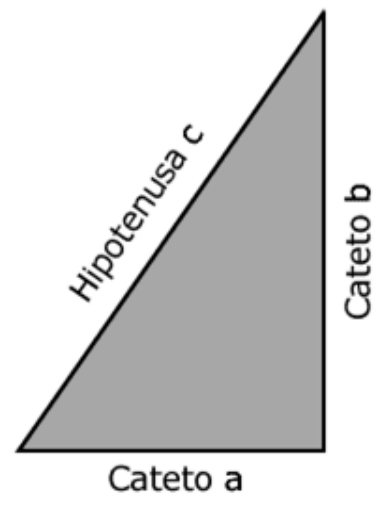

Figura 34

- Construcción geométrica del **Teorema de Pitágoras** a partir de un triángulo isósceles rectángulo, donde de forma fácil se puede ver la relación existente entre las áreas de los tres cuadrados al estar conformadas todas ellas por el triángulo isósceles y rectángulo inicial.

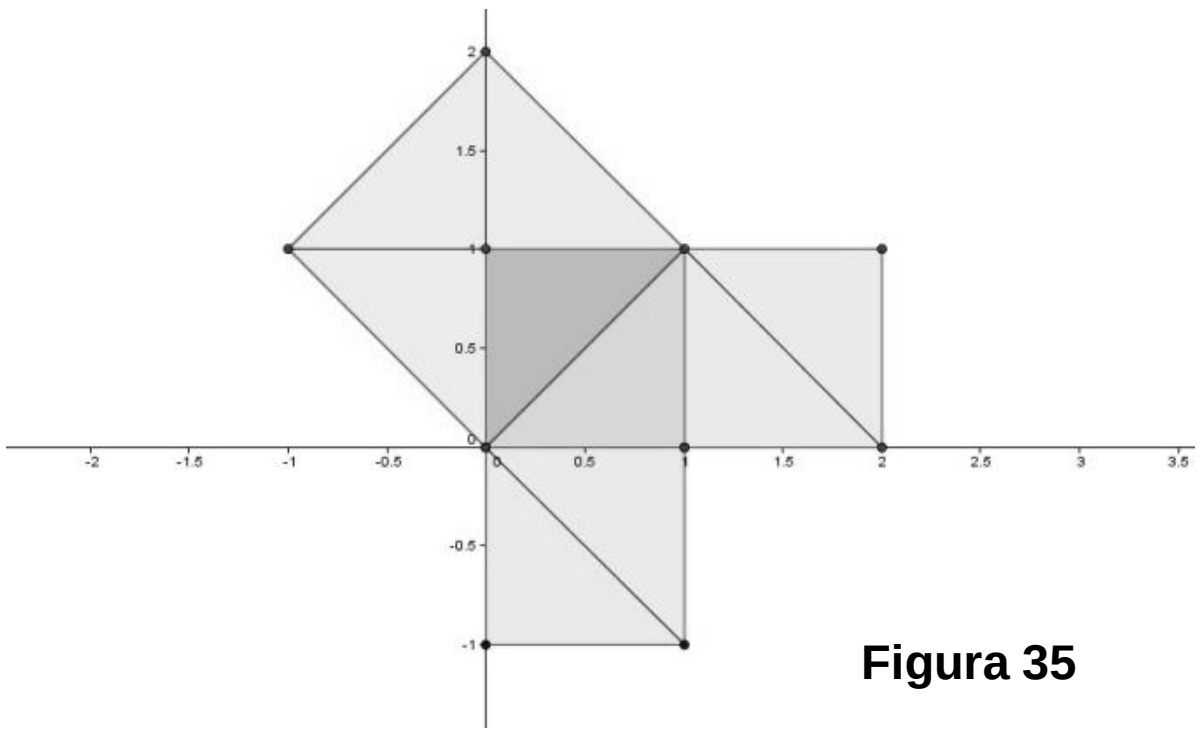

Figura 35

- Pregunta sobre si esa relación, entre las áreas de los cuadrados, se cumple para cualquier triángulo rectángulo.

- <u>Problema 2.5.1.</u>

  Investigación sobre por qué el área de un paralelogramo no varía si se desliza uno de sus lados (o los dos lados que sean paralelos) sobre la propia recta en la que los dos lados se hallan inmersos. La demostración sobre el área de un paralelogramo debería ser conocida por los alumnos.

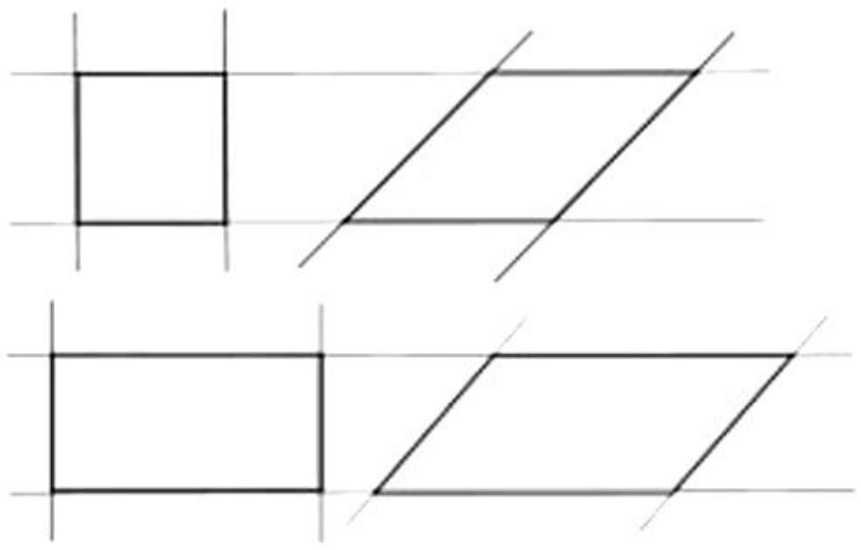

Figura 36

[3] Esta sesión tuve la oportunidad de llevarla a cabo con un grupo de 2º de ESO durante las prácticas del Máster correspondientes al módulo específico.

- Proyección de animación con la demostración de la **generalización del teorema de Pitágoras por parte de Pappus** (Legendre, 1849, p. 138). Descripción de la misma paso a paso.

  https://www.geogebratube.org/student/m224705

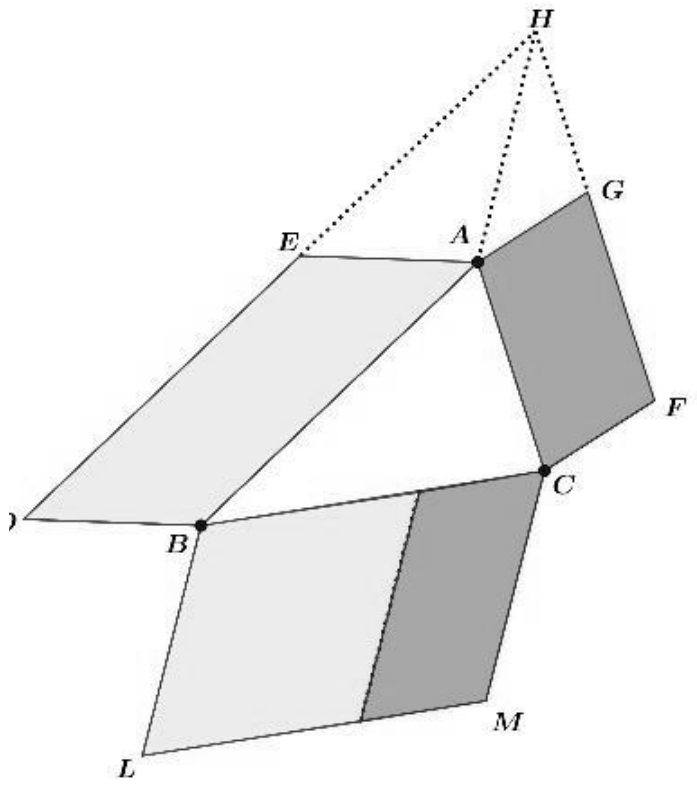

Figura 37

- Problema 2.5.2.

  Propuesta de tarea consistente en la aplicación de la demostración de Pappus al caso de un triángulo rectángulo y la construcción de dos cuadrados respectivamente sobre cada uno de sus catetos. Tendrán que dibujarlo y demostrarlo aplicando el razonamiento de la demostración anterior (Madonna, 2015, p. 251).

- Proyección de imágenes con los pasos a seguir para realizar la demostración.

- Proyección de animación gráfica con la demostración.

  https://www.youtube.com/watch?v=W6lBalZMCAI

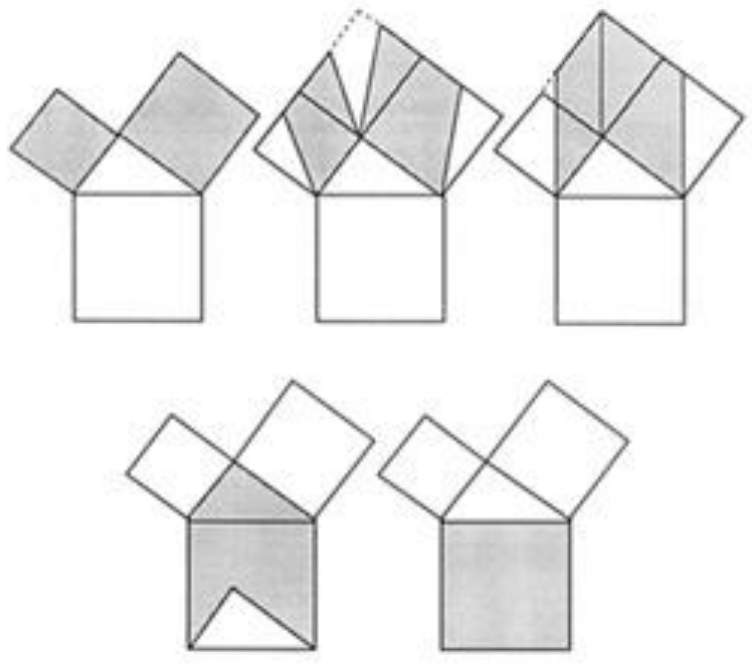
Figura 38

- Expresión aritmética del teorema, y mención a que conocida la longitud de dos de los tres lados de un triángulo rectángulo, podemos ,gracias al teorema, fácilmente conocer el tercero,

$$a^2 = b^2 + c^2 \iff \begin{matrix} a = \sqrt{b^2 + c^2} \\ b = \sqrt{a^2 - c^2} \\ c = \sqrt{a^2 - b^2} \end{matrix}$$

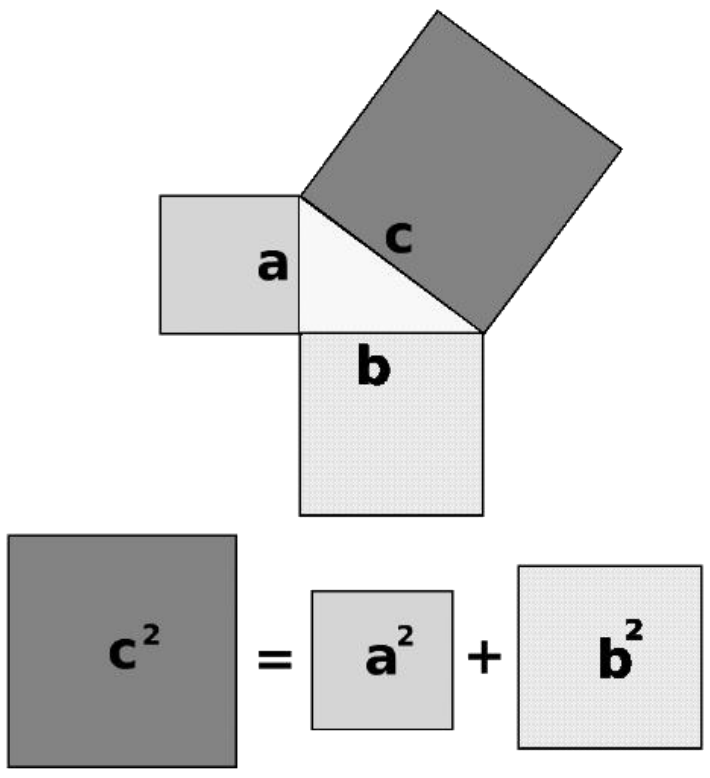

Figura 39

- Mención, a que el inverso del teorema también es cierto, es decir, si dicha relación se cumple en un triángulo, entonces el triángulo es rectángulo.

  $Si\ a^2 + b^2 = c^2\ \ entonces\ el\ triángulo\ es\ rectángulo$

- Problema 2.5.3.

  ¿Qué sucede en los casos en los que el Teorema de Pitágoras no se cumpa, es decir en el caso de que $a^2 + b^2 \neq c^2$?

  Sol.
  $a^2 + b^2 > c^2 \iff el\ triángulo\ es\ acutángulo$
  $a^2 + b^2 < c^2 \iff el\ triángulo\ es\ obtusángulo$

- Proyección de vídeo, con una demostración experimental, donde se ha construido un triángulo y tres recipientes cúbicos sobre sus lados. Los recipientes de los catetos se llenan de agua y al girar el mecanismo, todo el agua se deposita en el recipiente construido sobre la hipotenusa hasta llenarse por completo.

  https://www.youtube.com/watch?v=1er3cHAWwIM

- Problema 2.5.4.

  Resolver el puzle interactivo encontrado en la dirección:

  http://www.etudes.ru/en/etudes/pifagor/

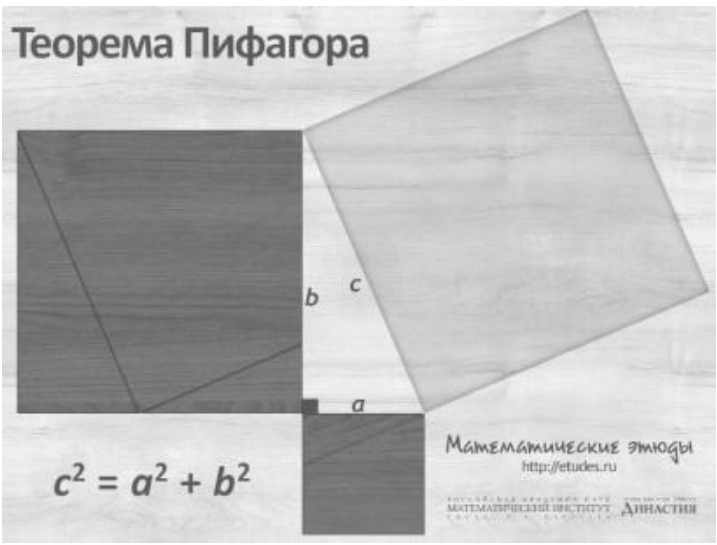

**Figura 40**

- Alusión, con refuerzo gráfico, de que el teorema también se cumple para la construcción de cualquier otra figura, manteniendo la semejanza, sobre cada uno de los lados del triángulo.

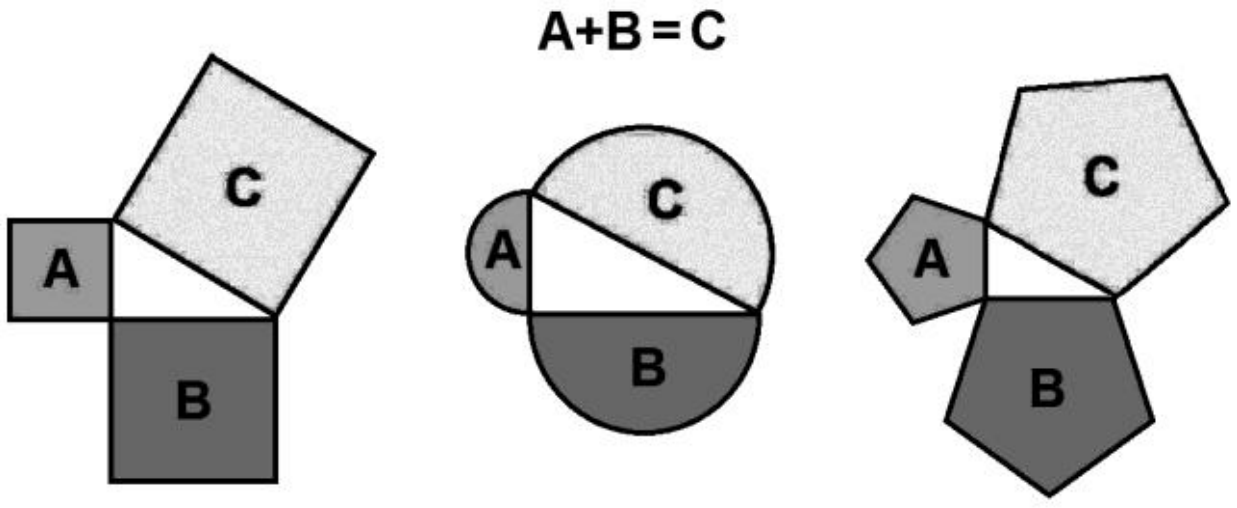

**Figura 41**

**Propuesta de investigación**

Se propone la resolución de dos problemas donde se tendrá que aplicar el teorema de Pitágoras, pero de forma que no se haga evidente, es decir, habrá que buscar estrategias para poder conducir el problema hasta una situación en la que se pueda aplicar el teorema. El tiempo estimado para su realización es de una semana.

- Problema 2.5.5.

  En la siguiente figura, halla la longitud del lado a, teniendo en cuenta que ABDE es un cuadrado.

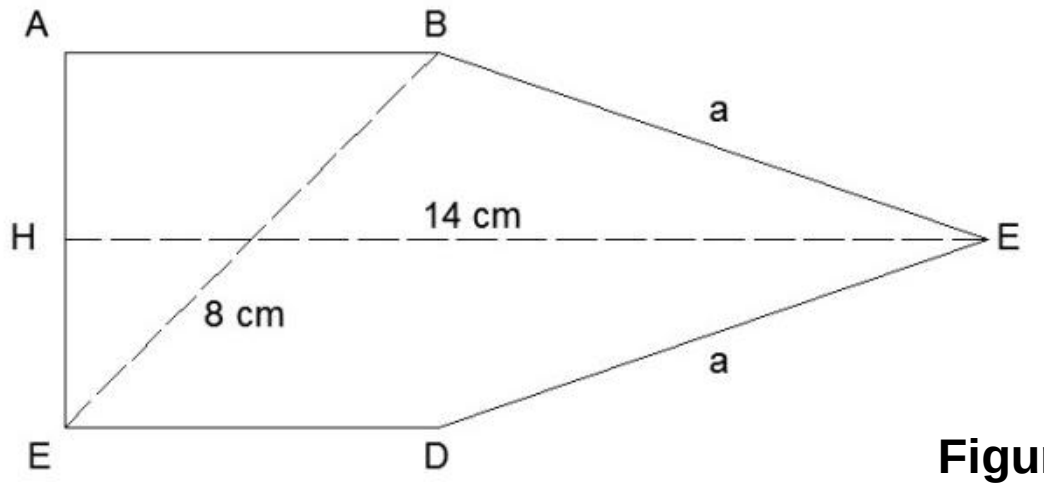

**Figura 42**

Sol.

Como ABDE es un cuadrado entonces $\overline{AB} = \overline{BD} = \overline{ED} = \overline{EA} = l$

$$8^2 = \overline{BD}^2 + \overline{ED}^2 \iff 64 = 2 \cdot l^2 \iff l^2 = 32 \iff l = \sqrt{32} \approx 5{,}657\text{cm}$$

$$\text{como} \quad \frac{\overline{BD}}{2} = \frac{\sqrt{32}}{2} = \frac{\sqrt{2^5}}{2} = \frac{2\sqrt{2^3}}{2} = \sqrt{8} \quad \text{y} \quad \overline{HE} - \sqrt{32} = 14 - \sqrt{32}$$

$$a^2 = \sqrt{8}^2 + (14 - \sqrt{32})^2 \iff a^2 = 8 + 196 - 28\sqrt{32} + 32 \iff a^2 = 236 - 112\sqrt{2} \iff$$
$$\iff a^2 \approx 77{,}61 \iff a \approx \sqrt{77{,}61} \approx 8{,}81$$

- Problema 2.5.6.

  Se rodea con un trozo de cuerda una barra circular. La cuerda da exactamente cuatro vueltas a la barra. La circunferencia de la barra es de 4 centímetros y su longitud, de 12 centímetros. Calcula la largura del trozo de cuerda. Muestra todos los pasos (ABC, 2015).

  Sol.

  En los extremos se forman dos triángulos rectángulos. Como sabemos que uno de los catetos mide 4 centímetros y el otro, 3 (una cuarta parte de la longitud de la barra), podemos hallar cuánto mide la hipotenusa. Este resultado nos basta para saber la longitud de la cuerda, ya que hay 4 trozos iguales alrededor de la barra.

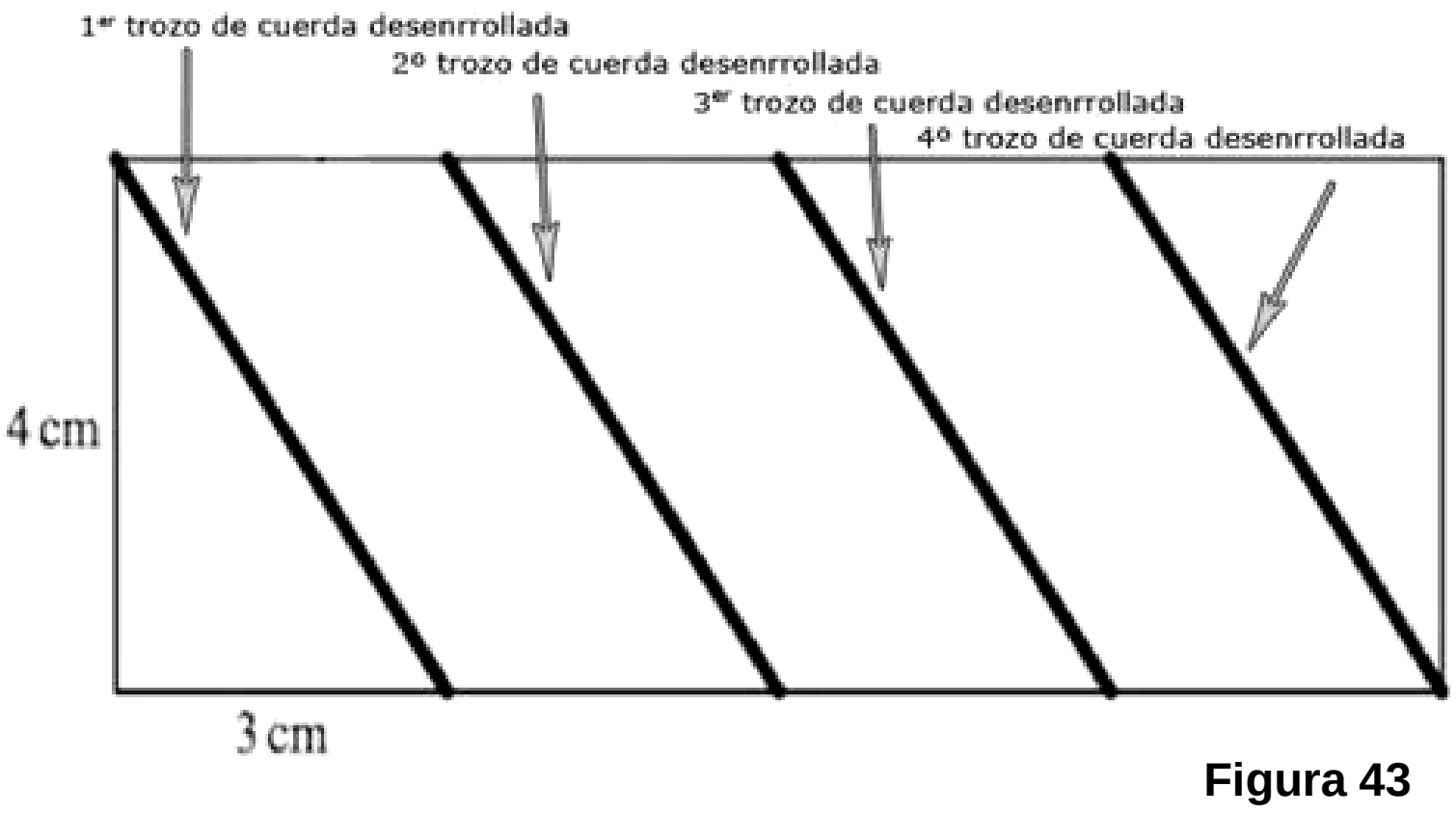

**Figura 43**

  Aplicando el Teorema de Pitágoras:

  $$4^2 + 3^2 = c^2 \Leftrightarrow \; c^2 = 25 \Leftrightarrow c = \sqrt{25} \Leftrightarrow c = 5 \; (c = -5 \;\; no\ resulta\ aplicable)$$

  Por lo tanto, 5 cm es la longitud de cada uno de los trozos de cuerda. Como hay 4 trozos iguales, la cuerda mide: 4·5 = 20 centímetros.

**Objetivos**

- Empleo del método de demostración directa
- Aplicar el Teorema de Pitágoras en la resolución de problemas
- Demostrar el Teorema de Pitágoras a partir de la generalización del teorema de Pitágoras por parte de Pappus.
- Averiguar la clasificación de un triángulo según sus ángulos utilizando el teorema de Pitágoras.
- Plantear problemas a partir de otros problemas.
- Expresar de forma algebraica las conclusiones geométricas halladas.

## A.3. Sesiones de problemas.

### A.3.1. Sesión de problemas Tipo SP1.

#### Introducción

Esta sesión de problemas se propone para el curso de 1º de ESO, y será desarrollada dentro de la Unidad Didáctica en la que se incluya el tema de Geometría Plana.

Se ha diseñado la sesión pensando en que los problemas no tienen que hacer referencia necesariamente a contenidos de la Unidad Didáctica que se hayan visto en sesiones anteriores, sino que pueden anticipar algunos posteriores e incluso optar por incluir alguno de mayor dificultad a los que se hayan venido realizando hasta el momento.

Se proponen 3 problemas bajo el tema indicado y 1 problema extra basado en algún tema correspondiente a Unidades Didácticas vistas anteriormente, a modo de refuerzo, o de la siguiente Unidad Didáctica prevista, con carácter introductorio de la siguiente hoja de problemas. En este caso se ha optado por un problema de Aritmética (Enteros) formulado como un problema de périmetros, en consonancia con el tema propuesto.

De manera adicional, se dispondrá de 4 problemas alternativos, también de diferente dificultad y que puedan ser utilizados de diversos modos:

- para sustituir o facilitar el camino hacia la resolución de alguno de los propuestos durante la sesión
- para quienes hayan terminado con los problemas propuestos
- como trabajo adicional para casa para quién los solicite

Los problemas no se recogerán, sino que se podrá seguir trabajando en ellos hasta que sean comentados finalmente en una sesión posterior. En la presente sesión se indicará que problemas se tiene previsto que se comenten en la próxima.

#### Objetivos

| | | |
|---|---|---|
| Problema 3.1.1. | Áreas, congruencia, divisibilidad | Convertir un problema en otro, visión espacial |
| Problema 3.1.2. | Polígonos regulares | Búsqueda de regularidades que impliquen otras regularidades |
| Problema 3.1.3. | Ángulos | Encontrar estrategia para aplicar conocimiento sobre ángulos |
| Problema 3.1.4. | Perímetro, enteros, divisibilidad | Convertir un problema en otro, propiedad de los números enteros, visión espacial |
| Problema 3.1.5. | Áreas, aritmética | Convertir un problema en otro, números impares, potencias, visión espacial |
| Problema 3.1.6. | Polígonos regulares | Búsqueda de regularidades que impliquen otras regularidades |
| Problema 3.1.7. | Ángulos | Encontrar estrategia para aplicar conocimiento sobre ángulos |
| Problema 3.1.8. | Paridad | Demostración por reducción a lo absurdo |

## Problemas principales

Problema 3.1.1.

Un propietario poseía un terreno ABCD con la forma exacta de un cuadrado. Vendió ¼ del mismo, y ese ¼ AGFE tenía también la forma de un cuadrado. La parte restante debía ser repartida en cuatro partes que fueran iguales en forma y tamaño. ¿Cómo resolver el problema? (Tahan, 2009, p. 50).

Pista: asociar un número entero a la parte restante

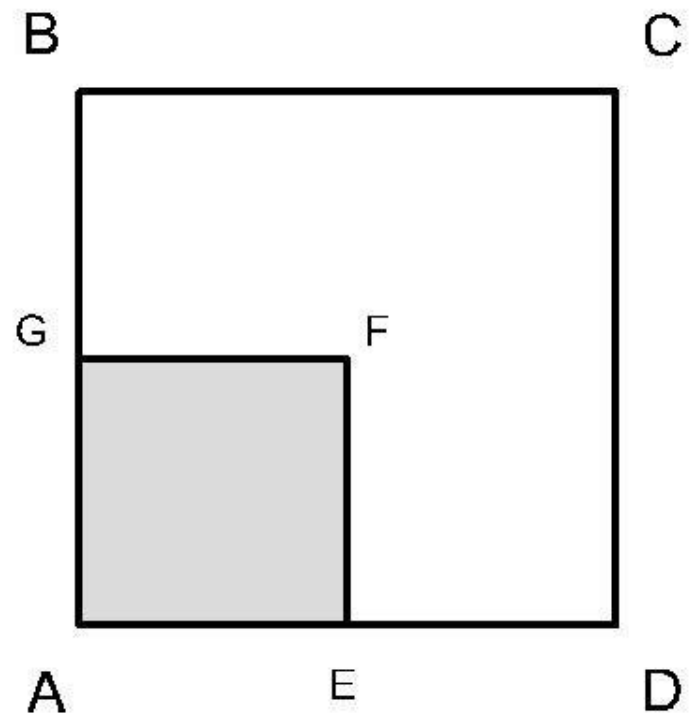

**Figura 44**

Problema 3.1.2.

Prueba en un octógono regular, cuyos vértices están denominados de forma consecutiva de 1 a 8, que uniendo mediante segmentos los vértices pares de forma consecutiva (o los impares) se construye un cuadrado (nota: probando solamente que los cuatro lados sean iguales no es suficiente para demostrar que se trate de un cuadrado) (Lehoczky y Rusczyk, 1993, p. 143).

Pista: identificar elementos y figuras congruentes

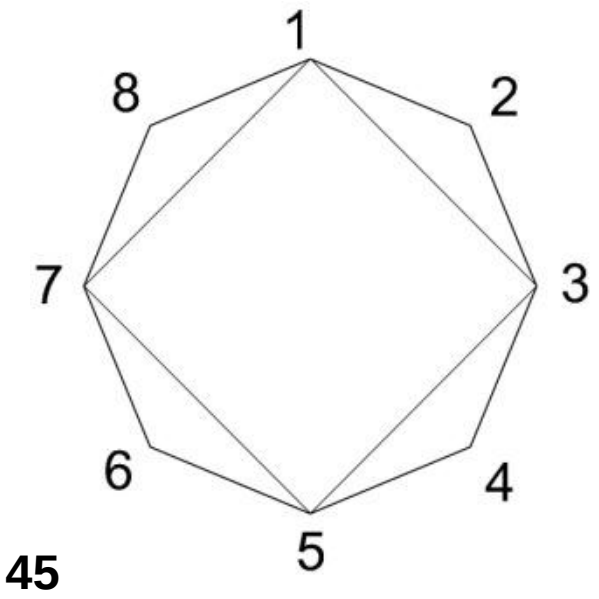

**Figura 45**

Problema 3.1.3.

Teniendo en cuenta que la recta *r* es paralela a la recta *s* y el valor de los ángulos que se indican en la figura, determina el valor de los ángulos Â, Ê y Û.

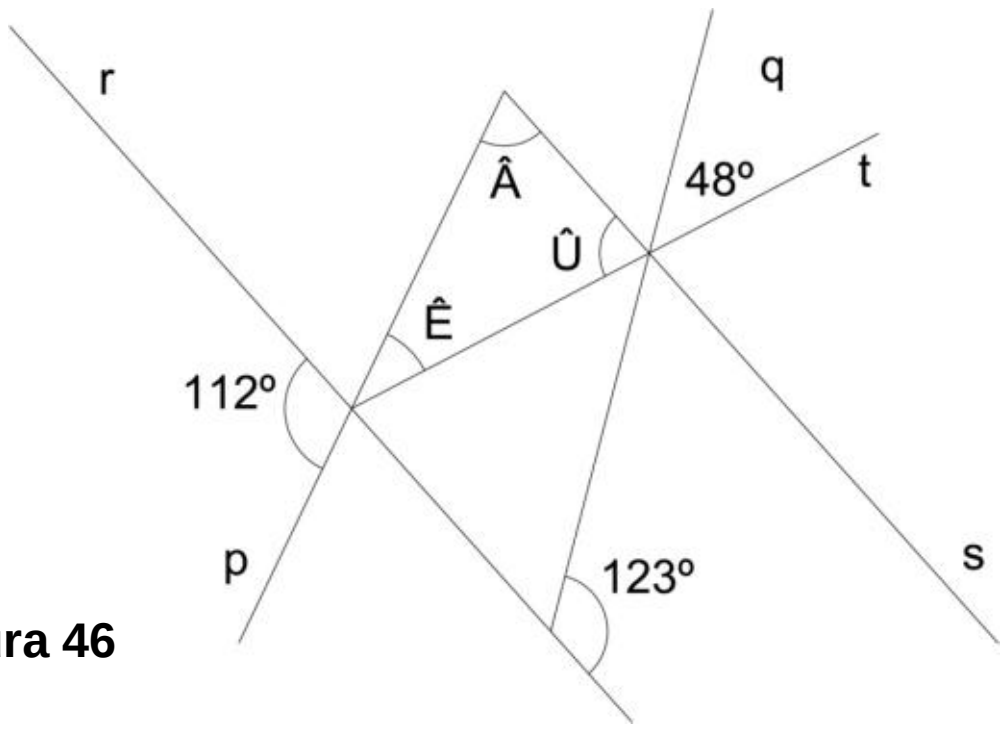

**Figura 46**

## Problema adicional:

Problema 3.1.4

Se ha dividido un rectángulo, mediante 6 líneas verticales y 6 horizontales, en 49 nuevos rectángulos (ver la figura). De esta división resulta que el perímetro de cada uno de los rectángulos resultantes es un número entero de metros. ¿El perímetro del rectángulo original tiene que ser necesariamente también un número entero de metros? (Yashchenko, 2010, p. 24).

Pista: sombrear los rectángulos de una diagonal (Yashchenko, 2010, p. 75).

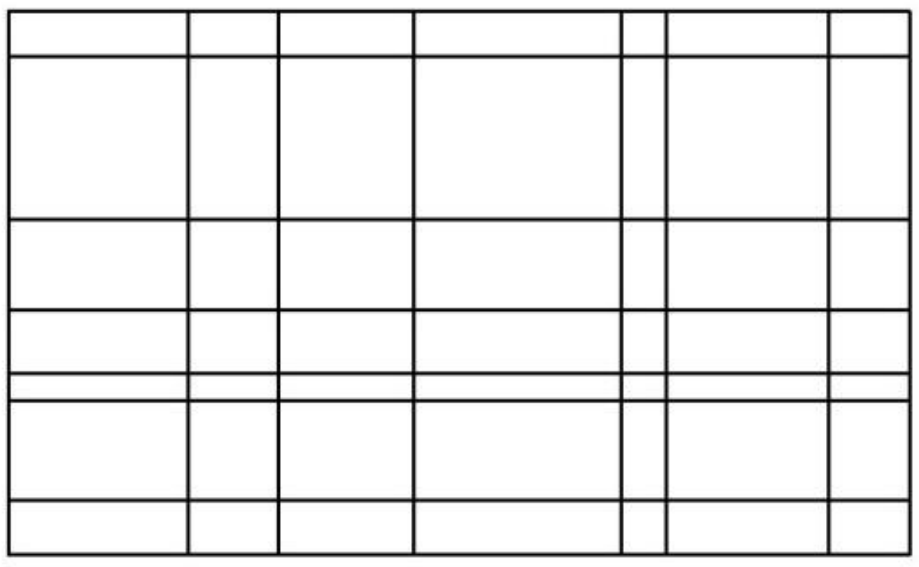

**Figura 47**

## Soluciones

<u>Problema 3.1.1.</u>

Trazando una cuadrícula que divida la parte restante en un número entero (n) y divisible por 3, el problema se reduce a buscar 4 figuras congruentes. En este caso n=12 por lo que 12/4=3 y resulta posible encontrar cuatro figuras congruentes con forma de trominó.

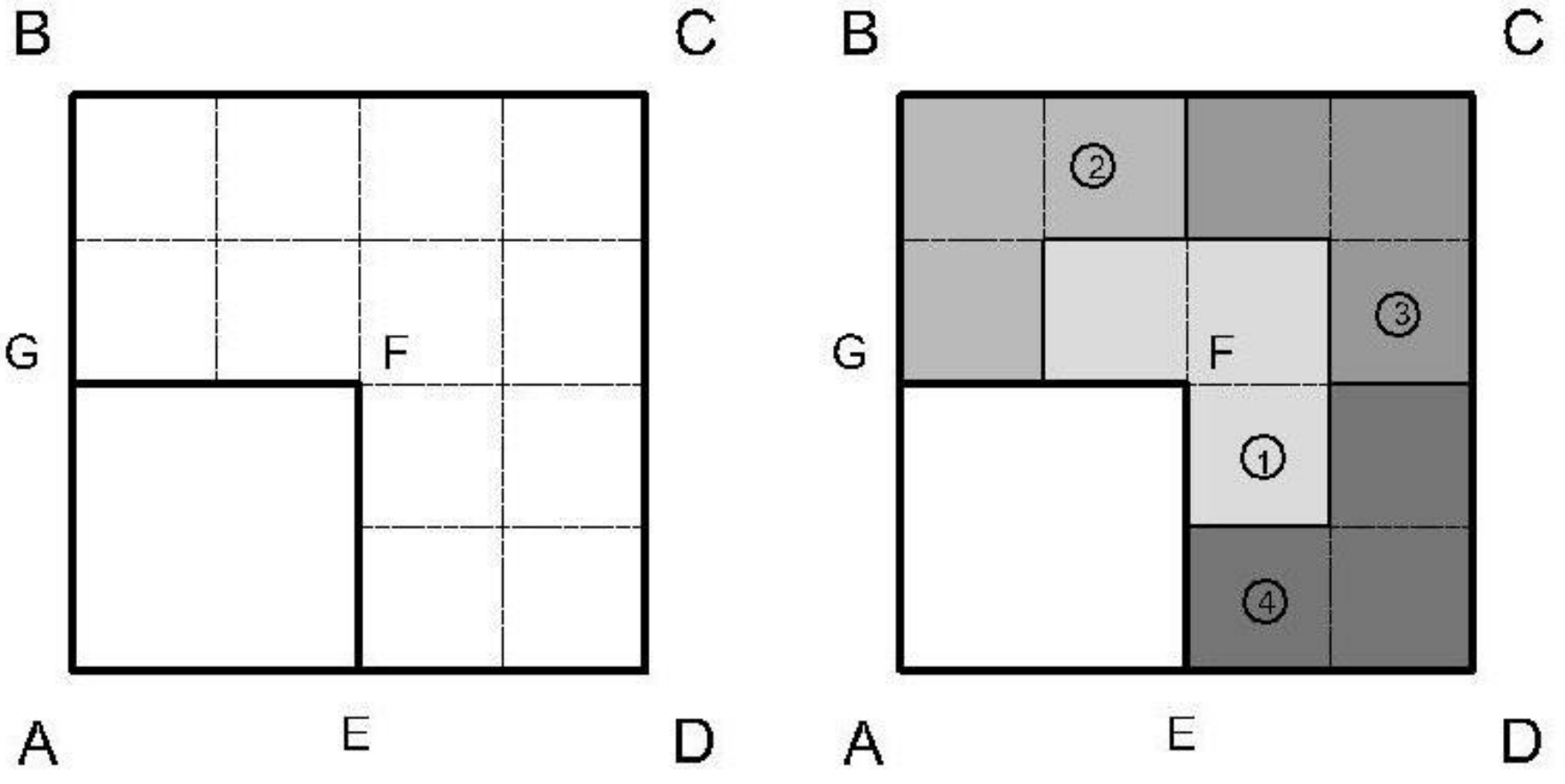

**Figura 48**

<u>Problema 3.1.2.</u>

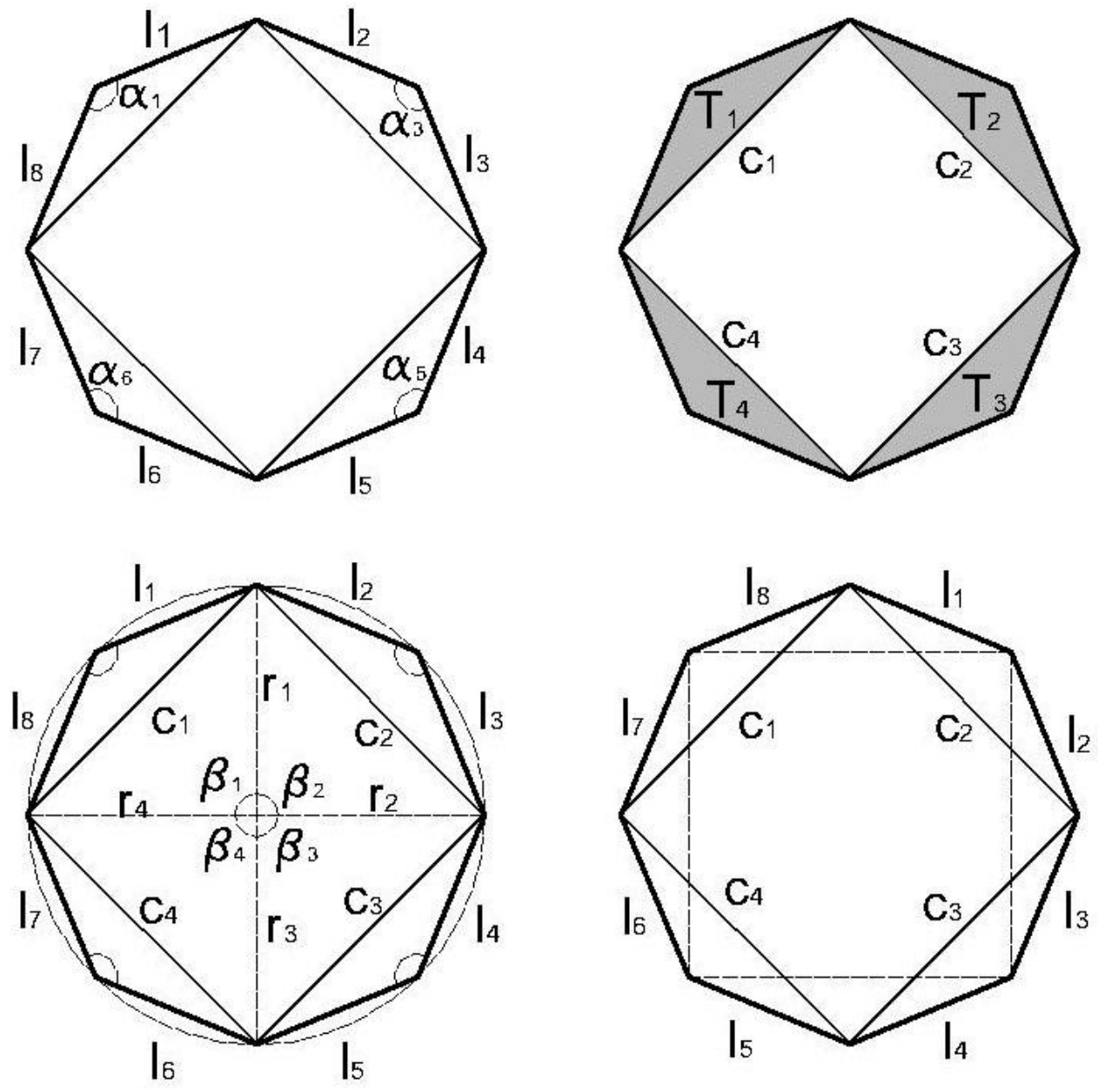

**Figura 49**

Los lados de un polígono regular son iguales. Lo mismo sucede con sus ángulos interiores. Luego:

$l_1 = l_2 = l_3 = l_4 = l_5 = l_6 = l_7 = l_8$ y

$\alpha_1 = \alpha_2 = \alpha_3 = \alpha_4$

Por ello, necesariamente:

$\mathrm{T}_1 = \mathrm{T}_2 = \mathrm{T}_3 = \mathrm{T}_4 \Longrightarrow \mathrm{c}_1 = \mathrm{c}_2 = \mathrm{c}_3 = \mathrm{c}_4$

Trazamos los radios que unen el centro de la circunferencia en la que se inscribe el polígono regular con los vértices del cuadrilátero. Luego:

$\mathrm{r}_1 = \mathrm{r}_2 = \mathrm{r}_3 = \mathrm{r}_4$ y $\mathrm{c}_1 = \mathrm{c}_2 = \mathrm{c}_3 = \mathrm{c}_4 \Longrightarrow$

$\beta_1 = \beta_2 = \beta_3 = \beta_4 = 90°$

Como las diagonales del cuadrilátero son iguales, perpendiculares y se cortan en su punto medio, necesariamente ha de tratarse de un cuadrado.

La otra posible situación (polígono trazado sobre los vértices pares), se puede transformar en la ya demostrada, girando la figura el ángulo central del polígono regular.

<u>Problema 3.1.3.</u>

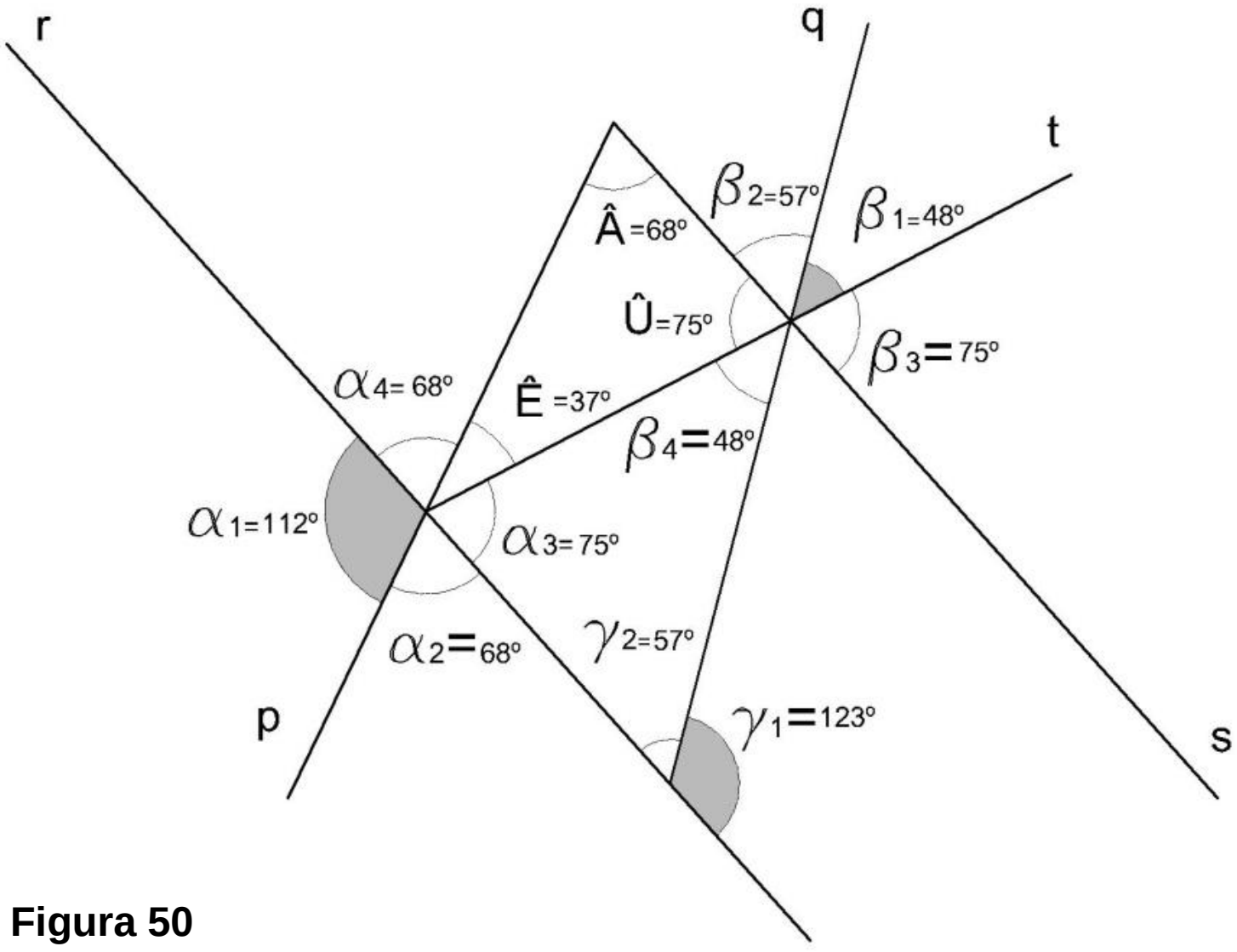

Figura 50

Por ser ángulos suplementarios:
$\alpha_1 = 112°\ ;\ \alpha_1 + \alpha_2 = 180°$
$\Longrightarrow = 68°$
$\gamma_1 = 123°\ ;\ \gamma_1 + \gamma_2 = 180°$
$\Longrightarrow \gamma_2 = 57°$

Por ser ángulos correspondientes:
$\gamma_2 = \beta_2 = 57°$
$\alpha_2 = \hat{A} = \mathbf{68}°$

Por estar a un lado de una línea recta:
$\beta_1 = 48°\ ;\ \beta_2 = 57°;\ \beta_1 + \beta_2 + \beta_3 = 180°$
$\Longrightarrow \beta_3 = 75°$

Por ser ángulos opuestos por el vértice:
$\beta_3 = \hat{U} = \mathbf{75}°$
$\alpha_2 = \alpha_4 = 68°$
$\beta_1 = \beta_4 = 48°$

Por ser ángulos correspondientes:
$\beta_3 = \alpha_3 = 75°$

Por estar a un lado de una línea recta:
$\alpha_3 = 75°\ ;\ \alpha_4 = 68°;\ \alpha_3 + \hat{E} + \alpha_4 = 180°\ \Longrightarrow \hat{E} = \mathbf{37}°$

Comprobación, por ser los ángulos de un triángulo:
$\hat{A} + \hat{E} + \hat{U} = 68 + 37 + 75 = 180°$

<u>Problema 3.1.4.</u>

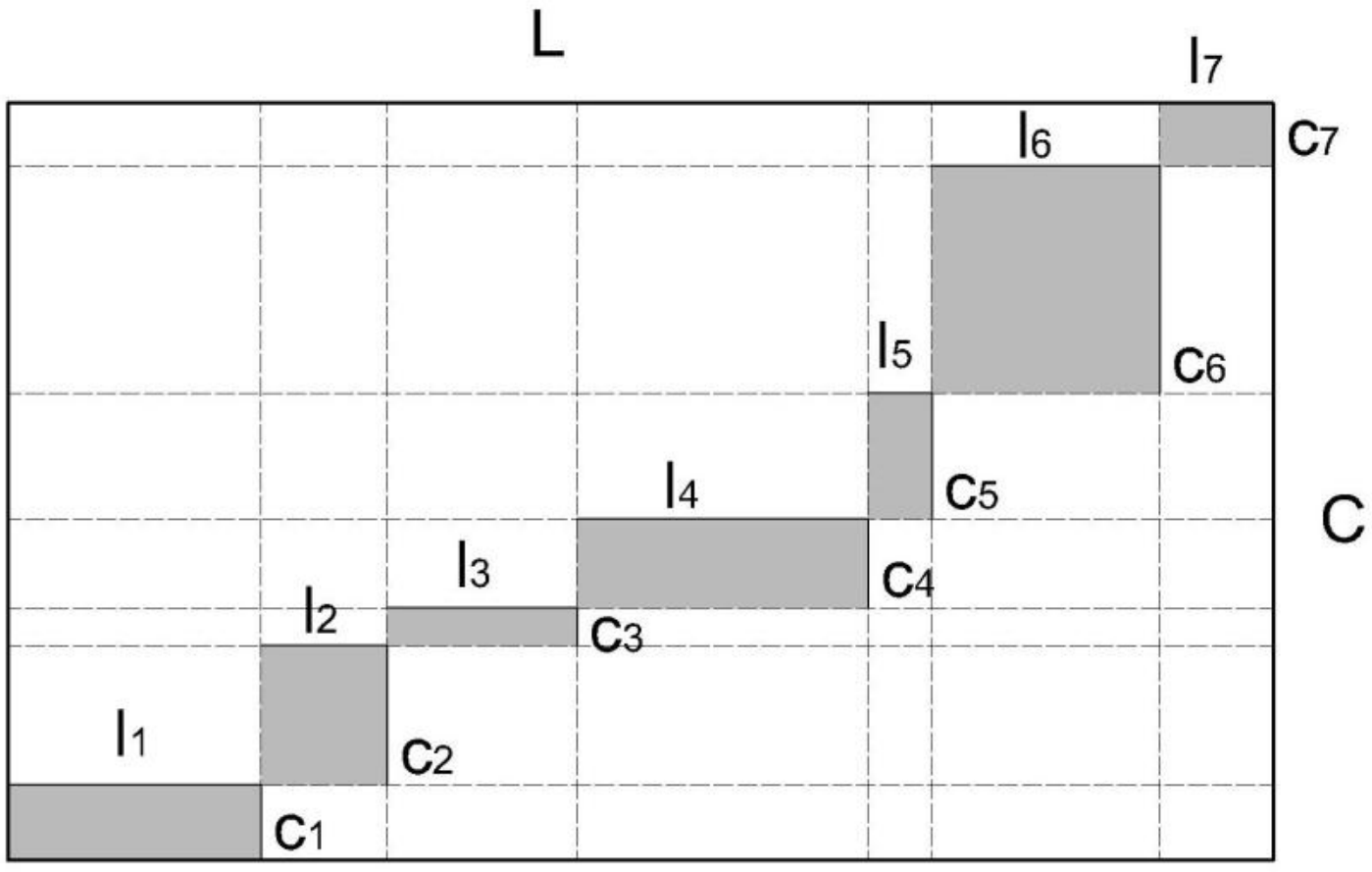

Figura 51

P, el perímetro del rectángulo mayor, se puede expresar como:

$P = 2(L + C) = 2(l_1 + l_2 + l_3 + l_4 + l_5 + l_6 + l_7 + c_1 + c_2 + c_3 + c_4 + c_5 + c_6 + c_7) =$
$2(l_1 + c_1) + 2(l_2 + c_2) + 2(l_3 + c_3) + 2(l_4 + c_4) + 2(l_5 + c_5) + 2(l_6 + c_6) + 2(l_7 + c_7) =$

$= p_1 + p_2 + p_3 + p_4 + p_5 + p_6 + p_7$

Dado que el perímetro de cada uno de los rectángulos pequeños tiene que ser un número entero, entonces, P, resulta ser necesariamente la suma de 7 números enteros, luego P es un número entero.

**Problemas de reserva principales**

Problema 3.1.5.

Divide un cuadrado en cuadrados de 2 tamaños diferentes, de tal forma, que el número de de cuadrados grandes y pequeños sea el mismo (Yashchenko, 2010, p. 10).

Problema 3.1.6.

Prueba en un hexágono regular, cuyos vértices están denominados de forma consecutiva de 1 a 6, que uniendo mediante segmentos los vértices pares (o los impares) se forma un triángulo equilátero (Lehoczky y Rusczyk, 1993, p. 143).

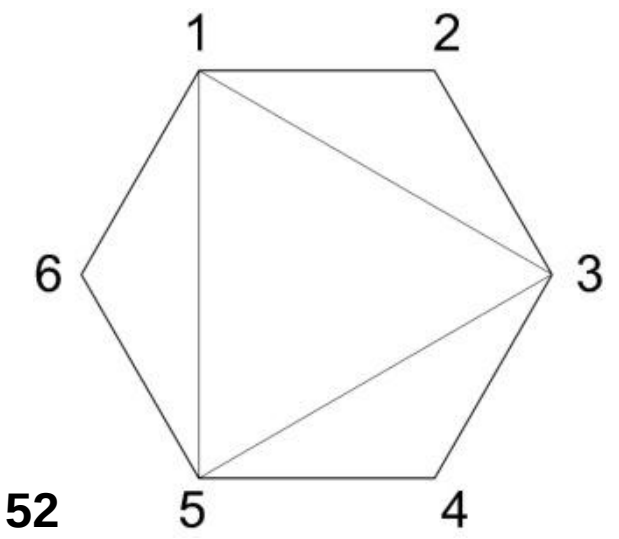

**Figura 52**

Problema 3.1.7.

Sin utilizar el transportador indica el valor del ángulo Â, teniendo en cuenta, que las rectas *r* y *s* son paralelas.

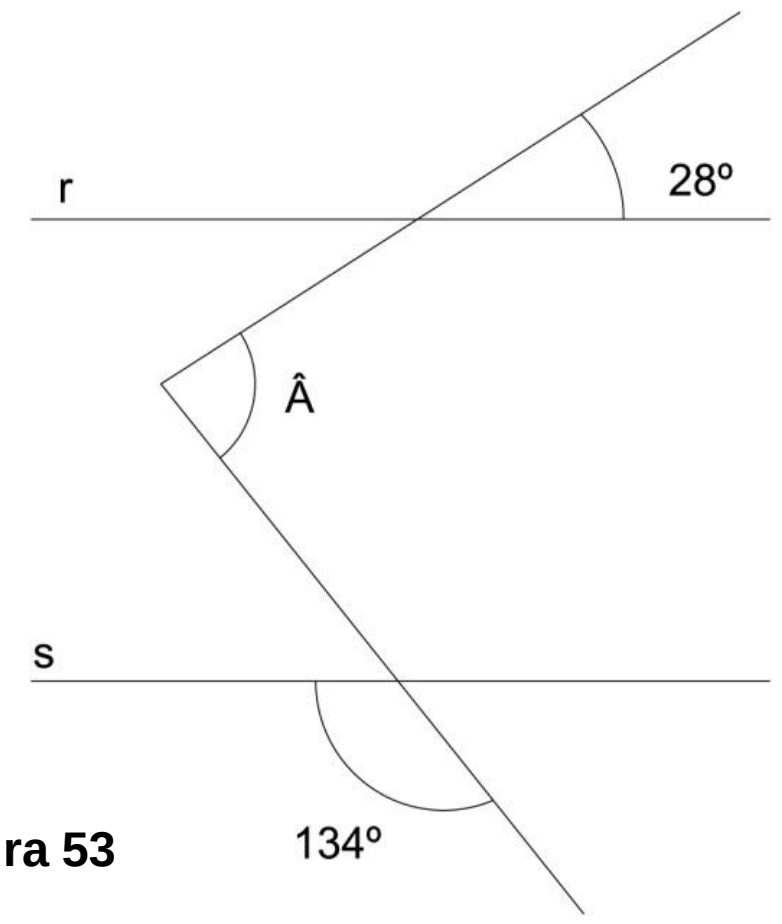

**Figura 53**

**Problemas de reserva adicional:**

Problema 3.1.8.

11 ruedas dentadas situadas en el mismo plano, forman una cadena. ¿Pueden todas ellas rotar simultáneamente? (Fomin, Genkin, y Itenberg, 1996, p. 5).

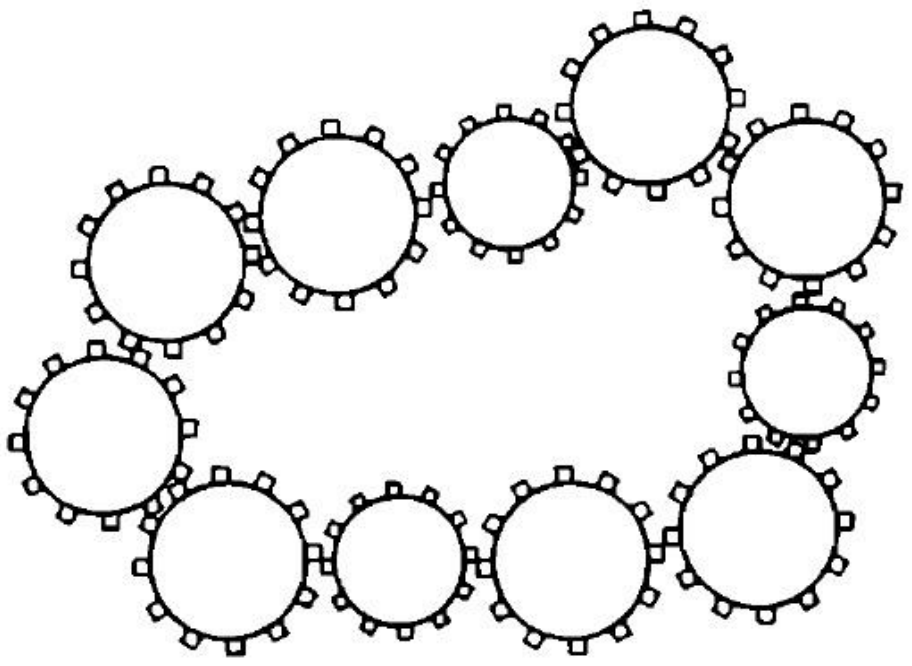

**Figura 54**

**Soluciones**

Problema 3.1.5.

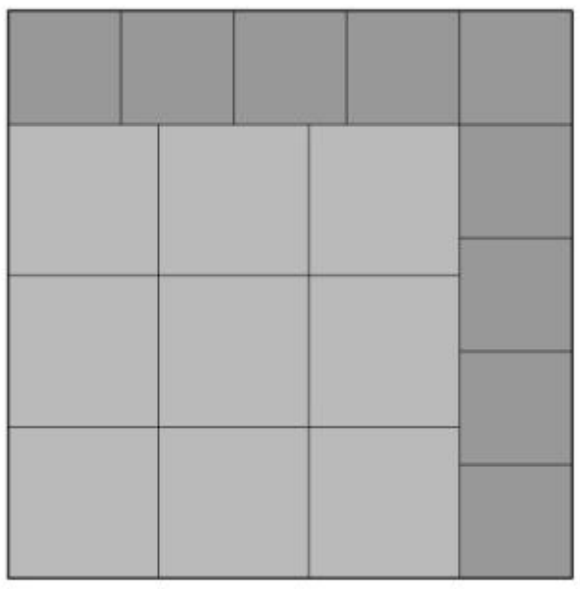

Dividiendo el lado del cuadrado en un número impar $n \in \mathbb{N}$ tal que $2n - 1 = r^2 \; con \; r \geq 3 \in \mathbb{N} \; e \; impar$, y donde r es el número de divisiones del lado de un cuadrado interior al dado, en la misma posición y construido sobre uno de sus vértices, tal y como muestra la figura, resulta sencillo obtener el mismo número de cuadrados, para dos cuadrados de tamaños diferentes.

En la Fig. 55, se ha tomado $r = 3 \;\Longrightarrow 2n - 1 = 9 \Longrightarrow n = 5$

**Figura 55**

Problema 3.1.6.

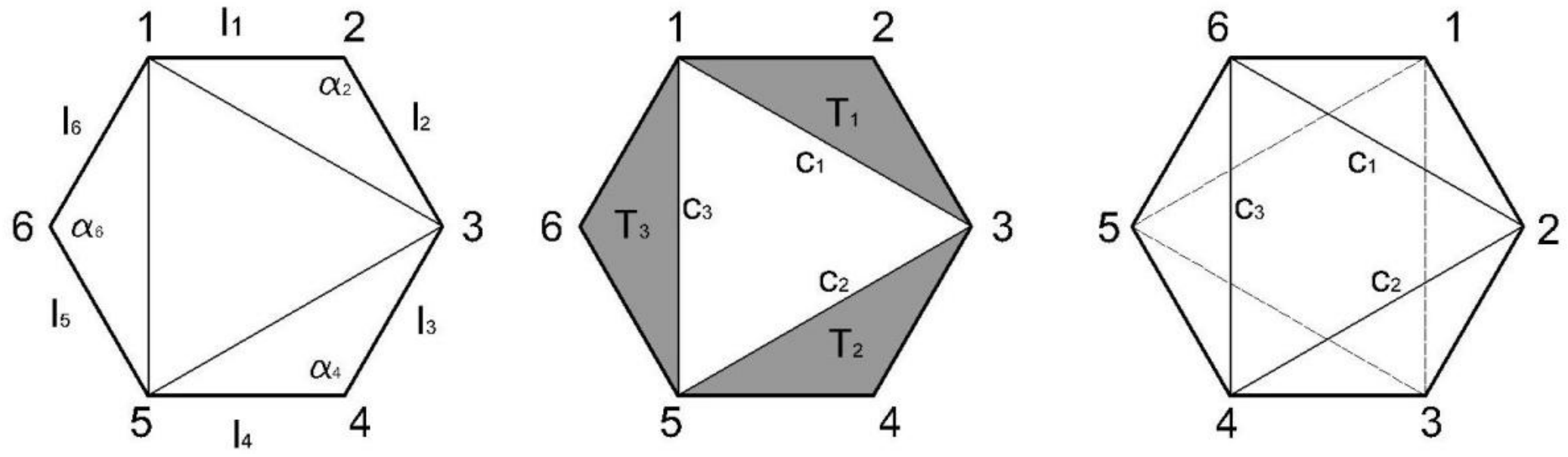

**Figura 56**

Los lados de un polígono regular son iguales. Lo mismo sucede con sus ángulos interiores. Luego: $l_1 = l_2 = l_3 = l_4 = l_5 = l_6$ y $\alpha_2 = \alpha_4 = \alpha_6$

Por ello, necesariamente: $\mathrm{T}_1 = \mathrm{T}_2 = \mathrm{T}_3 \Longrightarrow \mathrm{c}_1 = \mathrm{c}_2 = \mathrm{c}_3$ lo que necesariamente implica que el triángulo sea equilátero

La otra posible situación (polígono trazado sobre los vértices pares), se puede transformar en la ya demostrada, girando la figura el ángulo central del polígono regular.

Problema 3.1.7.

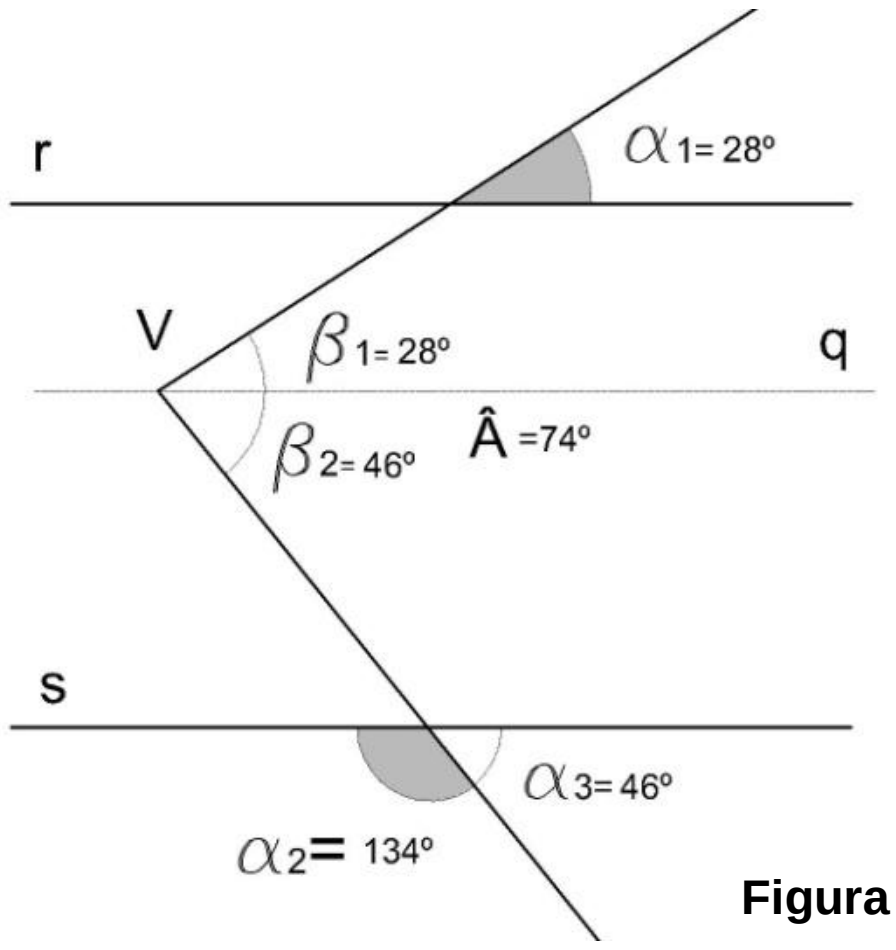

**Figura 57**

Por ser ángulos suplementarios:
$\alpha_2 = 134° \;;\; \alpha_2 + \alpha_3 = 180°$
$\Longrightarrow \alpha_3 = 46°$

Si se traza una recta paralela q a las rectas r y s, por el vértice V, entonces:

Por ser ángulos correspondientes:
$\alpha_1 = \beta_1 = 28°$
$\alpha_3 = \beta_2 = 46°$

Por suma de ángulos adyacentes:

$\hat{\mathrm{A}} = \beta_1 + \beta_2 = 28 + 46 = 74°$

Problema 3.1.8.

Hay 11 ruedas dentadas. Las numeramos arbitrariamente. Suponiendo que todas las ruedas giren, entonces si las ruedas impares giran hacia la derecha, las pares deben hacerlo hacia la izquierda (o viceversa). La rueda 1 y la rueda 11, giran en el mismo sentido, pero como forman una cadena están en contacto y deberían girar una en un sentido y otra en el otro. Luego obtenemos una contradicción y por tanto, no es posible que todas ellas giren simultáneamente tal y como habíamos supuesto.

## A.3.2. Sesión de problemas Tipo SP2.

### Introducción

Esta sesión de problemas se propone para el grupo de 3º y 4º de ESO, con libre elección de tema para cada uno de los problemas.

Se han propuestos temas diversos para los problemas de forma que permitan desarrollar el pensamiento matemático en toda su amplitud y alcanzar varios objetivos, a la vez que el estudiante se entrena para resolver cualquier tipo de problema matemático, sin esperar que se le presente alguno de unas características determinadas.

Se proponen 3 problemas y se dispondrá de otros 3 problemas de reserva, todos ellos de diferente dificultad y que puedan ser utilizados de diversos modos:

- para sustituir o facilitar el camino hacia la resolución de alguno de los propuestos durante la sesión
- para quienes hayan terminado con los problemas propuestos
- como trabajo adicional para casa para quién los solicite

Los problemas no se recogerán, sino que se podrá seguir trabajando en ellos hasta que sean comentados finalmente en una sesión posterior. En la presente sesión se indicará que problemas se tiene previsto que se comenten en la próxima.

### Objetivos

| | | |
|---|---|---|
| Problema 3.2.1. | Polígonos, bisectriz | Empleo de elementos auxiliares: ángulo central e inscrito en una circunferencia |
| Problema 3.2.2. | Aritmética | Tanteo |
| Problema 3.2.3. | Patrones, recurrencia | Convertir un problema en otro, iteración |
| Problema 3.2.4. | Triángulos | Visión espacial |
| Problema 3.2.5. | Triángulos, bisectriz | Empleo de elementos auxiliares: igualdad de ángulos, semejanza de triángulos |
| Problema 3.2.6. | Circunferencias, tangencias | Convertir un problema en otro, homotecia, recurrencia |

**Problemas**

Problema 3.2.1.

Sea P el centro de un cuadrado construido sobre la hipotenusa de un triángulo rectángulo ABC, prueba que BP es la bisectriz del ángulo recto del triángulo (Schoenfeld, 1985, p. 267).

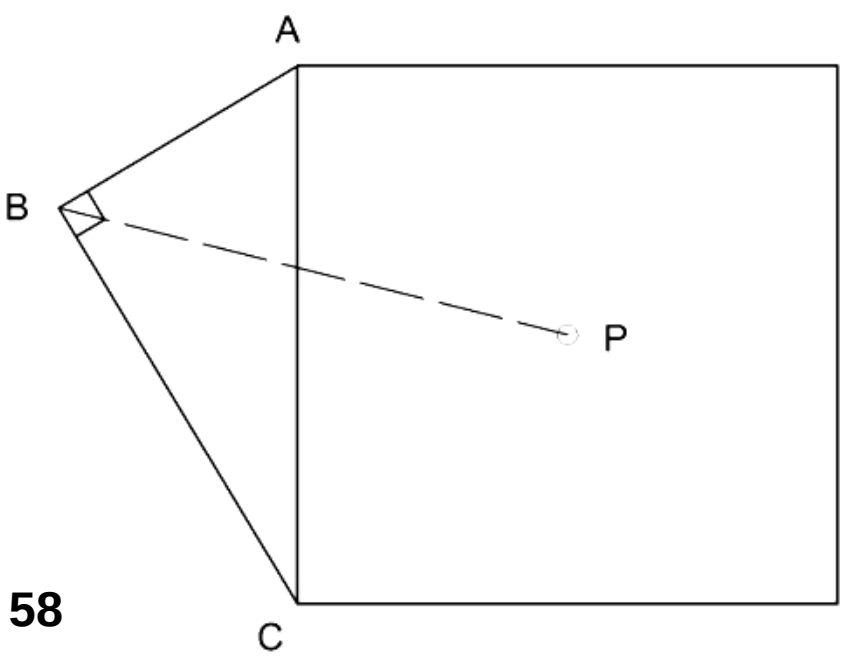

Figura 58

Problema 3.2.2.

La fuente y las jarras. Nos encontramos frente a una fuente y disponemos de dos jarras de 3 y 5 litros.¿Cómo podemos conseguir exactamente cuatro litros de agua sin utilizar ningún otro recipiente? (Epsilones: La fuente y las jarras, 2002).

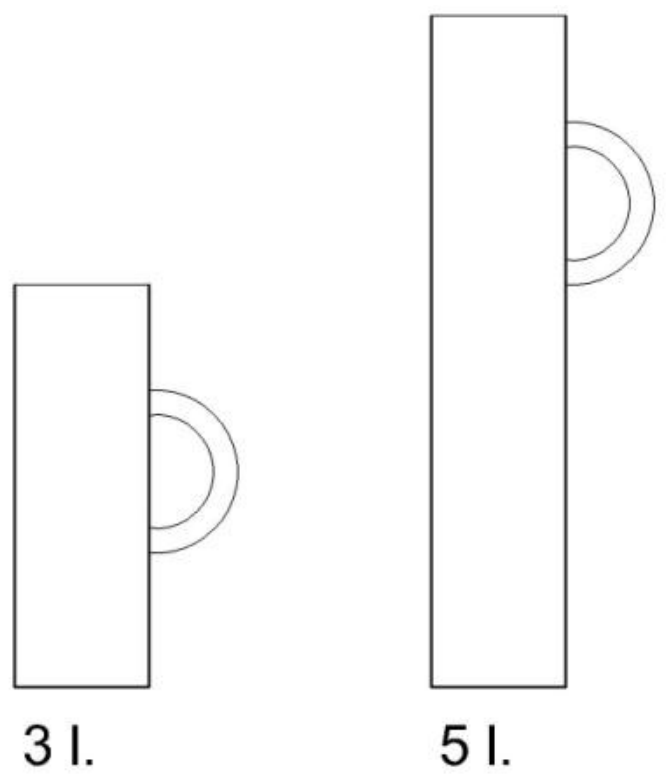

Figura 59

Problema 3.2.3.

Dada la siguiente sucesión de figuras que se muestra en la Fig. 60, y Sea $c_n$ el número de palillos para la construcción de la n-ésima figura. Se pide:

- Averiguar si $c_n$ está definido por recurrencia
- Determinar expresión explícita para $c_n$

(Madrid, 2015)

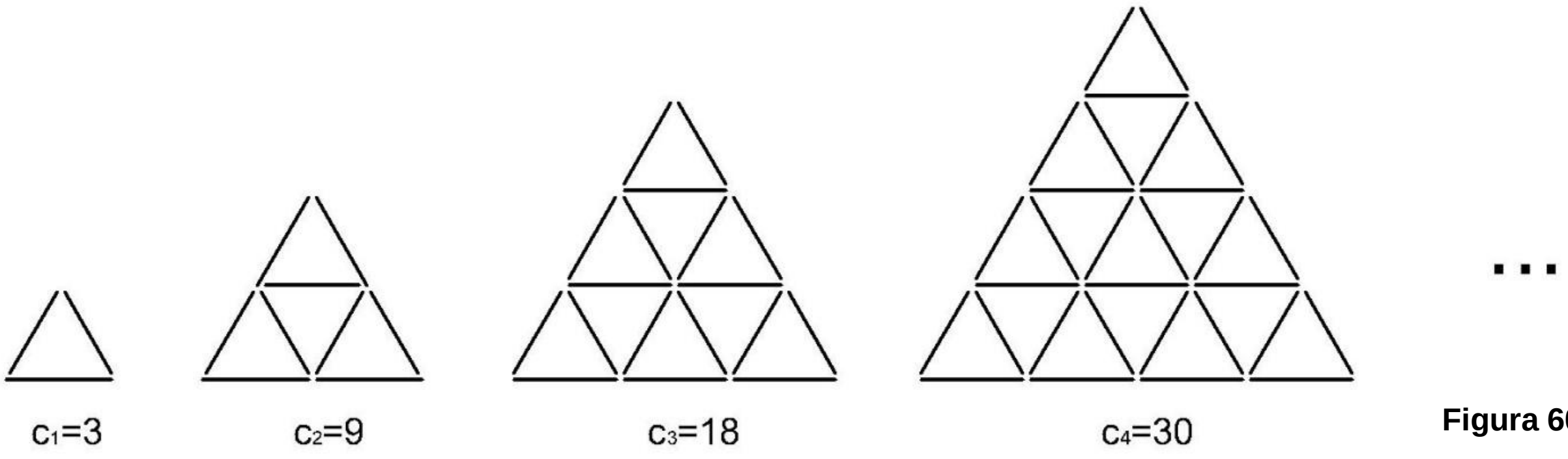

Figura 60

## Soluciones

Problema 3.2.1.

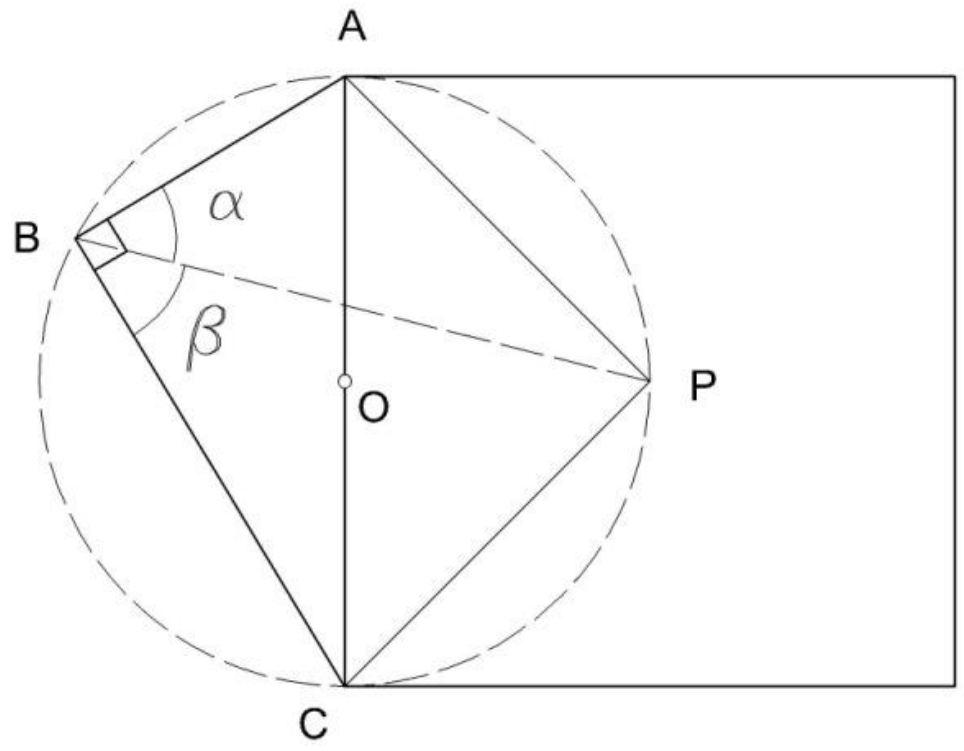

**Figura 61**

Trazando un circunferencia de centro O, siendo el punto medio del lado AC y radio OA, la circunferencia pasa necesariamente por B pues el ángulo ABC (α+β) es un ángulo recto y abarca al diámetro AC de la circunferencia.

Considerando que el triángulo rectángulo es exterior al cuadrado, el vértice B puede ocupar únicamente los puntos sobre la semicircunferencia externa al mismo.

Por tratarse de un cuadrado AP=CP. De los posibles lugares que puede ocupar B, el ángulo α siempre abarcará el segmento AP y el ángulo β siempre abarcará el segmento CP, luego α=β, y como α+ β=90º entonces: α= β=45º

Problema 3.2.2.

Llamaremos J3 a la jarra de 3 litros y J5 a la jarra de 5 litros. Hay varias posibilidades.
Una posibilidad:

- $(\frac{0}{3},\frac{5}{5})$: Llenamos J5.
- $(\frac{3}{3},\frac{2}{5})$: Pasamos 3 litros de J5 a J3. Luego en J5 hay 2 litros.
- $(\frac{0}{3},\frac{2}{5})$: Vaciamos J3.
- $(\frac{2}{3},\frac{0}{5})$: Pasamos 2 litros de J5 a J3. Luego J5 está vacía.
- $(\frac{2}{3},\frac{5}{5})$: Llenamos J5.
- $(\frac{3}{3},\frac{4}{5})$: Pasamos 1 litro de J5 a J3. Luego en J5 hay 4 litros.
- Número de pasos: 6. Desperdicio de agua: 6 litros.

Otra posibilidad:

- $(\frac{3}{3},\frac{0}{5})$: Llenamos J3.
- $(\frac{0}{3},\frac{3}{5})$: Pasamos los 3 litros de J3 a J5.
- $(\frac{3}{3},\frac{3}{5})$: Llenamos J3.
- $(\frac{1}{3},\frac{5}{5})$: Pasamos 2 litros de J3 a J5. Luego J5 está llena y en J3 hay 1 litro.
- $(\frac{1}{3},\frac{0}{5})$: Vaciamos J5
- $(\frac{0}{3},\frac{1}{5})$: Pasamos 1 litro de J3 a J5. Luego J3 está vacía.
- $(\frac{3}{3},\frac{1}{5})$: Llenamos J3
- $(\frac{0}{3},\frac{4}{5})$: Pasamos 3 litros de J3 a J5 . Luego en J5 hay 4 litros y J3 está vacía.
- Número de pasos: 7. Desperdicio de agua: 5 litros.

(Epsilones: La fuente y las jarras, 2002)

Problema 3.2.3.

Si en vez de contar cerillas, contamos únicamente los triángulos necesarios para disponer de ellas, habremos transformado el problema en otro más sencillo, de la resolución del cual resulta evidente la respuesta al problema inicial.

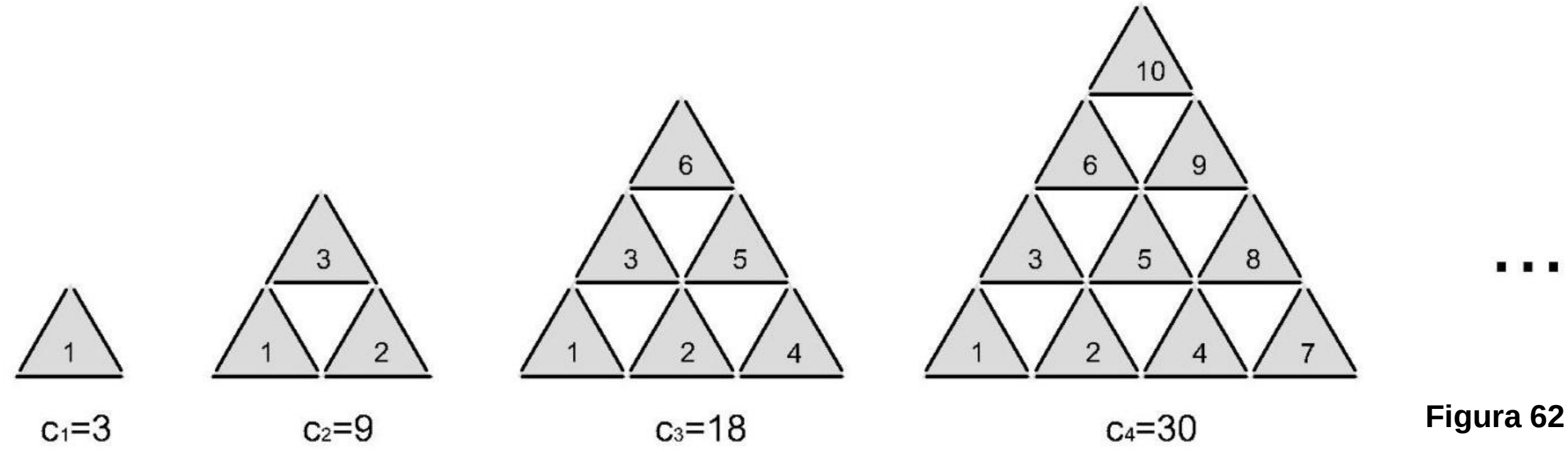

**Figura 62**

La sucesión recurrente queda definida por: $\begin{cases} c_1 = 3 \\ c_n = c_{n-1} + 3n \end{cases}$

Mediante iteración resulta sencillo llegar a la formulación de la suma de los n términos de una progresión aritmética:

$$c_n = 3 + 3 \cdot 2 + 3 \cdot 3 + \cdots + 3 \cdot (n-1) + 3 \cdot n = \sum_{i=1}^{n} 3 \cdot i$$

> Se llama progresión aritmética una sucesión de números, en la cual cada término se obtiene del anterior sumando a este un mismo número denominado diferencia de la progresión. […]
>
> Si designamos la suma de los n primeros términos de la progresión aritmética como:
>
> $S_n = a_1 + a_2 + \cdots + a_{n-1} + a_n$
>
> Si los sumandos del segundo miembro de la igualdad los escribimos en el orden inverso, la suma $S_n$ no variará por ello:
>
> $S_n = a_n + a_{n-1} + \cdots + a_2 + a_1$
>
> Sumando miembro a miembro ambas igualdades, obtendremos:
>
> $2S_n = (a_1 + a_n) + (a_2 + a_{n-1}) + \cdots + (a_{n-1} + a_2) + (a_n + a_1)$
>
> En cada paréntesis tenemos una suma de dos términos equidistantes de la progresión; por lo tanto todas estas sumas entre paréntesis son iguales entre sí. […] En total son n parétesis, es decir, tántos como términos de la progresión. Por eso tendremos:
>
> $2S_n = n(a_1 + a_n)$
>
> (Kalnin, 1978, pp. 280-283)

Luego, podemos expresar la suma de los n términos de una progresión aritmética como:

$$\sum_{i=1}^{n} a_i = \frac{n(a_1 + a_n)}{2}; \textit{donde } a_1 \textit{ es el primer término y } a_n \textit{ el último}$$

De esta forma resulta, hallar una expresión explícita para $c_n = \sum_{i=1}^{n} 3 \cdot i$

$$c_n = \frac{n(3+3n)}{2} = \frac{3}{2} n(n+1)$$

**Problemas de reserva**

Problema 3.2.4.

Demostrar el **Teorema del ángulo bisector de un triángulo**: si $\overline{AX}$ es la bisectriz de $\widehat{CAB}$ del triángulo $ABC$, entonces $\frac{\overline{AC}}{\overline{CX}} = \frac{\overline{AB}}{\overline{BX}}$. (Lehoczky y Rusczyk, 1993, p. 103).

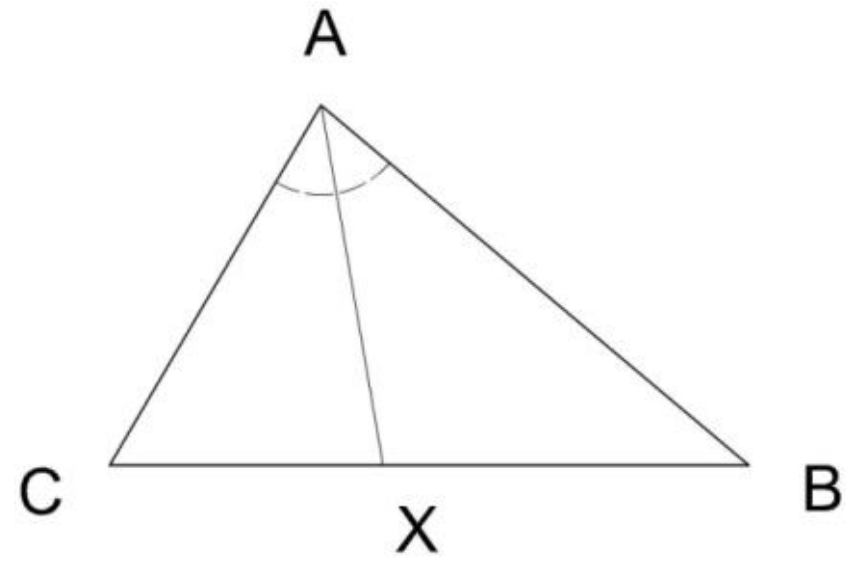

**Figura 63**

Problema 3.2.5.

Dado un triángulo formado con tres cerillas, ¿cómo conseguir seis triángulos añadiendo solo tres cerillas y sin mover las anteriores? (Epsilones: Triángulos y cerillas (II), 2002). ¿Cómo conseguir ocho triángulos añadiendo solo tres cerillas más y sin mover las anteriores? (Epsilones: Triángulos y cerillas (III), 2002).

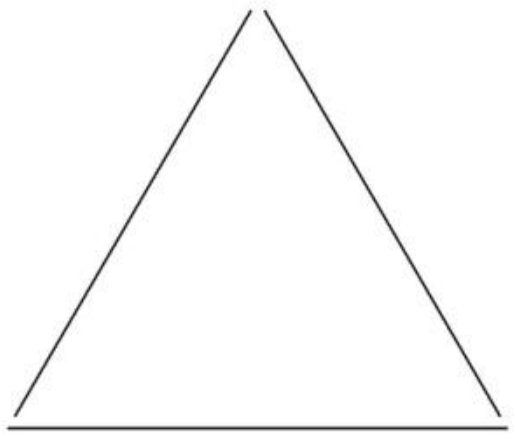

**Figura 64**

Problema 3.2.6.

En la figura adjunta, las cinco circunferencias son tangentes de forma consecutiva entre ellas, y todas ellas son tangentes a las dos rectas. Si el radio de la circunferencia menor es 3 y el de la mayor 48, ¿cuál es el radio de la circunferencia que se encuentra situada en el medio? (Maurer y Berzsenyi, 1997, p. 5).

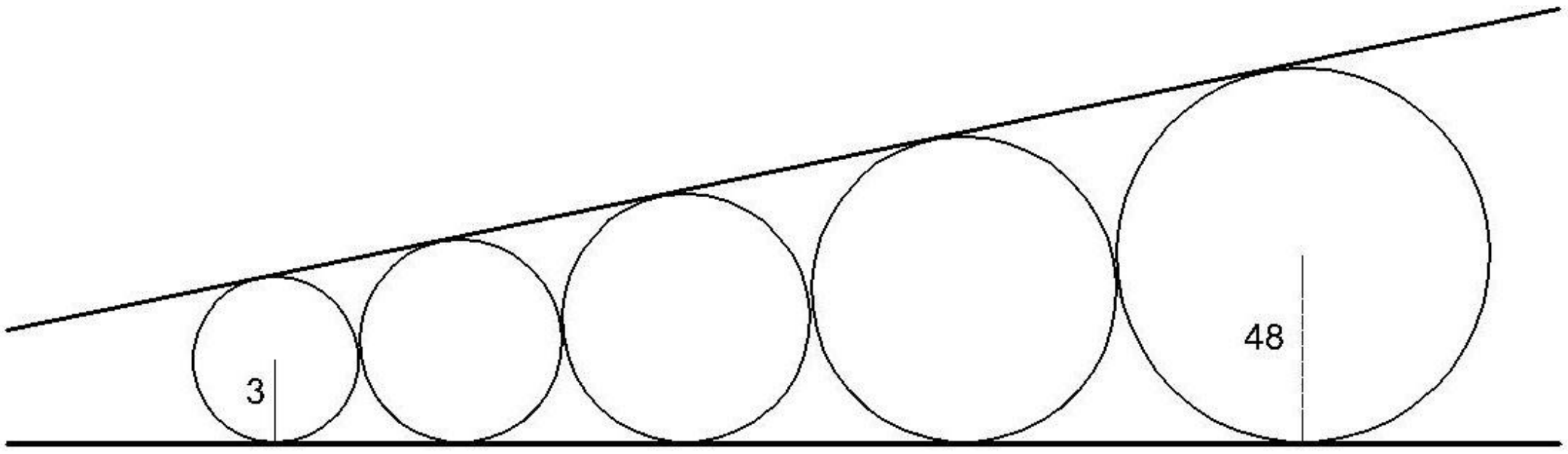

**Figura 65**

## Soluciones

Problema 3.2.4.

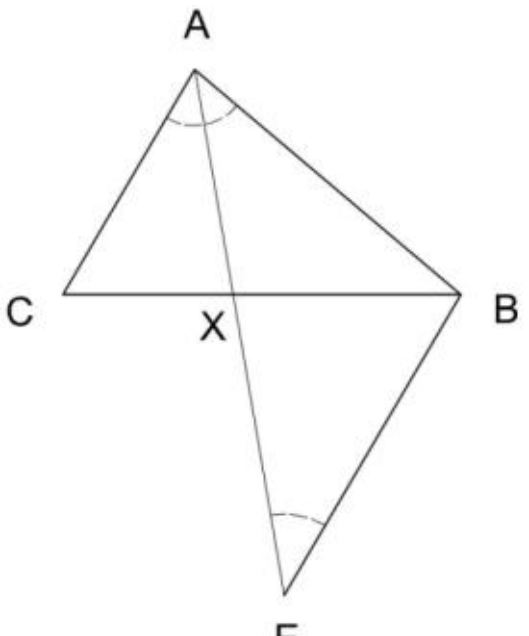

Figura 66

Trazamos una recta paralela a $\overline{AC}$ por $B$. Prolongamos $\overline{AX}$ hasta obtener el punto de intersección E con dicha recta.

$\widehat{CAE} = \widehat{AEB}$ por ser ángulos alternos internos y $\widehat{CAX} = \widehat{XAB}$ ya que $\overline{AX}$ es la bisectriz de $\widehat{CAB}$ luego $\widehat{EAB} = \widehat{AEB}$

Como $\widehat{EAB} = \widehat{AEB}$ entonces el triángulo ABE es isósceles y necesariamente $\overline{AB} = \overline{BE}$

Como $\widehat{CAX} = \widehat{XEB}$ por ser ángulo alternos internos y $\widehat{AXC} = \widehat{BXE}$ por ser ángulos opuestos por el vértice, entonces el triángulo BXE es semejante con el triángulo AXC, y por ello: $\frac{\overline{AC}}{\overline{CX}} = \frac{\overline{BE}}{\overline{BX}} = \frac{\overline{AB}}{\overline{BX}}$

Problema 3.2.5.

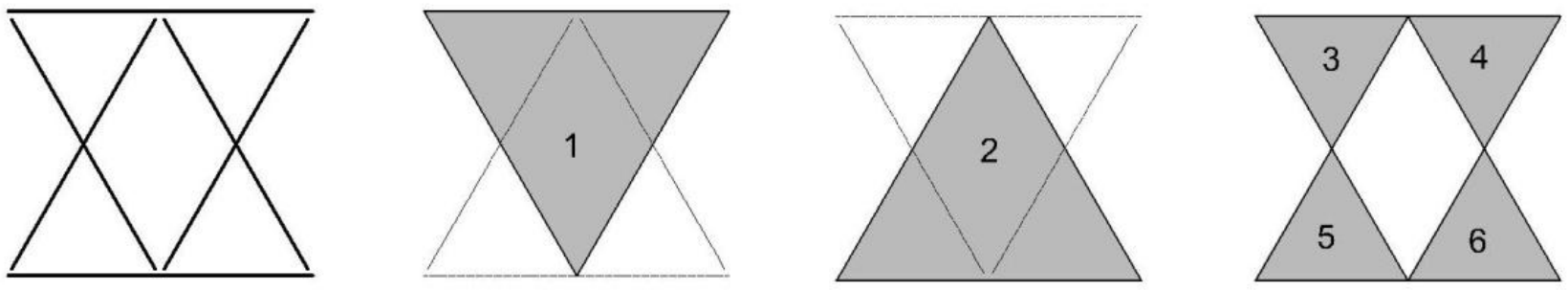

Figura 67

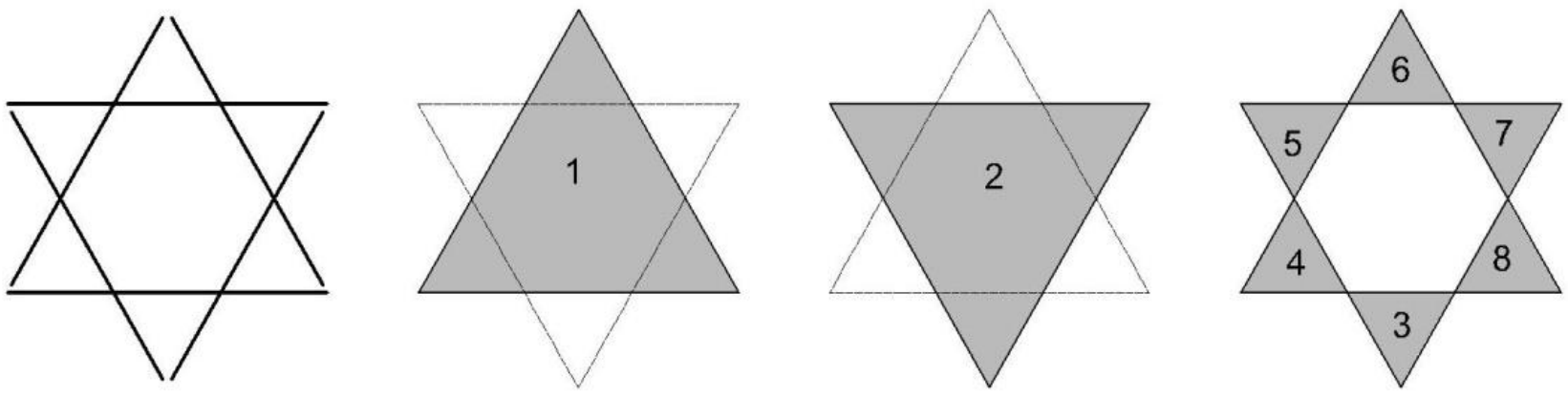

Figura 68

Problema 3.2.6.

Si dilatamos la circunferencia $C_1$ hasta convertirla en $C_2$, podemos pensar en una homotecia de centro O y razón k, donde la longitud el radio de $C_2$ vendrá dada por: $r_2 = k \cdot r_1$

Una vez que hemos obtenido $C_1$ y $C_2$, si volvemos a dilatarlas conjuntamente por la misma homotecia, entonces la longitud de sus radios vendrá dada por:

$r_2 = k \cdot r_1$ y $r_3 = k \cdot r_2 r_1$ luego $r_3 = k \cdot k \cdot r_1 = k^2 \cdot r_1$

Siguiendo ese mismo proceso de manera recursiva, podemos obtener fácilmente una expresión explícita para el valor del radio de la circunferencia n-ésima: $r_n = k^{n-1} r_1$ que constituye una progresión geométrica de razón k.

Luego para n=5: $r_5 = k^4 r_1 \Leftrightarrow k^4 = \frac{r_5}{r_1} = \frac{48}{3} = 16 \Leftrightarrow k = \sqrt[4]{16} = 2$

Conocida la razón $k = 2$ es sencillo obtener el radio de cualquiera de las circunferencias de la progresión. Para la circunferencia del medio, n=3: $r_3 = k^2 r_1 = 2^2 \cdot 3 = 4 \cdot 3 = 12$

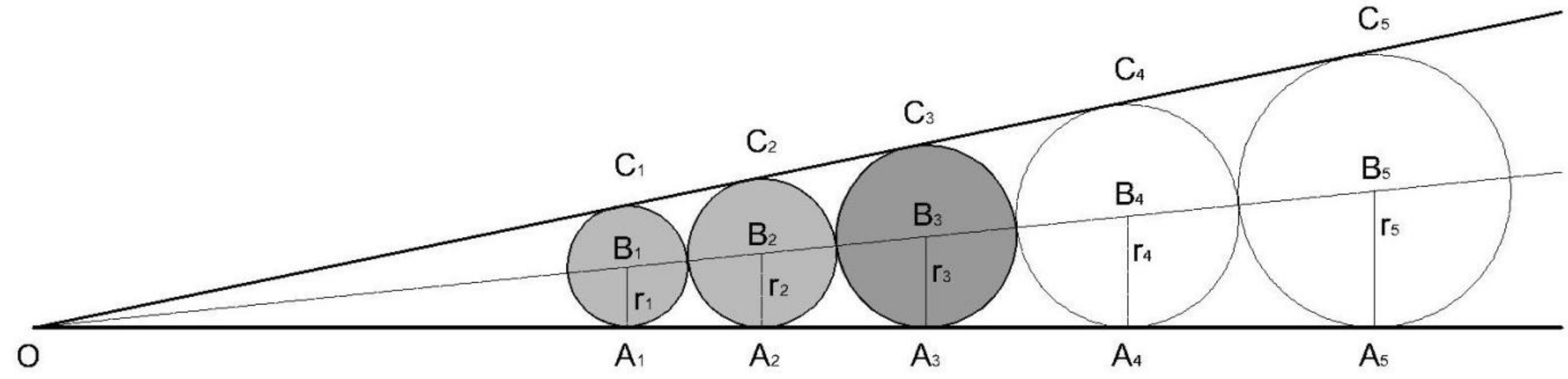

Figura 69

## Relación de figuras utilizadas

## Bibliografía

## Otro material consultado